\documentclass[10pt]{amsart}
\usepackage{amsmath}
\usepackage{amssymb}
\usepackage{amsfonts}
\usepackage{amsthm}
\usepackage{enumerate}
\usepackage[all]{xy}
\usepackage[latin1]{inputenc}
\usepackage{graphicx}
\usepackage{latexsym}
\usepackage{mathdots}
\usepackage{mathrsfs}
\input xy

\makeatletter

\newtheoremstyle{definition}
        {5pt}
        {3pt}
        {}
        {0pt}
        {\scshape}
        {.}
        {5pt}
        {\thmname{#1} \thmnumber{#2} \thmnote{[#3]}} 

\newtheoremstyle{theorems}
        {5pt}
        {3pt}
        {\itshape}
        {0pt}
        {\scshape}
        {.}
        {5pt}
        { \thmname{#1} \thmnumber{#2}\thmnote{[#3]}} 

\swapnumbers

\renewcommand\section{\@startsection{section}{1}{\z@}%
        {-3.5ex \@plus -1ex \@minus -.2ex}%
        {2.3ex \@plus .2ex}%
        {\centering\reset@font\scshape}}

\makeatother

\theoremstyle{theorems}
\newtheorem{Theo}{Theorem}[section]
\newtheorem{Prop}[Theo]{Proposition}
\newtheorem{Cor}[Theo]{Corollary}
\newtheorem{Lemma}[Theo]{Lemma}

\theoremstyle{definition}
\newtheorem{Defn}[Theo]{Definition}
\newtheorem{Notn}[Theo]{Notation}

\newcommand{\rep}{{\rm rep}}
\newcommand{\Rep}{{\rm Rep}}
\newcommand{\Hom}{{\rm Hom}}
\newcommand{\Ext}{{\rm Ext}}
\newcommand{\End}{{\rm End}}

\newcommand{\N}{\mathbb{N}}
\newcommand{\Z}{\mathbb{Z}}

\newcommand{\A}{\mathbb{A}}
\newcommand{\D}{\mathbb{D}}
\newcommand{\DTr}{{\rm D\hskip -1pt Tr}\hskip 0.6 pt}
\newcommand{\TrD}{{\rm Tr \hskip -0.3 pt D}\hskip 0.6 pt}

\def\id{\hbox{1\hskip -3pt {\sf I}}}
\def\Ta{\hbox{$\mathit\Theta$}}
\def\Ga{\hbox{$\mathit\Gamma$}}
\def\Da{\hbox{${\mathit\Delta}$}}
\def\Sa{\hbox{${\mathit\Sigma}$}}
\def\Oa{\hbox{${\mathit\Omega}$}}

\def\taud{\hbox{$\tau_{_{\hspace{-0.4pt} D}}$}}

\def\H{\hbox{$\mathcal{H}$}}

\begin{document}

\title{\sc Representation Theory of an Infinite Quiver}

\author{\tiny Raymundo Bautista, Shiping Liu, and Charles Paquette}

\address{Raymundo Bautista, Centro de Ciencias Matematicas, UNAM,
Apartado Postal 61-3, 58089 Morelia, Mexico}
\email{raymundo@matmor.unam.mx}
\address{Shiping Liu, D\'epartement de math\'ematiques, Universit\'e de Sherbrooke, Sherbrooke, Qu\'ebec, Canada, J1K 2R1}
\email{shiping.liu@usherbrooke.ca}
\address{Charles Paquette, D\'epartement de math\'ematiques, Universit\'e de Sherbrooke, Sherbrooke, Qu\'ebec, Canada, J1K 2R1}
\email{charles.paquette@usherbrooke.ca}

\maketitle

\begin{abstract}

\vspace{-25pt}

This paper deals with the representation theory of a locally
finite quiver in which the number of paths between any two given
vertices is finite. We first study some properties of the finitely
presented or co-presented representations, and then construct in
the category of locally finite dimensional representations some
almost split sequences which start with a finitely co-presented
representation and end with a finitely presented representation.
Furthermore, we obtain a general description of the shapes of the
Auslander-Reiten components of the category of finitely presented
representations and prove that the number of regular
Auslander-Reiten components is infinite if and only if the quiver
is not of finite or infinite Dynkin type. In the infinite Dynkin
case, we shall give a complete list of the indecomposable
representations and an explicit description of the
Auslander-Reiten components. Finally, we apply these results to
study the Auslander-Reiten theory in the derived category of
bounded complexes of finitely presented representations.

\end{abstract}

\bigskip
\smallskip

\centerline{\sc Introduction}

\medskip

As the best understood and the most stimulating part of the
representation theo\-ry of finite dimensional algebras, the theory
of representations of a finite quiver without oriented cycles has
been extensively studied over the last four decades; see, for
example, \cite{BGP, Gab, GR, Kac, Ker, R1, R3}. On the other hand,
the Auslander-Reiten theory of irreducible morphisms and almost
split sequences provides an indispensable powerful tool for the
representation theory and it appears in many other areas such as
algebraic geometry and algebraic topology; see \cite{AuR, Aus2,
Jo}. The impact of these two theories to other branches of
mathematics is best illustrated by their recent interaction with
the theory of cluster algebras via the cluster category; see, for
example, \cite{BMRRT, FoZ}. Now, new developments require the
study of representations of an infinite quiver. For instance, in
order to classify the noetherian Ext-finite heredi\-tary abelian
categories with Serre duality, Rei\-ten and Van den Bergh
investigated the category of finitely presented representations of
a locally finite quiver without left infinite paths; see
\cite{RVDB}. In particu\-lar, they showed that this category has
right almost split sequences and obtained some partial description
of its Auslander-Reiten components. Later, by considering
representations of ray quivers which are infinite in general,
Ringel provided an alternative construction of the noetherian
Ext-finite hereditary abelian categories which have Serre
duality and non-zero projective objects; see \cite{R2}. More
recently, Holm and J{\o}rgensen studied a cluster category of
infinite Dynkin type, which can be constructed from the category
of finite dimensional representations of a quiver of type
$\A_\infty$ and whose Auslander-Reiten quiver is of shape
$\Z\A_\infty$; see \cite{HoJ}. Finally, the bounded derived
category of a finite dimensional algebra with radical squared zero
is determined by the category of finitely co-presented
representations of a covering of the ordinary quiver of the
algebra which is usually infinite; see \cite{BaL}, and homogeneous
vector bundles over certain algebraic variety are characterized by
the finite dimensional representations of some infinite quiver;
see \cite{Hille}.

\medskip

All these motivate us to study representations of a quiver in the
most general setting. Indeed, we shall work with an arbitrary base
field $k$ and a quiver $Q$ which is assumed only to be locally
finite such that the number of paths between any given pair of
vertices is finite. The main objective of this paper is to present a
complete picture of the Auslander-Reiten theory in the category of
finitely presented $k$-representations of $Q$ and in its bounded
derived category. Our results yield many interesting examples of
Ext-finite, but not necessarily noetherian, hereditary abelian
categories which have (left, right) almost split sequences. We
outline the content of the paper section by section as follows.

\medskip

In Section $1$, we study almost finitely presented and almost
finitely co-presented representations; see
(\ref{defn1.5_shiping}), which are slightly more general than
finitely presented and finitely co-presented representations,
respectively. We shall give some combinatorial characterizations of
these representations; see (\ref{theo1.10_shiping}). These will
be useful for us to relate almost finitely presented or
co-presented representations of $Q$ to
representations of its finite subquivers; see
(\ref{restriction}).

\medskip

In Section $2$, we study the Auslander-Reiten theory in the
category $\rep(Q)$ of locally finite dimensional
$k$-representations of $Q$. We first construct almost split
sequences which start with a finitely co-presented representation
and end with a finitely presented one; see (\ref{AR-sequence}),
and then study some important properties of the Auslander-Reiten
orbits in $\rep(Q)$; see (\ref{orbits}). Finally, we develop some
tools to relate almost split sequences and irreducible morphisms
in $\rep(Q)$ to those of locally finite dimensional representations of
subquivers of $Q\,;$ see (\ref{ARS-extension}) and
(\ref{prop2.9}).

\medskip

In Section $3$, we start to concentrate on the study of the
Auslander-Reiten theory in the category $\rep^+(Q)$ of finitely
presented $k$-representations of $Q$. We shall show that
irreducible morphisms between indecomposable representations in
$\rep^+(Q)$ are irreducible in $\rep(Q)$; see (\ref{irr+rep}), and
almost split sequences in $\rep^+(Q)$ have finite dimensional
starting term and remain almost split in $\rep(Q)$; see
(\ref{ARrep+}). This will enable us to find necessary and
sufficient conditions so that $\rep^+\hspace{-1pt}(Q)$ has (left,
right) almost split sequences; see (\ref{RARcat}) and
(\ref{AR-category}).

\medskip

In Section $4$, we shall give a general description of the
Auslander-Reiten quiver $\Ga_{\rep^+(Q)}\vspace{1pt}$ of
$\rep^+\hspace{-1pt}(Q)$, which is defined since $\rep^+(Q)$ is
Krull-Schmidt; see \cite[(2.1)]{L3}. In case $Q$ is connected,
$\Ga_{\rep^+(Q)}\vspace{1pt}$ has a unique preprojective component
which is a predecessor-closed subquiver of
$\mathbb{N}Q^{\rm\,op}$; see (\ref{prep-component}). The
preinjective component correspond to the connected components of the subquiver of
$Q$ generated by the vertices which are not ending point of any left infinite paths, and they
are finite or infinite successor-closed subquiver of
$\mathbb{N}^-\hspace{-2pt}Q^{\rm\,op}$; see (\ref{preinjcomp}).
In particular,  the
number of preinjective components varies from zero to infinity.
The other connected components of $\Ga_{\rep^+(Q)}\vspace{0pt}$,
called {\it regular components}, are of shape
$\mathbb{N}^-\hspace{-1pt}\mathbb{A}_{\infty}$,
$\mathbb{N}\mathbb{A}_{\infty}$, $\mathbb{Z}\mathbb{A}_{\infty},$
or finite wings; see (\ref{reg-cpt}). As a consequence,
$\Ga_{\rep^+(Q)}\vspace{1pt}$ is always symmetrically valued. At
the end of this section, some conditions on $Q$ will be given so
that at most one type of regular components will appear; see
(\ref{no-infpath-cpt}) and (\ref{res-reg-cpt}).

\medskip

In Section $5$, we study the representation theory of infinite
Dynkin types $\A_\infty$, $\A_\infty^\infty,$ and $\D_\infty$. As
in the finite Dynkin case, we shall obtain a complete list of the
indecomposable representations in $\rep^+(Q)$; see (\ref{5.1}) and
(\ref{indec-D}). In case $Q$ is of type $\mathbb{A}_{\infty}$ or
$\mathbb{A}_{\infty}^{\infty}$, we shall describe the irreducible
morphisms between indecomposable representations and determine the
almost split sequences with an indecomposable middle term; see
(\ref{irr-mor}) and (\ref{5.3}). This allows us to give an
explicit description of the connected components of
$\Ga_{\rep^+(Q)}$ for each of the three infinite Dynkin types; see
(\ref{A-infinity}), (\ref{A-double-infity}), and
(\ref{D-inf-main}). In summary, $\Ga_{\rep^+(Q)}$ has at most four
connected components of which at most two are regular.

\medskip

In Section $6$, we prove that $\Ga_{\rep^+(Q)}$ has infinitely many
regular components in case $Q$ is not of finite or infinite Dynkin
type; see (\ref{maintheonumber}). Moreover, we will show that each
of the four types of regular components could appear infinitely many
times in certain particular situation; see (\ref{inf-num-cpt}).

\medskip

In Section $7$, we study the Auslander-Reiten theory in the
bounded derived category $D^b(\rep^+\hspace{-1pt}(Q))$ of
$\rep^+(Q)$. Indeed, for an arbitrary hereditary abelian category
$\H$, we shall describe the so-called {\it connecting almost split
triangles} in $D^b(\H)$ and prove that all other almost split
triangles are induced from the almost split sequences in $\H$; see
(\ref{ART}). Combining this with previously obtained results, we
get a general description of the Auslander-Reiten quiver of
$D^b(\rep^+\hspace{-1pt}(Q))$; see (\ref{proparrowsDQ}),
(\ref{conn}) and (\ref{Der-cpt}). We conclude the paper with some
necessary and sufficient conditions on $Q$ so that
$D^b(\rep^+(Q))$ has (left, right) almost split triangles; see
(\ref{fin-theo}).

\vspace{-5pt}

{\center \section{Almost finitely presented representations}}

The objective of this section is to investigate two classes of
representations, called {\it almost finitely presented} and {\it
almost finitely co-presented}, of a quiver over a field. This will
yield several hereditary abelian categories, some of them are
Krull-Schmidt in which we shall be able to study the
Auslander-Reiten theory. Since a direct sum of injective modules
is not necessarily injective, we shall concentrate on almost
finitely co-presented representations while the corresponding
results for almost finitely presented representations will follow
dually.

\medskip

We start by introducing some combinatorial terminology and
notation. Throughout this paper, $Q=(Q_0, Q_1)$ stands for a
quiver, where $Q_0$ is the set of vertices and $Q_1$ is the set of
arrows. Let $\alpha: x\to y$ be an arrow in $Q$. We call $x$ the
{\it starting point} and $y$ the {\it ending point} of $\alpha$,
and write $s(\alpha)=x$ and $e(\alpha)=y$. One introduces a formal
inverse $\alpha^{-1}$ with $s(\alpha^{-1})=y$ and
$e(\alpha^{-1})=x$. An {\it edge} in $Q$ is an arrow or the
inverse of an arrow. To each vertex $x$ in $Q$, one associates a
{\it trivial path}, also called {\it trivial walk},
$\varepsilon_x$ with $s(\varepsilon_x)=e(\varepsilon_x)=x$. A {\it
non-trivial walk} $w$ in $Q$ is a finite or infinite product
$$\cdots c_{i+1}c_i\cdots ,$$ where the $c_i$ are edges such that
$e(c_i)=s(c_{i+1})$ for all $i$, whose {\it inverse} $w^{-1}$ is the product
$$\cdots c_i^{-1}c_{i+1}^{-1}\cdots .$$
Such a walk $w$ is called {\it reduced} if
$c_{i+1}\ne c_i^{-1}$ for every $i$, and {\it acyclic}
if it passes through any given vertex in $Q$ at most once. Clearly an acyclic walk is reduced.
If $w$ is a finite or
infinite product
$$\cdots c_i\cdots c_2c_1,$$ then we call $c_1$ the {\it initial
edge} of $w$, we define $s(w)=s(c_1)$, the {\it starting point} of $w$, and we write $w \cdot
\varepsilon_{s(w)}=w$. Dually, if $w$ is a finite or infinite
product $$c_1c_2\cdots c_i\cdots ,$$ then we call $c_1$ the {\it
terminal edge} of $w$, we define $e(w)=e(c_1)$, the {\it ending point} of $w$, and
we write $\varepsilon_{s(w)} \cdot w=w$. An infinite walk $w$ is called {\it
right infinite} if $s(w)$ is defined,  {\it left infinite} if $e(w)$
is defined, and {\it double infinite} if neither $s(w)$ nor $e(w)$
is defined. Clearly, $w$ is finite if and only if both $s(w)$ and
$e(w)$ are defined, and in this case, we call $w$ a walk from $s(w)$
to $e(w)$.

\medskip

A non-trivial walk in $Q$ is called a {\it path} if all of its
edges are arrows. A {\it middle point} of a path is a vertex which
appears in the path but is neither the starting point nor the
ending point. A path is called {\it maximal} if it is not a proper
subpath of any other path in $Q$. Let $x, y$ be vertices in $Q$.
If $Q$ contains a path $p$ from $x$ to $y$, then we say that $x$
is a {\it predecessor} of $y$ which is {\it trivial} or {\it
immediate} if $p$ is a trivial path or an arrow, respectively;
while $y$ is a {\it successor} of $x$ which is {\it trivial} or
{\it immediate} if $p$ is a trivial path or an arrow,
respectively.

\medskip

For $x\in Q_0$, we denote by $x^+$ and $x^-$ the set of arrows
starting in $x$ and the set of arrows ending in $x$, respectively.
We say that $x$ is a {\it sink vertex} or a {\it source vertex} if
$x^+=\emptyset$ or $x^-=\emptyset$, respectively. Moreover, one
says that $Q$ is {\it locally finite} if $x^+$ and $x^-$ are
finite for any $x\in Q_0$, and in this case, one defines the {\it
weight} of $x$ to be the sum of the cardinalities of $x^-$ and
$x^+$. For $x, y\in Q_0$, let $Q(x, y)$ stand for the set of paths
in $Q$ from $x$ to $y$. One says that $Q$ is {\it interval-finite}
if $Q(x, y)$ is finite for any $x, y\in Q_0$. For short, we say
that $Q$ is {\it strongly locally finite} if it is locally finite
and interval-finite. Note that $Q$ contains no oriented cycle in
case it is interval-finite.

\medskip

Let $\it\Sigma$ be a subquiver of $Q$. We shall say that $\Sa$ is {\it full} if, for any vertices $x, y$ in $\Sa$, all the arrows in $Q$ from $x$ to $y$ lie in $\Sa$; {\it convex} if, for any vertices $x, y$ in $\Sa$, all the paths in $Q$ from $x$ to $y$ lie in $\Sa$;
\emph{predecessor-closed} if every path in $Q$ ending in
some vertex in $\Sa$ lies entirely in $\Sa$; and \emph{successor-closed} if every path in $Q$ starting in some vertex in $\Sa$ lies
entirely in $\Sa$.

\medskip

From now on, $k$ denotes an arbitrary field and $Q$ is a strongly locally finite quiver. We shall compose
morphisms in any category from the right to the left, and all tensor products are over $k$.
A {\it representation} $M$ of $Q$ over $k$, or simply a {\it $k$-representation},
consists of a family of
$k$-spaces $M(x)$ with $x\in Q_0,$ and a family of $k$-maps
$M(\alpha): M(x)\to M(y)$ with $\alpha: x\to y$ in $Q_1$. For each
path $\rho$ in $Q$, we write $M(\rho)=\id_{M(x)}$ if
$\rho=\varepsilon_x$ and $M(\rho)=M(\alpha_r)\cdots M(\alpha_1)$ if
$\rho=\alpha_r \cdots \alpha_1$ with $\alpha_1, \ldots, \alpha_r\in
Q_1$. Recall that a morphism $f: M\to N$ of $k$-representations of
$Q$ consists of a family of $k$-maps $f(x): M(x)\to N(x)$ with $x\in
Q_0$ such that $f(y)M(\alpha)=N(\alpha)f(x)$, for every arrow
$\alpha: x\to y$. Let $M$ be a $k$-representation of $Q$. The
\emph{socle} of $M$, written as ${\rm soc}M$, is the
sub-representation of $M$ so that $({\rm soc \hskip 0.5pt}M)(x)$,
for any $x\in Q_0$, is the intersection of the kernels of the maps
$M(\alpha)$ with $\alpha\in x^+$; the {\it radical} of $M$, written as
${\rm rad \hskip 0.5pt}M$, is the sub-representation of $M$ such
that $({\rm rad \hskip 0.5pt}M)(x)$, for any $x\in Q_0$, is the sum
of the images of the maps $M(\beta)$ with $\beta\in x^-$; and the
\emph{top} of $M$, written as ${\rm top\hskip 0.5pt}M$, is the
quotient $M/{\rm rad \hskip 0.5pt}M$. We shall say that ${\rm top
\hskip 0.5pt}M$ is \emph{essential} over $M$ if ${\rm rad \hskip
0.5pt}M$ is superfluous in $M$. Furthermore, the \emph{support} of a
representation $M$, written as ${\rm supp\hskip 0.5pt}M$, is the
full subquiver of $Q$ generated by the vertices $x$ for which
$M(x)\ne 0$. One says that $M$ is {\it sincere} if ${\rm supp\hskip
0.5pt}M=Q$, {\it finitely supported} if ${\rm supp\hskip
0.5pt}M$ is finite, and {\it supported by a subquiver $\Sa$} of $Q$ if ${\rm supp}{\hskip 0.5pt}M \subseteq \Sa$.
Finally, $M$ is called \emph{locally finite
dimensional} if $M(x)$ is of finite $k$-dimension for all $x\in
Q_0$; and \emph{finite dimensional} if $\sum_{x\in Q_0}\, {\rm
dim}M(x)$ is finite. We shall denote by ${\rm Rep}(Q)$ the abelian
category of all $k$-representations of $Q$, which is known to be
hereditary, that is, ${\rm Ext}^2(-, -)$ vanishes; see
\cite[(8.2)]{GR}. Moreover, ${\rm Rep}^b(Q)$, $\rep(Q)$ and
$\rep^b(Q)$ will stand for the full subcate\-gories of ${\rm
Rep}(Q)$ generated by the finitely supported representations, by the
locally finite dimensional representations, and by the finite
dimensional representations, respectively.

\medskip

\begin{Lemma}\label{lemma1.1} Let $M$ be an object in $\Rep(Q)$. If
$\;{\rm supp\hskip 1pt}M$ has no right infinite path, then ${\rm
soc}M$ is essential in $M$. If $\,{\rm supp\hskip 0.5pt}M$ has no
left infinite path, then ${\rm top}\hskip 1pt M$ is essential over
$M$.
\end{Lemma}

\noindent{\it Proof.} We shall prove only the first part. Let $N$ be
a non-zero sub-representation of $M$. Let $x$ be a vertex in ${\rm
supp\hskip 0.5pt}M$ such that $N(x)$ has a non-zero element $v$. If
$\,{\rm supp}\hskip 0.5pt M$ has no right infinite path, then ${\rm
supp\hskip 0.5pt}M$ has a maximal path $\rho$ starting with $x$ such
that $N(\rho)(v) \ne 0$. Note that $N(\rho)(v)$ lies in $N \cap {\rm
soc \hskip 0.5pt}M$. The proof of the lemma is completed.

\medskip

We shall introduce more notation which will be used throughout the
paper. Let $a\in Q_0$. The \emph{simple representation} $S_a$ at
$a$ is defined by $S_a(a)=k\,\varepsilon_a$ and $S_a(x)=0$ for all
vertices $x\ne a$. Moreover, let $P_a$ be the $k$-representation
such that $P_a(x)$, for any $x\in Q_0$, is the $k$-space spanned
by $Q(a, x)$; and $P_a(\alpha): P_a(x)\to P_a(y)$, for $\alpha:
x\to y\in Q_1$, is the $k$-map sending every path $\rho$ to
$\alpha \rho$. Finally, $I_a$ is the $k$-representation such that
$I_a(x)$, for $x\in Q_0$, is the $k$-space spanned by $Q(x, a)$;
and $I_a(\alpha): I_a(x)\to I_a(y)$, for $\alpha: x\to y\in Q_1$,
is the $k$-map sending $\rho \alpha$ to $\rho$ and vanishing on
the paths which do not factor through $\alpha$. Since $Q$ is
interval-finite, the representations $P_a, I_a$ are locally finite
dimensional.

\medskip

\begin{Lemma} \label{lemma1.2} Let $I=I_a\otimes V$ with $a\in Q_0$ and $V$ a
non-zero $k$-space. If $x\in Q_0$ and $\alpha_1, \ldots, \alpha_s$, $s\ge 1$,
are the arrows in ${\rm supp\hskip 0.5pt} I$ starting in $x$, then
$I(x)=W_1\oplus \cdots \oplus W_s$ such that $I(\alpha_i)(W_i)\ne 0$ and $I(\alpha_i)(W_j)=0$,
for $1\le i, j\le n$ with $i\ne j$. As a consequence, $I$ has an essential socle $S_a\otimes V$.
\end{Lemma}

\noindent{\it Proof.} Let $x\in Q_0$ and $\alpha_1, \cdots,
\alpha_s$ be the arrows in ${\rm supp} \hskip 0.5pt I$ starting with
$x$. Then $I_a(x)=N_1\oplus \cdots\oplus N_s$, where $N_i$ is
spanned by the paths $x \rightsquigarrow a$ factoring through
$\alpha_i$. By definition, $I_a(\alpha_i)(N_i)\ne 0$ and
$I_a(\alpha_i)(N_j)=0$, for $1\le i, j\le s$ with $i\ne j$. Since
$I(x)=(N_1\otimes V)\oplus \cdots \oplus (N_n\otimes V)$ and
$I(\alpha_i)=I_a(\alpha_i)\otimes \id_V$, the first part of the
lemma follows. Since $Q$ has no oriented cycle,
$I_a(a)=k\varepsilon_a$, and for any vertex $x$ in ${\rm supp}\hskip
0.5pt I$, we see that ${\rm supp}\hskip 0.5pt I$ has no arrow
starting with $x$ if and only if $x=a$. Thus ${\rm soc}\hskip 0.5pt
I=S_a\otimes V$. Furthermore, since $Q$ is interval-finite, ${\rm
supp\hskip 0.5pt}I$ has no right infinite path. By Lemma
\ref{lemma1.1}, ${\rm soc}\hskip 0.5pt I$ is essential in $I$. The
proof of the lemma is completed.

\medskip

\begin{Prop} \label{prop1.3_shiping} Let $M$ be an object in $\Rep(Q)$ and
$V$ be a $k$-space. For each $a\in Q_0$, there exist
$k$-linear isomorphisms which are natural in $M$ as follows$\,:$
$$\phi_M: {\rm Hom}_{{\rm Rep}(Q)}(M,I_a\otimes V) \longrightarrow
{\rm Hom}_k(M(a),V), \vspace{-2pt}$$ 
$$\vspace{-2pt} \psi_M: {\rm Hom}_{{\rm
Rep}(Q)}(P_a\otimes V, M)\longrightarrow {\rm Hom}_k(V, M(a)).$$
\end{Prop}

\vspace{3pt}

\noindent{\it Proof.} We shall prove only the first part of the
statement. Fix $a\in Q_0$. For each $x\in Q_0$, we have $(I_a\otimes
V)(x)=\oplus_{\rho \in Q(x,a)}(k\rho \otimes V)$. Since $Q$ has no
oriented cycle, $I_a(a)=k\varepsilon_a$. Let $e_{_V}:
I_{a}(a)\otimes V\to V$ be the $k$-isomorphism such that
$e_{_V}(\varepsilon_a\otimes v)=v$, for all $v\in V$. Define
$$\phi_M: {\rm Hom}_{{\rm Rep}(Q)}(M,I_a\otimes V)\to
{\rm Hom}_k(M(a),V): f\mapsto e_{_V}f(a),$$ which is clearly
$k$-linear and natural in $M$. If $f: M\to I_a\otimes V$ is a
morphism such that $\phi _{M}(f)=0$, then $f(a)=0$. Hence ${\rm
soc}(I_a\otimes V) \cap {\rm Im}(f)=0$. Since ${\rm soc \hskip
0.5pt}(I_a\otimes V)$ is essential in $I_a\otimes V$, we get $f=0$. That is, $\phi_M$
is a monomorphism. For proving the surjectivity, let $g: M(a)\rightarrow V$ be a
$k$-map. For each $x\in Q_0$, since $Q(x, a)$ is finite, we
have a $k$-linear map $$f(x): M(x)\to I_a(x)\otimes V: v\mapsto {\Sigma}_{\rho\in
Q(x, a)}\, \rho\otimes g(M(\rho)(v)).$$
Let $\alpha: x\to y$ be an arrow in $Q$. We claim that
$(I_a(\alpha)\otimes \id_V)f(x)=f(y) M(\alpha)$. Indeed, if $Q(x,
a)=\emptyset$, then $Q(y, a)=\emptyset$, and hence $f(x)=0$ and
$f(y)=0$. Assume now that $Q(x, a)\ne \emptyset$. Fix $v\in M(x)$.
We have $f(x)(v)=\sum_{\rho \in Q(x, a)}\, \rho\otimes
g(M(\rho)(v))$. If $\rho=\sigma \alpha$ for some $\sigma\in Q(y,
a)$, then $I_a(\alpha)(\rho)=\sigma$, and otherwise,
$I_a(\alpha)(\rho)=0$. Therefore,
$$(I_a(\alpha)\otimes
\id_V)(f(x)(v))={\Sigma}_{\sigma\in Q(y, a)}\, \sigma\otimes
g(M(\sigma)M(\alpha)(v))=f(y)(M(\alpha)(v)).$$ This establishes our
claim. As a consequence, the $f(x)$ with $x\in Q_0$ form a morphism
$f: M\rightarrow I_{a}\otimes V$ in ${\rm Rep}(Q)$. By definition,
$\phi_M(f)=e_{_V} f(a)=g$. The proof of the proposition is
completed.

\medskip

It follows from the preceding result that $P_a\otimes
V$ is projective and $I_a\otimes V$ is injective in $\Rep(Q)$. Moreover, $P_a$ and $I_a$
are indecomposable. Let ${\rm Inj}\,(Q)$ be the full additive subcategory of ${\rm Rep}(Q)$
generated by the objects isomorphic to
$I_a\otimes V_a$ with $a\in Q_0$ and $V_a$ some $k$-space, and let ${\rm Proj}\,(Q)$ be the one generated by the objects isomorphic to
$P_a\otimes U_a$ with $a\in Q_0$ and
$U_a$ some $k$-space. Moreover, we denote by ${\rm inj}\,(Q)$ and ${\rm proj}\,(Q)$ the full additive subcategories of ${\rm Inj}\,(Q)$
and ${\rm Proj}\,(Q)$, respectively, generated by the locally finite dimensional representations.

\medskip

\begin{Cor} \label{Cor1.4_Shiping} If $M$ is a representation in $\Rep(Q)$, then

\begin{enumerate}[$(1)$]

\item $M\in {\rm Inj}\,(Q)$ if and only if $M$ is injective in
$\Rep(Q)$ with ${\rm soc}\hskip 0.5pt M$ being finitely supported
and essential in $M;$

\vspace{2pt}

\item $M\in {\rm Proj}(Q)$
if and only if $M$ is projective in $\Rep(Q)$ with ${\rm top}\hskip 0.5pt M$ being finitely supported
and essential over $M$.

\end{enumerate}

\end{Cor}

\noindent{\it Proof.} We shall prove only Statement (1). The
necessity follows from Lemma \ref{lemma1.2} and Proposition
\ref{prop1.3_shiping}. For the sufficiency, let $M$ be a non-zero
injective object in ${\rm Rep}(Q)$ such that ${\rm soc}\hskip
0.5pt M$ is finitely supported and essential in $M.$ Then ${\rm
soc \hskip 0.5pt}M=(S_{a_1}\otimes V_1)\oplus\dots\oplus
(S_{a_s}\otimes V_s),$ where the $a_i$ are vertices in $Q$ and the
$V_i$ are non-zero $k$-spaces. Set $I=(I_{a_1}\otimes V_1)\oplus
\cdots \oplus (I_{a_s}\otimes V_s).$ Then ${\rm soc \hskip 0.5pt}
M ={\rm soc \hskip 0.5pt} I$. Observing that  $I$ is injective by
Proposition \ref{prop1.3_shiping}, we have a morphism $f: M\to I$
which acts identically on ${\rm soc \hskip 0.5pt} M$. Since ${\rm
soc \hskip 0.5pt}M$ is essential, $f$ is a monomorphism. Since $M$
is injective, $f$ is a section. Since ${\rm soc}\hskip 0.5pt I$ is
essential, $f$ is an isomorphism. Thus $M\in {\rm Inj}\,(Q)$. The
proof of the corollary is completed.

\medskip

\begin{Defn} \label{defn1.5_shiping} Let $M$ be an object in ${\rm Rep}(Q)$. We say that
$M$ is {\it almost finitely co-presented} if it admits an
injective co-resolution \vspace{-1pt}
$$\xymatrixcolsep{20pt}\xymatrix{0 \ar[r] &  M \ar[r] & I_0 \ar[r] & I_1 \ar[r] & 0}\vspace{-1pt}$$ with $I_0, I_1\in {\rm Inj}(Q)$, and
{\it finitely co-presented} if, in addition, $I_0, I_1\in {\rm inj}(Q)$.
Dually, $M$ is called {\it almost finitely presented} if it admits a projective resolution
\vspace{-1pt}$$\xymatrixcolsep{20pt}\xymatrix{0 \ar[r] &  P_1 \ar[r] & P_0 \ar[r] & M \ar[r] & 0}\vspace{-1pt}$$ with $P_1, P_0\in {\rm Proj}(Q)$, and {\it finitely presented} if, in addition, $P_1, P_0\in {\rm proj}(Q)$.
\end{Defn}

\medskip

We call a quiver {\it top-finite} if every vertex is a successor of finitely many
pre-fixed vertices, and {\it socle-finite} if every vertex is a
predecessor of finitely many pre-fixed vertices.

\medskip

\begin{Lemma} \label{Lemma1.6_shiping} Let $M$ be a representation of $\Rep(Q)$.

\begin{enumerate}[$(1)$]

\vspace{-1pt}

\item If $M$ is almost finitely co-presented, then ${\rm soc}\,M$
is finitely supported and essential in $M$, and $\,{\rm supp\hskip
0.5pt}M$ is socle-finite with no right infinite path.

\vspace{1pt}

\item If $M$ is almost finitely presented, then ${\rm top}\hskip
0.5pt M$ is finitely supported and essential over $M$, and $\,{\rm
supp\hskip 0.5pt}M$ is top-finite with no left infinite path.

\vspace{1pt}
\end{enumerate}
\end{Lemma}

\noindent{\it Proof.} We shall prove only (1). Assume that $M$ is non-zero and
almost finitely co-presented. Let $M\to I$ be the  injective envelope of
$M$, where $I=(I_{a_1}\otimes U_1)\oplus \cdots \oplus
(I_{a_r}\otimes U_r)$ with $a_1, \ldots, a_r\in Q_0$ and $U_1,
\ldots, U_r$ some non-zero $k$-spaces. Then ${\rm soc \hskip 0.5pt}
M$ is essential in $M$ and supported by the vertices $a_1, \ldots,
a_r$. Let $x$ be a vertex in ${\rm supp\hskip 0.5pt}M$. Choose some
non-zero element $v\in M(x)$ and let $L$ be the sub-representation
of $M$ genera\-ted by $v$. Since $ L \cap {\rm soc \hskip 0.5pt}M
\ne 0$, ${\rm supp\hskip 0.5pt} L$ contains some $a_i$ with $1\le
i\le r$. Therefore, ${\rm supp\hskip 0.5pt} L$ has a path from $x$
to $a_i$. That is, $x$ is a predecessor of $a_i$ in ${\rm supp\hskip
0.5pt}M$. Therefore, $\,{\rm supp\hskip 0.5pt}M$ is socle-finite.
Since $Q$ is interval-finite, we see that ${\rm supp\hskip 0.5pt}M$
has no right infinite path. The proof of the lemma is completed.

\medskip

The following result states some useful combinatorial properties of the finitely presented or finitely co-presented  but infinite dimensional representations.

\medskip

\begin{Cor} \label{newcor1.9} Let $M$ be an infinite dimensional representation in ${\rm rep}(Q)$.

\vspace{-2pt}

\begin{enumerate}[$(1)$]

\item If $M$ is finitely presented, then
$\,{\rm supp\hskip 0.5pt}M$ contains a right infinite path.

\vspace{2pt}

\item If $M$ is finitely co-presented, then $\,{\rm supp\hskip
0.5pt}M$ contains a left infinite path.

\end{enumerate}
\end{Cor}

\noindent{\it Proof.} We shall prove only Statement (1). Assume that
$M$ is finitely presented. By Lemma \ref{Lemma1.6_shiping}, ${\rm
supp}\hspace{0.5pt}M$ is top-finite. Since ${\rm
supp}\hspace{0.5pt}M$ is infinite and locally finite,  by
K\"{o}nig's lemma, it has a right infinite path. The proof of the
corollary is completed.

\medskip

We denote by ${\rm Rep}^-(Q)$ the full subcategory of ${\rm Rep}(Q)$
generated by the almost finitely co-presented objects and by ${\rm
Rep}^+(Q)$ the one generated by the almost finitely presented
objects.

\medskip

\begin{Prop} \label{prop1.7_Shiping} The categories ${\rm Rep}^+(Q)$ and ${\rm
Rep}^-(Q)$ are abelian, hereditary, and closed under extensions in
${\rm Rep}(Q)$.
\end{Prop}

\noindent{\it Proof.} We shall consider only ${\rm Rep}^-(Q)$.
Given a morphism $f: I_1\to I_2$ in ${\rm Inj}\,(Q)$, set $I={\rm
Im}(f)$ and consider the short exact sequence
$$(*) \;\;\xymatrixcolsep{15 pt}\xymatrix{0 \ar[r] &  I \ar[r] & I_2 \ar[r] & J \ar[r] & 0}.$$
Since ${\rm Rep}(Q)$ is hereditary, $I, J$ are injective, and hence
the sequence $(*)$ splits. In particular, ${\rm soc \hskip
0.5pt}I_2\cong {\rm soc \hskip 0.5pt}I \,\oplus \,{\rm soc \hskip
0.5pt}J$. Since ${\rm soc \hskip 0.5pt}I_2$ is finitely supported
and essential in $I_2$, we see that ${\rm soc \hskip 0.5pt}J$ is
finitely supported and essential in $J$. By Corollary
\ref{Cor1.4_Shiping}, $J\in {\rm Inj}\,(Q)$. Now it follows from the
dual of Proposition 2.1 in \cite{Aus} that ${\rm Rep}^-(Q)$ is
abelian and closed under extensions in ${\rm Rep}(Q)$. Moreover,
${\rm Rep}^-(Q)$ is hereditary since ${\rm Rep}(Q)$ is so. The proof
of the proposition is completed.

\medskip

Let $\it\Sigma$ be a subquiver of $Q$. For each representation $M\in {\rm Rep}(Q)$, we define an object $M_{_{\hskip -1pt \it\Sigma}}\in \Rep(Q)$, called the \emph{restriction} of $M$ to $\it\Sigma$, by setting $M_{_{\hskip
-0.8pt \it\Sigma\hskip 0.3pt}}(\rho)=M(\rho)$ if $\rho\in
{\it\Sigma}$; and $M_{_{\hskip -1pt \it\Sigma}}(\rho)=0$ otherwise,
where $\rho$ ranges over $Q_0\cup Q_1$. For a morphism $f: M\to N$ in $\Rep(Q)$,
we define its \emph{restriction} to $\Sa$, written as
$f_{_{\it\Sigma}}: M_{_{\hskip -1pt \it\Sigma}} \to N_{_{\hskip -1pt
\it\Sigma}}$, by setting $f_{_{\it\Sigma}}\,(x)=f(x)$ if $x\in
{\it\Sigma}_0$, and $f_{_{\hskip -0.5pt\it\Sigma}}\,(x)=0$
otherwise, where $x$ ranges over $Q_0$.

\medskip

\begin{Defn} \label{rel-proj} Let $\Sa$ be a full subquiver of $Q$. A representation $M\in \Rep(Q)$ is called \emph{projective} or {\it injective restricted to} $\Sa$  if $M_{_{\hspace{-1pt}\it \Sigma}}\in {\rm Proj}\,(Q)$ or $M_{_{\hspace{-1pt}\it \Sigma}}\in {\rm Inj}\,(Q)$, respectively.
\end{Defn}

\medskip

Let $\Sa$ be a full subquiver of $Q$. The {\it complement} of $\Sa$
is the full subquiver of $Q$ generated by the vertices not in $\Sa$,
while the {\it augmented complement} of $\Sa$ is the full subquiver
of $Q$ generated by the vertices and the arrows not in $\Sa$. Note
that a vertex $x\in \Sa_0$ lies in the augmented complement of $\Sa$
if and only if there exists an edge $x$ --- $y$ with $y\not\in
\Sa_0$. Since $Q$ is  locally finite, the augmented complement of
$\Sa$ is finite if and only if the complement is finite. We shall
say that $\Sa$ is \emph{co-finite in} $Q$ if its complement in $Q$
is finite.

\medskip

\begin{Lemma} \label{lemma1.8_shiping} If $M\in \Rep(Q)$ is
injective restricted to some full subquiver $\Sa$ of $Q$, then it
is injective restricted to any co-finite predecessor-closed
subquiver of $\Sa$.
\end{Lemma}

\noindent{\it Proof.} Let $M\in \Rep(Q)$ such that
$M_{_{\hspace{-1pt}\it \Sigma}}\in {\rm Inj}\,(Q)$, where $\Sa$ is a full
subquiver of $Q$. With no loss of generality, we may assume that
$M_{_{\hspace{-1pt}\it \Sigma}}=I$, where $I=I_a\otimes U$ with $a\in {\Sa}_0$ and
$U$ a non-zero $k$-space. Then ${\rm supp\hskip 0.5pt} I \subseteq
\Sa$. Let $\it\Theta$ be a co-finite predecessor-closed subquiver of
$\Sa$. Then $M_{_{\it\Theta}}=I_{_{\it\Delta}}$, where
$\Da={\it\Theta}\cap {\rm supp\hskip 0.5pt} I.$ Observe that $\Da$
is co-finite and predecessor-closed in ${\rm supp\hskip 0.5pt} I$.
As a consequence, $I_{_{\it\Delta}}$ is a quotient of $I$. Since
$\Rep(Q)$ is hereditary, $I_{_{\it\Delta}}$ is injective. Note that
${\rm supp}\hskip 0.5pt I$ has no right infinite path since $Q$ is
interval-finite. By Lemma \ref{lemma1.1}, ${\rm soc}\hspace{0.6pt}
I_{_{\it\Delta}}$ is essential in $I_{_{\it\Delta}}$. If $a\in \Da$,
then $I_{_{\it\Delta}}=I_a\otimes U$. Otherwise, by Lemma
\ref{lemma1.2}, for any vertex $x$ in the support of ${\rm soc
\hskip 0.8pt} I_{_{\it\Delta}}$, there exists an arrow $x\to y$ in
${\rm supp\hskip 0.5pt} I$ with $y\not\in \Da$. Since $\Da$ is
co-finite in ${\rm supp\hskip 0.5pt} I$, we see that ${\rm soc
\hskip 0.8pt} I_{_{\it\Delta}}$ is finitely supported. By Corollary
\ref{Cor1.4_Shiping}, $I_{_{\it\Delta}}\in {\rm Inj}\,(Q)$. That is,
$M_{_{\it\Theta}} \in {\rm Inj}\,(Q).$  The proof of the lemma is
completed.

\medskip

\begin{Prop} \label{prop1.9_shiping} The intersection of $\,{\rm Rep}^{-\hspace{-2pt}}(Q)$ and ${\rm Rep}^{+\hspace{-1pt}}(Q)$ is  ${\rm Rep}^{b \hskip -1pt}(Q)$.
\end{Prop}

\noindent{\it Proof.} Let $M\in {\rm Rep}(Q)$. If $M$ is almost
finitely presented and almost finitely co-presented then, by Lemma
\ref{Lemma1.6_shiping}, ${\rm supp\hskip 0.5pt}M$ is both
socle-finite and top-finite. Since $Q$ is interval-finite, ${\rm
supp\hskip 0.5pt}M$ is finite. Conversely, assume that ${\rm
supp\hskip 0.5pt}M$ is finite. In particular, ${\rm soc \hskip
0.5pt}M=\oplus_{i=1}^n (S_{a_i}\otimes U_i)$, where $a_1, \ldots,
a_n\in Q_0$ and $U_1, \ldots, U_n$ are non-zero $k$-spaces.
Therefore, $M$ has an injective envelope $M\to I$, where
$I=\oplus_{i=1}^n (I_{a_i}\otimes U_i)\in {\rm Inj}(Q)$. Consider
the short exact sequence  $$\xymatrixcolsep{15 pt}(*)\;\;\xymatrix{0
\ar[r] &  M \ar[r] & I \ar[r] & J \ar[r] & 0}$$ where $J$ is
injective. By Lemma \ref{Lemma1.6_shiping}, ${\rm supp\hskip 0.5pt}
I$ is socle-finite. Since $Q$ is interval-finite, ${\rm supp\hskip
0.5pt} I$ has no right infinite path. Since ${\rm supp\hskip 0.5pt}
J \subseteq {\rm supp\hskip 0.5pt} I$, by Lemma \ref{lemma1.1},
${\rm soc \hskip 0.5pt}J$ is essential in $J$. Denote by $\Da$ the
successor-closed subquiver of ${\rm supp}I$ generated by ${\rm
supp}M$. Since  ${\rm supp\hskip 0.5pt} I$ is socle-finite, $\Da$
is finite. Let $\it\Sigma$ be the complement of $\Da$ in ${\rm
supp\hskip 0.5pt}I$. Then $\Sa$ is predecessor-closed in ${\rm
supp}\hskip 0.8pt I$. Restricting the sequence $(*)$ to $\it\Sigma$
yields a short exact sequence
$$ \xymatrixcolsep{15 pt}\xymatrix{0 \ar[r] &  M_{_{\hspace{-1pt}\it \Sigma}} \ar[r] & I_{_{\it\Sigma}} \ar[r] & J_{_{\it\Sigma}} \ar[r] & 0}. $$ Since $M_{_{\hspace{-1pt}\it \Sigma}}=0$,
we get $J_{_{\it\Sigma}} \cong I_{_{\it\Sigma}}.$ By Lemma
\ref{lemma1.8_shiping}, ${\rm
soc}\hspace{0.2pt}(\hspace{-0.2pt}J_{_{\it\Sigma}})$ is finite. Let
$x$ be a vertex in the support of ${\rm soc}\hskip 0.5pt J$. If
$x\in \Sa$, then $x$ lies in the support of ${\rm soc}\hskip 0.5pt
J_{_{\it\Sigma}}$. Otherwise, $x\in \Da$. This shows that ${\rm
soc}\hskip 0.5pt J$ is finitely supported. By Corollary
\ref{Cor1.4_Shiping}, $J\in {\rm Inj}\,(Q)$. That is, $M$ is almost
finitely co-presented. Dually, $M$ is almost finitely presented. The
proof of the proposition is completed.

\medskip

We are ready to have a criterion for a representation to be almost finitely presented or co-presented.

\medskip

\begin{Theo} \label{theo1.10_shiping} Let $Q$ be a strongly locally finite quiver, and let $M\in \Rep(Q)$.

\begin{enumerate}[$(1)$]

\vspace{-2.5pt}

\item $M$ is almost finitely presented if and only if $M$ is
projective restricted to some co-finite successor-closed subquiver
of ${\rm supp\hskip 0.5pt}M$.

\vspace{1pt}

\item $M$ is almost finitely
co-presented if and only if $M$ is injective restricted to some
co-finite predecessor-closed subquiver of ${\rm supp\hskip 0.5pt}M$.

\vspace{-1pt}

\end{enumerate}
\end{Theo}

\noindent{\it Proof.} We shall prove only (2). Firstly, assume that
$M_{_{\hspace{-1pt}\it \Sigma}}\in {\rm Inj}\,(Q)$, where $\it\Sigma$ is a
co-finite predecessor-closed subquiver of ${\rm supp\hskip 0.5pt}M$.
Let $\it\Omega$ be the complement of $\it\Sigma$ in ${\rm supp\hskip
0.5pt}M$. Then $\Oa$ is finite and successor-closed in ${\rm
supp\hskip 0.5pt}M$. Thus $M_{_{\hspace{-0.8pt}\it\Omega}}$ is a finitely supported
sub-representation of $M$, and we have a short exact sequence
$$\xymatrixcolsep{15 pt}\xymatrix{0 \ar[r] &  M_{_{\hspace{-0.8pt}\it\Omega}}  \ar[r] & M \ar[r] & M_{_{\hspace{-1pt}\it \Sigma}} \ar[r] & 0}$$
 in ${\rm Rep}(Q)$. By Propositions \ref{prop1.7_Shiping} and
\ref{prop1.9_shiping}, $M$ is almost finitely co-presented.
Conversely, assume that $M$ admits a minimal injective resolution
$$\xymatrixcolsep{15 pt}\xymatrix{
 0\ar[r] & M \ar[r]^f & I \ar[r]^g & J \ar[r] & 0}$$
with $I=\oplus_{i=1}^m\, (I_{a_i}\otimes U_i)$ and
$J=\oplus_{j=1}^n\, (I_{b_j}\otimes V_j)$, where $a_i, b_j\in Q_0$
and the $U_i, V_j$ are $k$-spaces. Let $\it\Omega$ be the convex
hull in ${\rm supp\hskip 0.5pt}I$ generated by $a_1, \ldots, a_m,
b_1, \ldots, b_n$. It is easy to see that $\Oa$ is finite and successor-closed in
${\rm supp\hskip 0.5pt}I$. Let $\Da$ be the complement of $\Oa$ in
${\rm supp\hskip 0.5pt}I$. Then $\Da$ is co-finite and
predecessor-closed in ${\rm supp\hskip 0.5pt}I$. We claim that the
short exact sequence
$$\xymatrixcolsep{15 pt}\xymatrix{
 0\ar[r] & M_{_{\it\Delta}} \ar[r]^{f_{_{\it\Delta}}} & I_{_{\it\Delta}}
 \ar[r]^{g_{_{\it\Delta}}} & J_{_{\it\Delta}} \ar[r] & 0}$$
splits. Indeed, since $I_{_{\it\Delta}}, \, J_{_{\it\Delta}}\in {\rm
Inj}\,(Q)$ by Lemma \ref{lemma1.8_shiping}, it suffices to show that
$g_{_{\it\Delta}}$ induces a surjective map from ${\rm soc \hskip
0.5pt}I_{_{\it\Delta}}$ to ${\rm soc \hskip 0.5pt}J_{_{\it\Delta}}$.
For this purpose, fix a vertex $x$ in the support of ${\rm soc}
J_{_{\it\Delta}}$ and a non-zero element $v$ in $({\rm soc \hskip
0.5pt}J_{_{\it\Delta}})(x)$. Since $g(x)$ is surjective, $v=g(x)(u)$
for some $u\in I(x)$. Let $\alpha_i: x\to x_i$, $i=1, \ldots, s$, be
the arrows in ${\rm supp\hskip 0.5pt} I$ starting with $x$, where
$x_1, \ldots, x_r\in \Oa$ with $0\le r\le s$ and $x_{r+1}, \cdots,
x_s\in \Da$. It follows from Lemma \ref{lemma1.2} that
$I(x)=W_1\oplus \cdots \oplus W_s$ such that $I(\alpha_i)(W_j)=0$
for $1\le i, j\le s$ with $i\ne j$. Write $u=u_1+\cdots +u_s$ with
$u_i\in W_i$. For $1\le j\le r$, since $x_j\not\in \Da$, we have
$I_{_{\it\Delta}}(\alpha_j)=0$, and hence
$I_{_{\it\Delta}}(\alpha_j)(u_1+\cdots+u_r)=0$. Since
$I(\alpha_j)(u_1+\cdots+u_r)=0$ for $r<j\le s$, we get
$u_1+\cdots+u_r\in ({\rm soc \hskip 0.5pt} I_{_{\it\Delta}})(x)$. If
$r=s$, then $u\in ({\rm soc \hskip 0.5pt}I_{_{\it\Delta}})(x)$.
Otherwise, consider $g(x)(u_{r+1}+\cdots+u_s)\in J(x)$. For $1\le
j\le r$, we have
$J(\alpha_j)(g(x)(u_{r+1}+\cdots+u_s))=g(x_j)(I(\alpha_j)(u_{r+1}+\cdots+u_s))=g(x_j)(0)=0$.
For $r<j\le s$, since $\alpha_j\in \Da$ and $v\in {\rm soc \hskip
0.5pt}J_{_{\it\Delta}}$, we have
$J(\alpha_j)(v)=J_{_{\it\Delta}}(\alpha_j)(v)=0$. This yields
\vspace{-2pt}$$\hspace{40pt}\begin{array}{rcl}
J(\alpha_j)(g(x)(u_{r+1}+\cdots +u_s))&=&
g(x_j)(I(\alpha_j)(u_{r+1}+\cdots+u_s))\\
&=&g(x_j)(I(\alpha_j)(u_1+\cdots+u_s))\\
&=& J(\alpha_j)(g(x)(u_1+\cdots+u_s))\\
&=&J(\alpha_j)(v)\\
&=&0.\end{array}$$ Since ${\rm supp\hskip 0.5pt}J\subseteq {\rm
supp\hskip 0.5pt}I$, we have $g(x)(u_{r+1}+\cdots+u_s)\in ({\rm soc
\hskip 0.5pt}J)(x)$, and hence $g(x)(u_{r+1}+\cdots+u_s)=0$ since
$x\ne b_j$ for $1\le j\le n$. This shows that $u_1+\cdots+u_r$ is a
pre-image of $v$ by $g_{_{\it\Delta}}$ in $({\rm soc \hskip
0.5pt}I_{_{\it\Delta}})(x)$. Our claim is established. As a
consequence, $M_{_{\it\Delta}}\in {\rm Inj}\,(Q)$. Finally,
$\Sa=\Da\cap \,{\rm supp\hskip 0.5pt}M$ is co-finite and
predecessor-closed in ${\rm supp\hskip 0.5pt}M$ such that
$M_{_{\hspace{-1pt}\it \Sigma}}=M_{_{\it\Delta}}$. That is, $M$ is injective
restricted to $\Sa$. The proof of the theorem is completed.

\medskip

Combined with Theorem \ref{theo1.10_shiping}, the following result
and its dual allow us to reduce the study of almost finitely presented or co-presented representations of
$Q$ to the study of representations of finite quivers.

\medskip

\begin{Theo} \label{restriction} Let $Q$ be a strongly locally finite quiver with a full subquiver $Q\hspace{0.2pt}'$. Let $M, N\in {\rm Rep}^-(Q)$ such that $M\oplus N$ is supported by $Q\hspace{0.2pt}'$ and injective restricted to a full subquiver $\Sa$ of $\,Q\hspace{0.2pt}'$, and let
$\Oa$ be a successor-closed subquiver of $Q\hspace{0.2pt}'$ which contains the
socle-support of $(M\oplus N)_{_{\hspace{-1pt}\it\Sigma}}$ and the augmented
complement of $\Sa$ in $Q\hspace{0.2pt}'$.

\begin{enumerate}[$(1)$]

\item  $M\cong N$ if and only if $M_{_{\hspace{-0.8pt}\it\Omega}}\cong N_{_{\hspace{-1.5pt}\it\Omega}}$.

\vspace{1pt}

\item  $M$ is indecomposable if and only if $M_{_{\hspace{-0.8pt}\it\Omega}}$ is indecomposable.

\vspace{1pt}

\item A morphism $f: M\to N$ is a section or a retraction if and only if $f_{_{\it\Omega}}$ is a section or a retraction, respectively.

\vspace{1pt}

\end{enumerate}
\end{Theo}

\noindent{\it Proof.} First of all, we show that the $k$-linear map
$$\phi: {\rm Hom}_{\,{\rm Rep}(Q)}(M, N)\to {\rm Hom}_{\,{\rm Rep}(Q)}(M_{_{\hspace{-0.8pt}\it\Omega}}, N_{_{\hspace{-1.5pt}\it\Omega}}): f\mapsto f_{_{\it\Omega}}$$
is bijective. Let $f: M\to N$ be a morphism such that
$f_{_{\it\Omega}}=0$. If $x$ is a vertex in the socle-support of
$N$, then $x$ lies in the socle-support of $N_{_{\it\Sigma}}$ or in
the complement of $\Sa$ in $Q\hspace{0.2pt}'$, and hence $x\in \Oa$. Thus
$f(x)=0$. Hence ${\rm Im}(f)\cap {\rm soc \hskip 0.5pt}N=0$. Since
${\rm soc \hskip 0.5pt}N$ is essential in $N$, we get $f=0$. That
is, $\phi$ is injective.

Next, let $\Da=\Sa \,\cap\, \Oa$, which is successor-closed in
$\Sa$. Hence $M_{_{\it\Delta}}$ and $N_{_{\it\Delta}}$ are
sub-representations of $M_{_{\hspace{-1pt}\it \Sigma}}$ and $N_{_{\it\Sigma}}$,
respectively. Let $g: M_{_{\hspace{-0.8pt}\it\Omega}}\to N_{_{\hspace{-1.5pt}\it\Omega}}$ be a
morphism. Consider the restriction $g_{_{\it\Delta}}:
M_{_{\it\Delta}}\to N_{_{\it\Delta}}$. Since $N_{_{\it\Sigma}}$ is
injective, we have a morphism $h: M_{_{\hspace{-1pt}\it \Sigma}}\to
N_{_{\it\Sigma}}$ such that $h_{_{\it\Delta}}=g_{_{\it\Delta}}.$ For
any $x\in Q_0$, set $f(x)=g(x)$ if $x\in \Oa$; $f(x)=h(x)$ if $x\in
\Sa$; and $f(x)=0$ if $x\not\in Q\hspace{0.2pt}'$. Since every arrow in $Q\hspace{0.2pt}'$ lies
in $\Sa$ or $\Oa$, we verify easily that $f=\{f(x)\mid x\in Q_0\}$
is a morphism from $M$ to $N$ in $\Rep(Q)$ such that
$f_{_{\it\Omega}}=g$. That is, $\phi$ is surjective.

Specializing to the case where $N=M$, we get ${\rm End}(M)\cong {\rm
End}(M_{_{\hspace{-0.8pt}\it\Omega}})$. Since ${\rm Rep}(Q)$ is abelian, an object
is indecomposable if and only if its endomorphism algebra has only
trivial idempotents. As a consequence, $M$ is indecomposable if and
only if $M_{_{\hspace{-0.8pt}\it\Omega}}$ is indecomposable. This establishes
Statement (2).

Let $f: M\to N$ be a morphism in $\Rep(Q)$. If $f$ is a section,
then it is evident that $f_{_{\it\Omega}}$ is a section. Assume that
$f_{_{\it\Omega}}$ is a section. Let $g:
N_{_{\hspace{-1.5pt}\it\Omega}} \to M_{_{\hspace{-0.8pt}\it\Omega}}$
be a morphism such that
$gf_{_{\it\Omega}}=1_{M_{_{\hspace{-0.8pt}\it\Omega}}}$. Since
$\phi$ is surjective, there exists a morphism $h: N\to M$ such that
$h_{_{\it\Omega}}=g$. This yields
$(hf)_{_{\it\Omega}}=(1_M)_{_{\it\Omega}}$, and thus $hf=1_M$. The
first part of Statement (3) is established, and the second part
follows in the same way. As a consequence, $f$ is an isomorphism if
and only if $f_{_{\it\Omega}}$ is an isomorphism. This proves
Statement (1). The proof of the theorem is completed.

\medskip

The rest of this section is devoted to study the finitely presented
or co-presented representations. Let ${\rm
rep}^+(Q)$ stand for the category of finitely presented
representations of $Q$ and ${\rm rep}^-(Q)$ for that of finitely co-presented ones.

\medskip

\begin{Lemma} \label{Prop1.13_shiping} Let $L, M$ be representations in ${\rm rep}(Q)$.

\begin{enumerate}[$(1)$]

\vspace{-1pt}

\item If $M\in {\rm rep}^+(Q)$, then ${\rm Ext}^i_{{\rm rep}(Q)}(M,
L)$ is finite dimensional for $i\ge 0$.

\vspace{1pt}

\item If $M\in {\rm rep}^-(Q)$, then ${\rm Ext}^i_{{\rm rep}(Q)}(L, M)$ is finite dimensional for $i\ge 0$.

\vspace{-0pt}

\end{enumerate} \end{Lemma}

\noindent{\it Proof.} We shall prove only Statement (2). Let $M$ be
an object in ${\rm rep}^-(Q)$ with a minimal injective
co-resolution $\xymatrixcolsep{15 pt}\xymatrix{
0\ar[r] & M \ar[r] & I \ar[r] & J \ar[r] & 0},$ where $I, J\in {\rm inj}(Q)$. Applying
${\rm Hom}_{{\rm rep}(Q)}(L, -)$ yields an exact sequence
\vspace{-3pt}$$\xymatrixcolsep{15 pt} \xymatrix{0\ar[r] & {\rm Hom}(L, M) \ar[r] & {\rm Hom}(L, I) \ar[r] & {\rm Hom}(L, J) \ar[r] & {\rm Ext}^1(L, M)\ar[r] & 0}.\vspace{-3pt}$$
Since $L\in \rep(Q)$, by Proposition \ref{prop1.3_shiping}, ${\rm
Hom}_{{\rm rep}(Q)}(L, I)$ and ${\rm Hom}_{{\rm rep}(Q)}(L, J)$ are
finite dimensional. The proof of the lemma is completed.

\medskip

One says that an additive $k$-category is {\it Hom-finite} if the Hom-spaces are of finite $k$-dimension and that
an abelian $k$-category is {\it Ext-finite} if the Ext-spaces are of finite $k$-dimension. Note that a Hom-finite additive $k$-category is
Krull-Schmidt if its idempotents split. In particular, a Hom-finite abelian $k$-category is Krull-Schmidt. The following result is an immediate consequence of Propositions \ref{prop1.7_Shiping}, \ref{prop1.9_shiping} and \ref{Prop1.13_shiping}.

\medskip

\begin{Prop} \label{cat-fin-pres}
The $k$-categories ${\rm rep}^+(Q)$ and ${\rm rep}^-(Q)$ are
Ext-finite, heredi\-tary, and  abelian. Moreover, they are
extension-closed in $\rep(Q)$ and their intersection is
$\rep^b(Q)$.

\end{Prop}

\medskip

\noindent{\sc Remark.} If $Q$ has no left infinite path or no
right infinite path,  then we have $\rep^-(Q)=\rep^b(Q)$ or
$\rep^+\hspace{-1pt}(Q)=\rep^b(Q)$, respectively. As a
consequence, if $Q$ has no infinite path, then
$\rep^+\hspace{-1pt}(Q)=\rep^b(Q)=\rep^-(Q)$.

\medskip

We shall describe the projective objects and the finite
dimensional injective objects in $\rep^+\hspace{-1pt}(Q)$. For
this purpose, denote by $Q^+$ the full subquiver of $Q$ generated
by the vertices which are not ending point of any left infinite
path.

\medskip

\begin{Prop}\label{cor1.14shiping} Let $M$ be an indecomposable representation in
$\rep^+\hspace{-1pt}(Q)$.

\begin{enumerate}[$(1)$]

\item $M$ is a projective in $\rep^+\hspace{-1pt}(Q)$  if and only
if $M\cong P_x$ for some $x\in Q_0$.

\vspace{1pt}

\item $M$ is finite dimensional and injective in
$\rep^+\hspace{-1pt}(Q)$ if and only if $M\cong I_x$ for some
$x\in Q^+$.

\vspace{1pt}

\item If $M\cong S_x$ with $x\in Q_0$, then it has an injective
hull in $\rep^+\hspace{-1pt}(Q)$ if and only if $x\in Q^+$.

\end{enumerate}\end{Prop}

\noindent{\it Proof.} (1)  By definition, $M$ has a projective
cover $f: P\to M$, where $P$ is an object in ${\rm
proj}\hspace{0.4pt}(Q)$. If $M$ is projective in
$\rep^+\hspace{-1pt}(Q)$, then $f$ is an isomorphism. Since $M$ is
indecomposable, we have $P\cong P_x$ for some $x\in Q_0$.

(2) Let $x\in Q^+$. Since  $Q$ is locally finite, $x$ admits only
finitely many predecessors in $Q$, and hence $I_x$ is finite
dimensional. Since $I_x$ is injective in $\rep(Q)$, it is
injective in $\rep^+\hspace{-1pt}(Q)$. Suppose conversely that $M$
is a finite dimensional injective object in
$\rep^+\hspace{-1pt}(Q)$. Then ${\rm
soc}\hspace{0.4pt}M=\oplus_{i=1}^r\,(S_{x_i}\otimes U_i)$ with
$x_i\in Q_0$ and $U_i$ some finite dimensional non-zero
$k$-spaces. We have an essential monomorphism $f: M\to I$ in
$\rep(Q)$, where $I=\oplus_{i=1}^r \,(I_{x_i}\otimes U_i)$. In
particular, ${\rm supp}\hspace{0.4pt}M\subseteq {\rm
supp}\hspace{0.4pt}I$. Fix a vertex $y$ in ${\rm
supp}\hspace{0.4pt}I$. Let $\Sa$ be the convex subquiver of ${\rm
supp}\hspace{0.4pt}I$ generated by ${\rm supp}\hspace{0.4pt}M$ and
$y$. Since ${\rm supp}\hspace{0.4pt}M$ is finite, so is $\Sa$. By
restriction, we get an essential monomorphism $f_{_{\it\Sigma}}:
M\to I_{_{\it\Sigma}}$ in $\rep(\Sa)$. Observing that $M$ is
injective in $\rep(\Sa)$, we see that $f_{_{\it\Sigma}}$ is an
isomorphism. In particular, $M(y)\cong I_{_{\it\Sigma}}(y)=I(y)\ne
0$. Thus ${\rm supp}\hspace{0.4pt}M = {\rm supp}\hspace{0.4pt}I$.
Therefore, $f$ is an essential monomorphism in
$\rep^+\hspace{-1pt}(Q)$. Since $M$ is injective in
$\rep^+\hspace{-1pt}(Q)$, $f$ is an isomorphism, and consequently,
$r=1$ and $x_1\in Q^+$.

\vspace{1pt}

(3) If $x\in Q^+$, then $I_x$ is the injective hull of $S_x$ in $\rep^+\hspace{-1pt}(Q)$. Conversely, suppose that $S_x$ has an injective hull $I$ in $\rep^+\hspace{-1pt}(Q)$. Since $S_x$ is essential in $I$, every vertex in ${\rm supp}\hspace{0.5pt}I$ is a predecessor of $x$ in $Q$. Thus ${\rm supp}\hspace{0.5pt}I$ is both top-finite and socle-finite. As a consequence, $I$ is finite dimensional. We now deduce from Statement (2) that $I=I_x$, and hence $x\in Q^+$. The proof of the proposition
is completed.

\medskip

\noindent{\sc Example.} Let $Q$ be the following infinite quiver

\vspace{-8pt}
$$\xymatrixcolsep{15 pt}\xymatrix{
1\ar[r] & 2 \ar[r] & 3 \ar[r] & \cdots \ar[r] & n \ar[r] &
\cdots}.\vspace{1pt}$$ Since $\rep^+\hspace{2pt}(Q)$ is
Krull-Schmidt, using the description of the indecomposable
representations in $\rep^+\hspace{-1pt}(Q)$ given in (\ref{5.1}),
one can verify that $P_1$ is an injective object in
$\rep^+\hspace{-1pt}(Q)$. It is clear that $P_1\not\in {\rm
inj}(Q)$.

\medskip

To conclude this section, we construct a duality between ${\rm
proj}(Q)$ and ${\rm inj}(Q)$. The proofs will be left out since
this is similar to the construction in the finite case. Consider
the opposite quiver $Q^{\rm op}$ of $Q$ which is defined in such a
way that every vertex $x$ in $Q$ corresponds to a vertex $x^{\rm
o}$ in $Q^{\hspace{0.4pt} \rm op}$ and every arrow $\alpha: x\to
y$ corresponds to an arrow $\alpha^{\rm o}: y^{\rm o}\to x^{\rm
o}$ in $Q^{\hspace{0.4pt} \rm op}$. If $p=\alpha_n\cdots \alpha_1$
is a path in $Q$ from $x$ to $y$, then we write $p^{\rm
o}=\alpha_1^{\rm o}\cdots \alpha_n^{\rm o}\vspace{1pt}$, the
corresponding path in $Q^{\hspace{0.4pt} \rm op}$ from $y^{\rm o}$
to $x^{\rm o}$. For any given object $M\in {\rm rep}(Q)$, we
define $DM\in {\rm rep}(Q^{\hspace{0.4pt} \rm op})$ by setting
$(DM)(x^{\rm o})={\rm Hom}_k(M(x), k)$ for each vertex $x^{\rm o}$
and $(DM)(\alpha^{\rm o})$ to be the transpose of $M(\alpha)$, for
each arrow $\alpha^{\rm o}$. For a morphism $f:M\to N$ in ${\rm
rep}(Q)$, we define a morphism $Df: DN\to DM$ in ${\rm
rep}(Q^{\hspace{0.4pt} \rm op})$ by setting $(Df)(x^{\rm o})$, for
each vertex $x^{\rm o}$ in $Q^{\hspace{0.4pt} \rm op}$, to be the
transpose of $f(x)$.

\medskip

\begin{Lemma} \label{piduality}

The functor $D: {\rm rep}(Q)\to {\rm rep}(Q^{\hspace{0.4pt} \rm
op})$ is a duality such that $DI_x\cong P_{x^{\rm o}}$ and
$DP_x\cong I_{x^{\rm o}}$, for all $x\in Q_0$.

\end{Lemma}

\medskip

For the rest of the paper, put $A=kQ$, the path algebra of $Q$
over $k$. Note that $A$ has a complete set of primitive orthogonal
idempotents $\{\varepsilon_x\mid \hspace{-2pt} x\in Q_0\}$. A left
$A$-module $M$ is called \emph{unitary} if $M=\oplus_{x\in Q_0}\,
\varepsilon_x M$. Let ${\rm Mod\hskip 0.5pt}A$ be the category of
left unitary $A$-modules. It is well known that there exists an
equivalence from ${\rm Rep}(Q)$ to ${\rm Mod\hskip 0.5pt}A$,
sending a representation $M$ to the module $\oplus_{x\in Q_0}\,
M(x)$. For convenience, we shall make the identification
$M=\oplus_{x\in Q_0}\, M(x)$. In this way, $P_x=A{\hskip
0.5pt}\varepsilon_x$ as a module, while $A=\oplus_{x\in Q_0}P_x$
as a representation. Note that there exists a contravariant
functor ${\rm Hom}_A(-, A)$ from the category of all left
$A$-modules to that of all right $A$-modules which, however, does
not necessarily send a unitary module to a unitary one.

\medskip

\begin{Lemma} \label{ppduality} The functor $\,\Hom_A(-, A): {\rm proj}(Q)\to
{\rm proj}\hspace{0.4pt}(Q^{\hspace{0.4pt} \rm op})\,$ is a duality such that
$\Hom_A(P_x, A)\cong P_{x^{\rm o}}$ for all $x\in Q_0$.
\end{Lemma}

\medskip

Composing the dualities in Lemmas \ref{piduality} and \ref{ppduality},
we get the following result.

\medskip

\begin{Prop}\label{piequiv} The functor $\nu=D\Hom_A(-, A): {\rm proj}\hspace{0.5pt}(Q)\to {\rm inj}(Q)$,
called the {\rm Nakayama functor}, is an equivalence such that
$\nu(P_x)\cong I_x$  for all $x\in Q_0$, whose quasi-inverse is
$\nu^-={\rm Hom}_{A^{\rm op}}(D-, A^{\rm op})$. \vspace{-2pt}
\end{Prop}

{\center \section{Almost split sequences}}

For the rest of this paper, $Q$ stands for a strongly locally finite quiver and $k$ for a field. Recall that $\rep(Q)$ denotes the category of locally finite dimensional $k$-representations of $Q$, while $\rep^+\hspace{-1pt}(Q)$ and $\rep^-(Q)$ denote its full subcategories generated by the finitely presented representations and by the finitely co-presented representations, respectively.
The main objective of this section is to study the Auslander-Reiten theory in $\rep(Q)$. Our major task is to construct an almost split sequence which ends with any given indecomposable non-projective representation in $\rep^+\hspace{-1pt}(Q)$, and one which starts with any given indecomposable non-injective representation in $\rep^-(Q)$. This is a more specific version of a result by Auslander; see \cite[Theorem 6]{A}. We shall also study some properties of the Auslander-Reiten translates, and show how to relate the Auslander-Reiten theory over $Q$ to that over its subquivers.

\medskip

We need to recall some basic notions. Let $\mathcal A$ be an
additive category. Recall that an object in $\mathcal{A}$ is {\it
strongly indecomposable} if it has a local endomorphism algebra. Let
now  $f: X\to Y$ a morphism in $\mathcal{A}$. One says that $f$ is
{\it irreducible} if $f$ is neither a section nor a retraction while
every factorization $f=gh$ implies that $h$ is a section or $g$ is a
retraction. Moreover, $f$ is called {\it right minimal} if any
morphism $h: X\to X$ such that $f=fh$ is an automorphism; {\it right
almost split} if $f$ is not a retraction and every non-retraction
morphism $g: Z\to Y$ factors through $f$; and {\it minimal right
almost split} if $f$ is right minimal and right almost split. In a
dual manner, one defines $f$ to be {\it left minimal}, {\it left
almost split}, and {\it minimal left almost split}; see \cite{AuR}.
Note that a minimal left or right almost split morphism is
irreducible if and only if it is non-zero. Furthermore, a sequence
of morphisms $X\stackrel{f}{\longrightarrow} Y
\stackrel{g}{\longrightarrow} Z$ in $\mathcal A$ with $Y\ne 0$ is
called {\it almost split} if  $f$ is minimal left almost split and a
pseudo-kernel of $g$, while $g$ is minimal right almost split and a
pseudo-cokernel of $f$. Such an almost split sequence is unique for
$X$ and for $Z$ if it exists, moreover, this definition coincides
with the classical one in case $\mathcal{A}$ is abelian; see
\cite[(1.4),(1.5)]{L3}. We shall say that $\mathcal A$ is {\it right
Auslander-Reiten} if every indecomposable object in $\mathcal A$ is
the co-domain of a minimal right almost split monomorphism or the
ending term of an almost split sequence; {\it left Auslander-Reiten}
if every indecomposable object in $\mathcal A$ is the domain of a
minimal left almost split epimorphism or the starting term of an
almost split sequence; and {\it Auslander-Reiten} if it is left and
right Auslander-Reiten; compare \cite[(2.6)]{L3}. The following
result is probably well known; compare \cite[(I.3.2)]{RVDB}.

\medskip

\begin{Lemma}\label{pc-ih} Let $\mathscr{C}$ be an abelian category with a short exact sequence
\vspace{-2pt} $$\xymatrixcolsep{15pt} \xymatrix{0\ar[r]& X
\ar[r]^q & Y \ar[r]^p & Z \ar[r] & 0.}\vspace{-2pt}$$

\begin{enumerate}[$(1)$]

\item  The morphism $q$ is a minimal right almost split
monomorphism in $\mathscr{C}$ if and only if $Z$ is simple and $p$
is its projective cover.

\vspace{1pt}

\item The morphism $p$ is a minimal left almost split epimorphism
in $\mathscr{C}$ if and only if $X$ is simple and $q$ is its
injective hull.

\end{enumerate}\end{Lemma}

\noindent{\it Proof.} We shall prove only Statement (1). Suppose
first that $Z$ is simple and $p: Y\to Z$ is its projective cover.
Since $q$ is the kernel of $p$, it is right minimal. Let $f: M\to
Y$ be a non-retraction morphism. Since $Y$ is projective, $f$ is
not an epimorphism, and since $p$ is superfluous, neither is $pu$.
Since $Z$ is simple, $pf=0$, and hence $f$ factors through $q$.
That is, $q$ is minimal right almost.

Conversely, suppose that $q$ is minimal right almost split. If $Y$
is not projective, then $\mathscr{C}$ has a non-retraction
epimorphism $f: M\to Y$. Since $f$ factors through $q$, we see
that $q$ is an epimorphism and hence an isomorphism, a
contradiction. Thus $Y$ is projective. If $g: L\to Y$ is not an
epimorphism, then $g$ factors through $q$. In particular, $pg=0$,
which is not an epimorphism. This shows that $p$ is a superfluous
epimorphism and hence a projective cover of $Z$. Finally, consider
an arbitrary  morphism $u: M\to Z$ in $\mathscr{C}$. Being
abelian, $\mathscr{C}$ admits a pullback diagram \vspace{-6pt}
$$\xymatrixrowsep{16pt} \xymatrixcolsep{20pt}\xymatrix@1{
0 \ar[r] & X \ar[r]^{q\hspace{0.2pt}'} \ar@{=}[d] & \, N \ar[d]^v \ar[r]^{p\hspace{0.2pt}'} & \,M \ar[d]^u \ar[r] & \,0\\
0 \ar[r] & X \ar[r]^q & \,Y \ar[r]^p & \, Z \ar[r] &
\,0.}\vspace{2pt}$$ If $v$ is an epimorphism, then so is $u$.
Otherwise, $v=qh$ for some $h: N\to X$, and consequently,
$up\hspace{0.4pt}'=p\hspace{0.4pt}v=0$. Since $p\hspace{0.4pt}'$
is an epimorphism, we get $u=0$. This shows that $Z$ is simple.
The proof of the lemma is completed.

\medskip

\noindent {\sc Remark.} (1) In the situation as in Lemma
\ref{pc-ih}(1),  $X$ is the greatest sub-object of $Y$. One writes
$X={\rm rad}\,Y$ and calls $Z$ the {\it top} of $Y$.

\vspace{1pt}

\noindent (2) In the situation as in Lemma \ref{pc-ih}(2), $X$ is
the smallest sub-object of $Y$, which is called the {\it socle} of
$Y$ and  written as ${\rm soc}\,Y$.

\medskip

The following statment is an immediate consequence of the
preceding lemma.

\medskip

\begin{Cor}\label{AR-ab} If $\mathscr{C}$ is an abelian category, then

\begin{enumerate}[$(1)$]

\item $\mathscr{C}$ is right Auslander-Reiten if and only if every
indecomposable non-projective object is the ending term of an
almost split sequence and every indecomposable projective object
has a simple top.

\item $\mathscr{C}$ is left Auslander-Reiten if and only if every
indecomposable non-injective object is the starting term of an
almost split sequence and every indecomposable injective object
has a simple socle.

\end{enumerate}

\end{Cor}

\medskip

We now begin to study the Auslander-Reiten theory in $\rep(Q)$.
Note that, although $\rep(Q)$ is not Hom-finite in general, its
indecomposable objects are strongly indecomposable; see
\cite[(3.6)]{GR}. The following result, which is an immediate
consequence of Lemma \ref{pc-ih}, will be used frequently.

\medskip

\begin{Lemma}\label{mrasm} If $x\in Q_0$, then the inclusion $q_x: {\rm rad}\,P_x\to P_x$
is a minimal right almost split monomorphism, and the projection $p_x: I_x\to I_x/\,{\rm soc}\,I_x$
is a minimal left almost split epimorphism in $\rep(Q)$.

\end{Lemma}

\medskip

The following construction is analogous to the classical one for modules over an artin algebra; see, for example, \cite{ARS}.

\medskip

\begin{Defn}\label{DTr} Let $M$ be a representation in $\rep(Q)$.

\vspace{-3pt}

\begin{enumerate}[$(1)$]

\item If $M$ has a minimal projective
resolution $\xymatrixcolsep{14 pt}\xymatrix{0\ar[r]& P_1 \ar[r]^f&
P_0 \ar[r] & M \ar[r]& 0}\vspace{-4pt}$ with $P_1, P_0\in {\rm proj}\hskip 0.5pt(Q)$, then $\DTr M$ denotes the kernel of $\xymatrixcolsep{14 pt}\xymatrix{\nu(f): \nu(P_1)\ar[r] & \nu(P_0)}\hspace{-2pt}.$

\vspace{-2pt}

\item If $M$ has a minimal injective
co-resolution $\xymatrixcolsep{14pt}\xymatrix{0\ar[r]& M \ar[r] &
I_0 \ar[r]^g & I_1 \ar[r]& 0}\vspace{-4pt}$ with $I_0, I_1\in {\rm inj}\hskip
0.5pt (Q)$, then $\TrD M$ denotes the co-kernel of $\xymatrixcolsep{14 pt}\xymatrix{\nu^-(g): \nu^-(I_0)\ar[r] & \nu^-(I_1)}\hspace{-2pt}.$

\end{enumerate} \end{Defn}

\medskip

\noindent{\sc Remark.} (1) $\DTr M$ is defined only up to isomorphism and only for $M\in \rep^+\hspace{-1pt}(Q)$, in such a way that
$\DTr M=0$ if and only if $M\in {\rm proj}\hspace{0.4pt}(Q)$.

\vspace{1pt}

\noindent (2) $\TrD M$ is defined only up to isomorphism and only for $M \in \rep^-(Q)$, in such a way that
$\TrD M=0$ if and only if $M\in {\rm inj}\hspace{0.4pt}(Q)$.

\medskip

The following lemma and its dual play an an important role in the construction of almost split sequences.

\medskip

\begin{Lemma} \label{AR-translate} Let $M$ be an indecomposable object in $\rep^+\hspace{-1pt}(Q)$ with a minimal projective
resolution

\vspace{-16pt} $$\xymatrixcolsep{18pt}\xymatrix{0\ar[r]& P_1
\ar[r]^f& P_0 \ar[r] & M \ar[r]& 0.}$$ If $M$ is not projective,
then $\DTr M\cong D{\rm Ext}^1_A(M, A)\vspace{1pt}$, which is an
indecomposable non-injective object in $\rep^-(Q)$ with a minimal
injective co-resolution

\vspace{-8pt}

$$\xymatrixcolsep{18pt}\xymatrix{0 \ar[r] & \DTr M \ar[r] & \nu(P_1) \ar[r]^{\nu(f)}&
\nu(P_0) \ar[r] & 0.}\vspace{0pt}$$
\end{Lemma}

\noindent{\it Proof.} Suppose that $M\vspace{-2pt}$ is not projective.
Since $\Rep(Q)$ is hereditary, we have $\Hom_A(M,A) = 0$. Applying $\Hom_A(-, A)$ to the minimal projective resolution stated in the lemma, we get a short exact
sequence of right $A$-modules as follows:

\vspace{-8pt}

$$\xymatrixcolsep{15 pt}\xymatrix{0 \ar[r] & \Hom_A(P_0,A) \ar[r]^{f^*} & \Hom_A(P_1,A) \ar[r] & \Ext^1_A(M,A) \ar[r] & 0,}$$

\vspace{-1pt}

\noindent where $\Ext^1_A(M,A)$ is unitary since $\Hom_A(P_0,A)$
and $\Hom_A(P_1,A)$ are unitary. Applying the duality $D:
\rep(Q^{\hspace{0.4pt} \rm op})\to \rep(Q)$ stated in Lemma
\ref{piduality}, we obtain a short exact sequence

\vspace{-12pt}

$$\xymatrixcolsep{20pt}\xymatrix{\eta: \qquad 0 \ar[r] & D\Ext^1_A(M,A)  \ar[r] & \nu(P_1) \ar[r]^-{\nu(f)}& \nu(P_0) \ar[r]& 0 \vspace{0pt}}$$
in $\rep(Q)$. By definition, $\DTr M \cong D\Ext^1_A(M,A)\in
\rep^-(Q)$. Furthermore, since ${\rm Im}(f)\subseteq {\rm rad}\hskip
0.5pt P_0$, we see that ${\rm Im}(f^*)$ is contained in the radical of
$\Hom_A(P_1,A)$, and hence the kernel of $\nu(f)$ contains the socle of $\nu(P_1)$. That is, $\eta$ is a minimal injective co-resolution of
$\DTr M$. In particular, $\DTr M$ is not injective. Finally, since
$\nu$ is an equivalence and $M$ is indecomposable,  $\DTr M$ is indecomposable. The proof of the lemma is
completed.

\medskip

As a consequence of Lemma \ref{AR-translate} and its dual, we have
the following result.

\medskip

\begin{Cor} \label{DTr-TrD}If $M, N$ are indecomposable objects in $\rep(Q)$,
then $N\cong \DTr M$ if and only if $M\cong \TrD N.$
\end{Cor}

\medskip

The following consequence of Lemma \ref{AR-translate} will be needed later; compare
\cite[(4.2)]{AuR}.

\medskip

\begin{Cor} \label{preservemono} Let $M, N$ be indecomposable non-projective objects in
$\rep^+\hspace{-1pt}(Q)$. If $\rep^+\hspace{-1pt}(Q)$ has a
monomorphism $f: M \to N$, then $\rep^-(Q)$ has a monomorphism $g:
\DTr M\to \DTr N$.
\end{Cor}

\noindent{\it Proof.} Since $\rep(Q)$ is hereditary,
$D\Ext^1_A(-,A): \rep^+\hspace{-1pt}(Q)\to \rep^-(Q)$ is a left exact functor.
If $f: M \to N$ is a  monomorphism in $\rep^+\hspace{-1pt}(Q)$, then

\vspace{-6pt} $$g=D{\rm Ext}_A^1(f, A): D{\rm Ext}_A^1(M, A)\to
D{\rm Ext}_A^1(N, A)$$ is a monomorphism in $\rep^-(Q)$. By Lemma
\ref{AR-translate}, $g$ is a monomorphism from $\DTr M$ to $\DTr
N$. The proof of the corollary is completed.

\medskip

We are ready to have the existence theorem for almost split sequences.

\medskip

\begin{Theo} \label{AR-sequence} Let $Q$ be a strongly locally finite quiver, and let $M$ be an indecomposable representation in $\rep(Q)$.

\begin{enumerate}[$(1)$]

\item If $M\in \rep^+\hspace{-1pt}(Q)$ is not projective, then
$\rep(Q)$ has an almost split sequence
$\xymatrixcolsep{15pt}\xymatrix{0\ar[r] & \DTr M \ar[r] & N \ar[r] &
M \ar[r] & 0,}$ where $\DTr M\in \rep^-(Q).$

\vspace{1pt}

\item If $M\in \rep^-(Q)$ is not injective, then
$\rep(Q)$ has  an almost split sequence
$\xymatrixcolsep{15pt}\xymatrix{0\ar[r] & M \ar[r] & N \ar[r] & \TrD
M \ar[r] & 0,}$ where $\TrD M\in \rep^+\hspace{-1pt}(Q)$.

\end{enumerate} \end{Theo}

\noindent{\it Proof.} We only prove Statement (1). Assume that $M$ is finitely presented and not projective. By Lemma \ref{AR-translate}, $\DTr
M$ is finitely co-presented, indecomposable, and not injective. Let $L\in \rep(Q)$. By Lemma \ref{Prop1.13_shiping}, $\Ext^1_{A}(M,L)$ and $\Hom_{A}(L, \DTr M)$ are of finite $k$-dimension.
We claim, for $P\in {\rm proj}(Q)\vspace{1pt}$, that there
exists a $k$-linear isomorphism, which is
natural in $P$ and $L$, as follows$\,:$
\vspace{-1.5pt}
$$\psi_{_{L,P}}: \Hom_A(P,L) \to D \Hom_A(L,\nu P).\vspace{-0.5pt}$$ Indeed, we may assume with no loss of generality that $P=P_x$ for some $x\in Q_0$. By Proposition \ref{prop1.3_shiping}, we have the following $k$-isomorphisms:
$$\Hom_A(P_x, L)\cong L(x) \cong D \Hom_k(L(x), k) \cong D \Hom_A(L, I_x),$$
each of which is natural in $P_x$ and $L$. This establishes our claim.

Let $\xymatrixcolsep{13 pt}\xymatrix{0 \ar[r] & P_1 \ar[r] & P_0 \ar[r]
& M \ar[r] & 0}$ be a minimal projective resolution of $M$, where $P_0, P_1\in {\rm proj}(Q)$. By Lemma \ref{AR-translate}, $\DTr M$ has
a minimal injective co-resolution
$\xymatrixcolsep{13 pt}\xymatrix{0
\ar[r] & \DTr M \ar[r] & \nu(P_1) \ar[r] & \nu(P_0) \ar[r] & 0,}$ where
$\nu(P_1), \nu(P_0)\in {\rm inj}(Q)$. Applying $\Hom_{A}(-,L)$ and $D\Hom_A(L, -)$, we get a commutative diagram with exact rows:
$$\xymatrixcolsep{13pt}\xymatrix{
\Hom_{A}(P_0, L) \ar[r]\ar[d]^{\psi_{_{P_0, \hskip 0.5pt L}}} & \Hom_{A}(P_1,L) \ar[r] \ar[d]^{\psi_{_{P_1, \hskip 0.5pt L}}} & \Ext^1_{A}(M,L) \ar[r] \ar[d]^{\phi_{_L}} & 0\\
D\Hom_{A}(L,\nu P_0) \ar[r] & D\Hom_{A}(L, \nu P_1) \ar[r]&
D\Hom_{A}(L,\DTr M) \ar[r]& \,0,}$$ where $\psi_{_{P_0, \hskip 0.5pt
L}}\vspace{1pt}$ and $\psi_{_{P_1, \hskip 0.5pt L}}$ are natural
isomorphisms. Thus there exists an isomorphism $\phi_{_L}:
\Ext^1_{A}(M,L) \to D\Hom_{A}(L, \DTr M),$ which is natural in $L$.
Since $\End(\DTr M)$ is finite dimensional, there exists a non-zero
$k$-linear map $\theta: {\End}(\DTr M)\to k,$ which vanishes on
${\rm rad}({\End}(\DTr M))$. Consider the corresponding non-zero
element

\vspace{-8pt}

$$\xymatrixcolsep{20pt}\xymatrix{\eta=\phi_{_{{\rm DTr} M}}^{-1}(\theta): \quad 0 \ar[r] & \DTr M \ar[r]^-f & N \ar[r]^g & M \ar[r] & 0}$$ in $\Ext^1(M, \DTr M).$
Let $u: \DTr M \to L$ be a non-section morphism $\rep(Q)$. For any
$v: L\to \DTr M$, since $vu\in {\rm rad}(\End(\DTr M))$, we have $\theta(vu)=0$. This shows that $D\Hom(u, \DTr
M)(\theta)=0$. In view of the commutative diagram

\vspace{-8pt}

$$\xymatrixcolsep{20 pt}\xymatrix{\Ext^1(M, \DTr M)\ar[rr]^{\Ext^1(M, u)}\ar[d]^{\phi_{_{\DTr M}}} && \Ext^1(M,
L) \ar[d]^{\phi_{_L}}\\
D\End(\DTr M) \;\ar[rr]^-{D\Hom(u, \DTr M)} &&\; D\Hom(L, \DTr M),
}$$ we get $\Ext^1(M, u)(\eta)=0$, that is, $u$ factors through $f$.
Thus $\eta$ is an almost split sequence in $\rep(Q)$; see \cite[(2.14)]{AuR}. The
proof of the theorem is completed.

\medskip

\noindent{\sc Remark.} It is shown in \cite{Paq} that every almost split sequence in $\rep(Q)$ is of the form stated in Theorem \ref{AR-sequence}.

\medskip

The following result is a consequence of Theorem \ref{AR-sequence} and Proposition \ref{cat-fin-pres}.

\medskip

\begin{Cor}\label{cor_ARS_exist} Let $M$ be an indecomposable representation in $\rep^b(Q)$.

\begin{enumerate}[$(1)$]

\item If $M$ is not projective, then $\rep(Q)$ has an almost split sequence
ending with $M$, which is also an almost split sequence in
$\rep^-\hspace{-1pt}(\hspace{-1pt}Q\hspace{-1pt})\hspace{0pt}$.

\vspace{1pt}

\item If $M$ is not injective, then $\rep(Q)$ has an almost split sequence starting with $M$, which is also an almost split sequence in $\rep^+\hspace{-1pt}(Q).$

\end{enumerate} \end{Cor}

\medskip

Next, we shall study the Auslander-Reiten translates. To this end, the following easy result is useful.

\medskip

\begin{Lemma}\label{resol-supp} Let $M$ be an indecomposable representation in $\rep^+\hspace{-1pt}(Q)$ with
a minimal projective resolution \vspace{-8pt}
$$\qquad \xymatrixrowsep{15pt}
\xymatrixcolsep{20pt}\xymatrix{0 \ar[r] &
\oplus_{j=1}^s P_{y_j} \ar[r]&
\oplus_{i=1}^r P_{x_i} \ar[r]^-f & M \ar[r]& 0.} \vspace{-2pt}$$  If $x\in Q_0$ is not in ${\rm supp}\hspace{0.4pt}M$, then $x=y_j$ for some $1\le j\le s$ if and only if $x$ is an immediate successor of some vertex in ${\rm supp}\hspace{0.4pt}M$.
\end{Lemma}

\noindent{\it Proof.} Let $N$ denote the kernel of $f$. Then $y_1, \ldots, y_s$ are the vertices in the top-support of $N$. Fix a vertex $x$ not lying in ${\rm supp}\hspace{0.4pt}M$. Suppose first that $Q$ has an arrow $\alpha: y\to x$ with $y\in {\rm supp}\hspace{0.4pt}M$. Since $f$ is surjective, there exists a path $p$ in $Q$ from some $x_i$ to $y$ such that $f(p)\ne 0$. Since $x\not\in {\rm supp}\hspace{0.4pt}M$, we have $f(\alpha p)=0$. Thus $\alpha p$ lies in $N$ but not in its radical, and hence, $x$ lies in the top-support of $N$.

Suppose conversely that $x$ lies in the top-support of $N$. Then there exists some $\rho=\lambda_1 p_1+\cdots +\lambda_t p_t$, where $\lambda_1 \ldots, \lambda_t\in k$, and $p_1, \ldots, p_t$ are paths in $Q$ from some of $x_1, \ldots, x_r$ to $x$, which lies in $N$ but not in its radical. Since $x\not\in {\rm supp}\hspace{0.4pt}M,$ we may write $p_i=\alpha_iq_i$, where $\alpha_i$ is an arrow ending in $x$ and $q_i$ is a path in $Q$, for $i=1, \ldots, t$. Since $\rho\not\in  {\rm rad}\hspace{0.4pt}N$, there exists some $1\le i_0\le t$ such that $f(q_{i_0})\ne 0$. In particular, $e(q_{i_0})\in {\rm supp}\hspace{0.5pt}M$. The proof of the lemma is completed.

\medskip

The following result is an immediate consequence of Lemmas \ref{AR-translate} and \ref{resol-supp}.

\medskip

\begin{Cor}\label{DTr-supp} Let $M$ be an indecomposable representation in $\rep^+\hspace{-1pt}(Q)$. If $x$ is a vertex in $Q$ not lying in ${\rm supp}\hspace{0.4pt}M$, then $x$ lies in the socle-support of $\,\DTr M$ if and only if $x$ is an immediate successor of some vertex in ${\rm supp}{\hskip 0.5pt}M$.
\end{Cor}

\medskip

Let $M$ be an indecomposable representation in $\rep(Q)$. By
convention, we write
$\DTr^0\hspace{-1.5pt}M=M=\TrD^0\hspace{-1.5pt}M$. If $M\in
\rep^+(Q)$, then $\DTr M\in \rep^-(Q)$; and if, moreover, $\DTr
M\in \rep^+(Q)$, then $\DTr^2\hspace{-1pt}M$ is defined and lies
in $\rep^-(Q)$. In general, if $n>0$ is such that
$\DTr^{n-1}\hspace{-1pt}M$ is defined and lies in $\rep^+(Q)$,
then $\DTr^n\hspace{-1pt}M$ is defined and lies in $\rep^-(Q)$. If
$\DTr^n\hspace{-1pt}M$ is defined and non-zero for some $n>0$,
then it follows from  Proposition \ref{cat-fin-pres} and Lemma
\ref{AR-translate} that $\DTr^i\hspace{-1pt}M$ is indecomposable
for $0\le i\le n$, and finite dimensional for $0<i<n$. We shall
say that $M$ is {\it $\DTr$-stable} if $\DTr^n\hspace{-1pt}M$ is
defined and non-zero for all $n\ge 0$, or equivalently,
$\DTr^n\hspace{-1pt}M$ is indecomposable of finite dimension for
all $n> 0$.

\medskip

\begin{Lemma} \label{stableorbit2}
Let $M\in \rep^+\hspace{-1pt}(Q)$ be indecomposable, and let $w$
be an infinite acyclic walk in $Q$ which starts with an arrow and
intersects ${\rm supp}\hspace{0.4pt}M$ only at $s(w)$. Then $M$ is
$\DTr$-stable or $\DTr^n\hspace{-1pt}M$ is infinite dimensional
for some $n\ge 0$. Furthermore,

\vspace{-1.5pt}

\begin{enumerate}[$(1)$]

\item  if all but finitely many edges in $w$ are inverses of
arrows, then $\DTr^n\hspace{-1pt}M$ is infinite dimensional for
some $n\ge 0;$

\vspace{1pt}

\item if $\,\DTr^m\hspace{-1pt}M$ with $m>0$ is defined, then its
support contains some vertex lying in $w$ but different from
$s(w)$.

\end{enumerate} \end{Lemma}

\noindent{\it Proof.} Write $w=\cdots w_n \cdots w_2w_1$, where
the $w_i$ are edges. Put $a_i=s(w_i)$ for $i\ge 1$. Set $s_0=1$.
Then $w_{s_0}$ is an arrow and $s_0$ is maximal such that
$a_{s_0}\in {\rm supp}\hspace{0.4pt}M$. Let $r\ge 0$ be an integer
such that $\DTr^r\hspace{-1pt}M\in \rep^+\hspace{-1pt}(Q)$ and
there exist integers $s_0< \cdots < s_r$ satisfying the following
property: $w_{s_i}$ is an arrow and $s_i$ is maximal for which
$a_{s_i}$ is in the support of $\DTr^i\hspace{-1.5pt}M$, for $i=0,
\ldots, r$. Since $w_{s_r}$ is an arrow $a_{s_r}\to a_{s_r+1}$ and
$a_{s_r+1}$ is not in the support of $\DTr^{r}\hspace{-1.5pt}M$,
by Corollary \ref{DTr-supp}, $a_{s_r+1}$ is in the support of
$\DTr^{r+1}\hspace{-1.5pt}M$. If
$\DTr^{r+1}\hspace{-1.5pt}M\not\in \rep^+\hspace{-1pt}(Q)$, then
$\DTr^{r+1}\hspace{-1.5pt}M$ is infinite dimensional and
$\DTr^i\hspace{-1pt}M$ is not defined for every $i>r+1$. In this
case, the lemma is proved and we stop the process. Otherwise,
$\DTr^{r+1}M$ is non-zero of finite dimension by Proposition
\ref{cat-fin-pres}. Therefore, there exists a maximal integer
$s_{r+1}>s_r$ such that $a_{s_{r+1}}$ is in the support of
$\DTr^{r+1}M$. Suppose that $w_{s_{r+1}}$ is the inverse of an
arrow $a_{s_{r+1}}\leftarrow a_{s_{r+1}+1}$. Since $a_{s_{r+1}+1}$
is not in the support of $\DTr^{r+1}M$, applying the dual of
Corollary \ref{DTr-supp} to $\DTr^{r+1}M$, we see that
$a_{s_{r+1}+1}$ is in the support of $\DTr^r\hspace{-1pt}M$,
contrary to the maximality of $s_r$. Therefore, $w_{s_{r+1}}$ is
an arrow $a_{s_{r+1}}\to a_{s_{r+1}+1}$. If this process never
stops, then we get an infinite increasing sequence  of integers
$$s_0<s_1<\cdots
<s_i<\cdots $$ satisfying the above-stated property. In
particular, $M$ is $\DTr$-stable, and $a_{s_i}$ lies in the
support of $\DTr^i\hspace{-1pt}M$ for all $i\ge 0$. Moreover,
since the $w_{s_i}$ are arrows, the hypothesis stated in Statement
(1) does not occur. The proof of the lemma is completed.

\medskip

\noindent{\sc Remark.} If $Q$ has no left infinite path, then $\rep^-(Q)=\rep^b(Q).$ In this case, Lemma \ref{stableorbit2}
provides a simple combinatorial condition for an
indecomposable representation to be $\DTr$-stable.

\medskip

Dually, if $M\in \rep^-(Q)$ is indecomposable, then $\TrD M\in
\rep^+(Q)$. If $n>0$ is such that $\TrD^{n-1}\hspace{-1pt}M$ is
defined and lies in $\rep^-(Q)$, then $\TrD^n\hspace{-1pt}M$ is
defined and lies in $\rep^+(Q)$. If $\TrD^n\hspace{-1pt}M$ is
defined and non-zero for some $n>0$, then it follows from
Proposition \ref{cat-fin-pres} and the dual of Lemma
\ref{AR-translate} that $\TrD^i\hspace{-1pt}M$ is indecomposable
for $0\le i\le n$, and finite dimensional for $0<i<n$. We shall
say that $M$ is {\it $\TrD$-stable} if $\TrD^n\hspace{-1pt}M$ is
defined and non-zero for all $n\ge 0$, or equivalently,
$\TrD^n\hspace{-1pt}M$ is indecomposable of finite dimension for
all $n> 0$. The following result is a dual statement of Lemma
\ref{stableorbit2}.

\medskip

\begin{Lemma} \label{stableorbit}
Let $M\in \rep^-(Q)$ be indecomposable, and let $w$ be an infinite
acyclic walk in $Q$ which ends with an arrow and intersects ${\rm
supp}\hspace{0.4pt}M$ only at $e(w)$. Then $M$ is $\TrD$-stable or
$\TrD^{n}M$ is infinite dimensional for some $n\ge 0$. Moreover,

\vspace{-0.5pt}

\begin{enumerate}[$(1)$]

\item if all but finitely many edges in $w$ are inverses of
arrows, then $\TrD^{n}M$ is infinite dimensional for some $n\ge
0;$

\vspace{1pt}

\item if $\,\TrD^mM$ with $m>0$ is defined, then its support
contains some vertex lying in $w$ but different from $e(w)$.
\end{enumerate}\end{Lemma}

\medskip

The preceding results yield some very useful consequences.

\medskip

\begin{Prop}\label{orbits} Suppose that $Q$ is infinite and connected.

\vspace{-1pt}

\begin{enumerate}[$(1)$]

\item For any $x, y\in Q_0$, there exists no integer $m \ge 0$
such that $\DTr^m \hspace{-1pt} I_y \cong P_x$.

\vspace{1pt}

\item If $M, N\in \rep(Q)$ are indecomposable such that $M\cong
\DTr^n \hspace{-1pt}N$ for some $n \ge 0$, then $\,{\rm supp}
\hskip 0.5 pt M = {\rm supp} \hskip 0.5pt N$ if and only if $n=0$.

\end{enumerate} \end{Prop}

\noindent{\it Proof.} (1) Let $x, y\in
Q_0$ be such that $P_x\cong \DTr^m\hspace{-1pt}I_y$ for some $m\ge 0$. If $m=0$,
then it is easy to see that $Q$ consists of a single path from $x$ to
$y$, a contradiction. Thus $m>0$. Since $P_x\in \rep^+\hspace{-1pt}(Q)$,
by Proposition \ref{cat-fin-pres}, $\DTr^i\hspace{-1pt}I_y$ is finite dimensional,
for $i=1, \ldots, m$. On the other hand, since $I_y=\TrD^mP_x$ by Corollary \ref{DTr-TrD}, $I_y$ is finite dimensional by Proposition \ref{cat-fin-pres}. Since $\DTr^{m+1}\hspace{-1pt}I_y=0$, we see that $I_y$ is not $\DTr$-stable and
$\DTr^i\hspace{-1pt}I_y$ is finite dimensional for all $i\ge 0$.
Since $Q$ is connected and infinite, applying K\"{o}nig's lemma to the complement of ${\rm supp}{\hskip
0.5pt}I_y$, we get a right infinite acyclic walk $w$ which
intersects ${\rm supp}{\hskip 0.5pt}I_y$ only at $s(w)$. Since ${\rm
supp}\hskip 0.5pt I_y$ is successor-closed, $w$ starts with an arrow.
By Lemma \ref{stableorbit2}, $I_y$ is $\DTr$-stable or $\DTr^r\hspace{-1pt}I_y$ is
infinite dimensional for some $r\ge 0$, a contradiction.

\vspace{1pt}

(2) Let $M, N\in \rep(Q)$ be indecomposable such that $M\cong \DTr^n \hspace{-1pt}N$ with $n \ge 0$. By Corollary \ref{DTr-TrD},
$N\cong \TrD^n \hspace{-1pt}M.$ Suppose that $n>0$ and that ${\rm supp}M={\rm supp}N=\Sa$. Then $N\in \rep^+\hspace{-1pt}(Q)$ and $M\in
\rep^-(Q)$. Moreover, by Lemma \ref{Lemma1.6_shiping}, $\Sa$ is top-finite and socle-finite,
and hence finite. Applying K\"{o}nig's lemma to the complement of $\Sa$, we get
a left infinite acyclic walk $w$ which intersects $\Sa$ only at
$e(w)$. If $w$ ends with an arrow, then we may apply Lemma \ref{stableorbit}(2) to $M$ to see that
${\rm supp}\hspace{0.4pt} N$, that is $\Sa$, contains some vertex lying in $w$ but different from
$e(w)$, a contradiction.
If $w$ ends with the inverse of an arrow, then $w^{-1}$ is an infinite
acyclic walk which starts with an arrow and intersects
${\rm supp}\hspace{0.4pt}N$ only at $s(w^{-1})$. It follows from Lemma
\ref{stableorbit2}(2) that ${\rm supp}\hspace{0.4pt}M$, that is $\Sa$, contains some
vertex lying in $w^{-1}$ but different from $s(w^{-1})$, a contradiction again.
The proof of the proposition is completed.

\medskip

\noindent{\sc Remark.} Proposition \ref{orbits}(1) is well known in
the finite non-Dynkin case. If $Q$ is infinite without left infinite
paths, Reiten and Van der Bergh proved this by using a highly
indirect argument to treat the infinite Dynkin case; see
\cite{RVDB}.

\medskip

Let $\it\Sigma$ be a subquiver of $Q$. If $M$ is a representation of $Q$ supported by $\Sa$, for the sake of convenience, we shall regard $M$ as a representation of $\Sa$ whenever no risk of confusion is possible. In particular, if $N$ is a representation of $Q$, then its restriction $N_{_{\it\Sigma}}$ will be regarded as a representation of $\it\Sigma$. On the other hand, every representation $M$ of $\Sa$ can be extended trivially to a representation of $Q$ which, by abuse of notation, is denoted again by $M$. In this way, we shall identify $\rep(\it\Sigma)$ with the full subcategory of $\rep(Q)$ generated by the representations supported by $\Sa$. One of our techniques in our later investigation is to  relate the almost split sequences and the irreducible morphisms in $\rep(Q)$ to those in $\rep(\Sa)$.

\medskip

\begin{Lemma}\label{rep+restrict} Let $\Sa$ be a convex subquiver of $Q$, and let $M$ be an object in $\rep^+\hspace{-1pt}(Q)$.

\vspace{-1.5pt}

\begin{enumerate}[$(1)$]

\item If $\Sa$ is predecessor-closed in $Q$, then $M_{_{\hspace{-1pt}\it \Sigma}} \in \rep^+(\Sa)$.

\vspace{1pt}

\item If $\Sa$ contains the trivial and the immediate successors of the vertices in ${\rm supp}\hspace{0.4pt}M$, then $M\in \rep^+\hspace{-1pt}(\Sa)$.

\end{enumerate}
\end{Lemma}

\noindent{\it Proof.} For $x\in Q_0$, let $P_x'$ denote the
restriction of $P_x$ to $\Sa$. Since $\Sa$ is convex, for $x\in
\Sa$, it is easy to see that $P_x'$ is isomorphic to the
indecomposable projective representation in $\rep(\Sa)$ at $x$. Now,
$M$ has a minimal projective resolution \vspace{-1pt} $$\eta: \quad
\xymatrixrowsep{15pt} \xymatrixcolsep{20pt}\xymatrix{0 \ar[r] &
\oplus_{j=1}^s P_{y_j} \ar[r]& \oplus_{i=1}^r P_{x_i} \ar[r] & M
\ar[r]& 0.} \vspace{-1pt} $$ Restricting $\eta$ to $\Sa$, we get a
short exact sequence in $\rep(\Sa)$ as follows$\,:$ \vspace{-1pt}
$$\eta_{_{\it\Sigma}}: \quad \xymatrixrowsep{15pt}
\xymatrixcolsep{20pt}\xymatrix{ 0 \ar[r] &
\oplus_{j=1}^s P_{y_j}' \ar[r]
& \oplus_{i=1}^r P_{x_i}' \ar[r] & M_{_{\hspace{-1pt} \it\Sigma}} \ar[r]& 0.}\vspace{-1pt}$$

Suppose that $Q$ is predecessor-closed in $Q$. Then $P'_x=0$, for $x\not\in \Sa$. This implies that $P'_x\in {\rm proj}(\Sa)$, for all $x\in Q_0$. In particular, $\eta_{_{\it\Sigma}}$ is a minimal projective resolution of $M_{_{\hspace{-1pt} \it\Sigma}}$ in $\rep(\Sa)$. That is, $M_{_{\hspace{-1pt}\it \Sigma}} \in \rep^+(\Sa)$.

Suppose next that $\Sa$ contains the trivial and the immediate successors of the vertices in ${\rm supp}\hspace{0.4pt}M$. Then $M_{_{\hspace{-1pt}\it \Sigma}}=M$, and by Corollary \ref{DTr-supp}, the $x_i$ and the $y_j$ all lie in $\Sa$. Hence, the $P'_{x_i}$ and the $P'_{y_j}\vspace{1pt}$ all lie in ${\rm proj}(\Sa)$. As a consequence, $\eta_{_{\it\Sigma}}$ is a minimal projective resolution of $M$ in $\rep(\Sa)$. That is, $M \in \rep^+(\Sa)$.
The proof of  the lemma is completed.

\medskip

\begin{Prop}\label{ARS-extension} Let $\Sa$ be a convex subquiver of $Q$, and let $N$ be an indecomposable object in $\rep^+\hspace{-1pt}(Q)$ such that
the predecessors of the trivial and the
immediate successors of the vertices in ${\rm supp}\, N$ are all contained in $\Sa$.

\vspace{-1pt}

\begin{enumerate}[$(1)$]

\item The almost split sequence in $\rep(\Sa)$ ending with $N$ is almost split in $\rep(Q)$.

\vspace{2pt}

\item Every irreducible morphism in $\rep(\Sa)$ ending in $N$ is irreducible in $\rep(Q)$.

\end{enumerate} \end{Prop}

\noindent{\it Proof.} First of all, $N$ can be considered to be a representation of $\Sa$, which is finitely presented by Lemma \ref{rep+restrict}(2). Moreover, since $\Sa$ contains the immediate successors
of the vertices in ${\rm supp}\hspace{0.4pt}N$, we see that $N$ is projective in $\rep(\Sa)$ if and only if it is projective in $\rep(Q)$.

\vspace{1pt}

(1) Let $\eta: \xymatrixcolsep{15pt}\xymatrix{0\ar[r]& L \ar[r]&
M\ar[r]& N \ar[r]& 0}$ be an almost split sequence in $\rep(\Sa)$. In particular, $N$ is not projective in $\rep(\Sa)$, and hence
not projective in $\rep(Q)$. Then,
$\rep(Q)$ has an almost split sequence $\zeta:
\xymatrixcolsep{15pt}\xymatrix{0\ar[r]& \DTr N \ar[r]& E \ar[r]& N
\ar[r]& 0}$. By Corollary \ref{DTr-supp}, the vertices
in the socle-support of $\DTr N$ are trivial or immediate successors
of the vertices in ${\rm supp}\,N$. Since the socle of $\DTr N$ is essential, the vertices in
the support of $\DTr N$ are all predecessors of the vertices in the socle-support of $\DTr N$, which lie in
$\Sa$ by the hypothesis stated in the proposition. Therefore, $\zeta$ lies entirely in $\rep(\Sa)$. Then $\zeta$ is an almost split sequence in $\rep(\Sa)$, and hence, it is isomorphic to $\eta$. In other words, $\eta$ is almost
split in $\rep(Q)\vspace{1pt}$.

(2) Let $f: M\to N$ be an irreducible morphism in $\rep(\Sa)$. Since $N\in \rep^+(\Sa)$, by Lemma \ref{mrasm} and Theorem \ref{AR-sequence}(1),
$\rep(\Sa)$ has a minimal right almost split morphism $g: L\to N$. Then $f=gs$ for some section $s: M\to L$. If $N$ is not
projective in $\rep(\Sa)$, then $g$ is minimal right almost split in
$\rep(Q)$ by Statement (1). Otherwise, $N=P_x$ for some $x\in \Sa_0$,
and $g$ is the inclusion map ${\rm rad}\hspace{0.5pt}P_x\to P_x$. In any case, $f$ is irreducible in $\rep(Q)$. The
proof of the proposition is completed.

\medskip

Conversely, we have the following result.

\medskip

\begin{Prop} \label{prop2.9} Let $\Sa$ be a convex subquiver of $Q$, and let $N$ be an indecomposable object in $\rep^+\hspace{-1pt}(Q)$ such that
the trivial and the immediate successors of the vertices
in ${\rm supp}\, N$ are all contained in $\Sa$.

\begin{enumerate}[$(1)$]

\item If $\xymatrixcolsep{13pt}\xymatrix{0\ar[r]& L \ar[r]& M \ar[r]& N \ar[r]& 0}\vspace{-1pt}$ is an almost split sequence
in $\rep(Q)$, then its restriction
$\xymatrixcolsep{13pt}\xymatrix{0\ar[r]& L_{_{\hspace{-0.8pt}\it\Sigma}} \ar[r]&
M_{_{\hspace{-1pt}\it \Sigma}} \ar[r]& N \ar[r]& 0}$ is an almost split sequence
in $\rep(\it\Sigma)$.

\vspace{1pt}

\item If $f: M\to N$ is an irreducible morphism in $\rep(Q)$, then $f_{_{\it\Sigma}}: M_{_{\hspace{-1pt}\it \Sigma}}\to N$ is an irreducible morphism in $\rep(\it\Sigma)$.

\vspace{-2pt}

\end{enumerate}\end{Prop}

\noindent {\it Proof.} (1) Let $\xi: \;
\xymatrixcolsep{15pt}\xymatrix{0\ar[r]& L \ar[r]& M \ar[r]^g& N
\ar[r]& 0}$ be an almost split sequence in $\rep(Q)$. Restricting
$\xi$ to $\Sa$, we get a short exact sequence  in $\rep(\Sa)$ as
follows$\,:$ \vspace{-4pt}
$$\xi_{_{\hspace{-1pt}\it\Sigma}}: \;
\xymatrixrowsep{15pt}
\xymatrixcolsep{18pt} \xymatrix{ 0 \ar[r] & L_{_{\hspace{-0.8pt}\it\Sigma}} \ar[r]
& M_{_{\hspace{-1pt}\it \Sigma}} \ar[r]^{g_{_{\hspace{-1pt}\it\Sigma}}} &  N \ar[r] & 0.} \vspace{-3pt}$$

If $g_{_{\hspace{-1pt}\it\Sigma}}$ is a retraction, then  $\rep(\hspace{-1.5pt}\it\Sigma)$ has a morphism $h': N\to
M_{_{\hspace{-1pt}\it \Sigma}}\vspace{1pt}$ such that
$g_{_{\hspace{-1pt}\it\Sigma}} h'=\id_N$.
Since $\Sa$
contains the immediate successors of the vertices in ${\rm supp}\,N$, we can extend $h'$ to a
morphism $h: N\to M$ in $\rep(Q)$ such that $gh=\id_N$, a
contradiction. If $u : X \to N$ is
a non-retraction morphism in $\rep(\Sa)$, then it is
not a retraction in $\rep(Q)$. Thus
$u=gv$, for some morphism $v: X\to M$ in $\rep(Q)$. Restricting the equation to $\Sa$ yields $u=g_{_{\hspace{-1pt}\it\Sigma}}v_{_{\hspace{-1pt}\it\Sigma}}$.
This shows that $g_{_{\hspace{-1pt}\it\Sigma}}$ is right almost split in $\rep(\Sa)$.

For $x\in Q_0$, let $P_x'$ and $I_x'$ denote the restrictions
of $P_x$ and $I_x$ to $\Sa$, respectively.  If $x\in \Sa$ then, since $\Sa$ is convex, $P_x'$ and $I_x'$ are isomorphic to the
indecomposable projective and injective representations in $\rep(\Sa)$ at $x$, respectively. Now, $N$ has a minimal projective resolution
\vspace{-2pt} $$\eta: \quad \xymatrixrowsep{15pt}
\xymatrixcolsep{20pt}\xymatrix{0 \ar[r] &
\oplus_{j=1}^s P_{y_j} \ar[r]^w&
\oplus_{i=1}^r P_{x_i} \ar[r] & N \ar[r]& 0}$$ in $\rep(Q).$ By Corollary \ref{DTr-supp}, the $x_i$ and the $y_j$ all lie in $\Sa$.
Thus, restricting $\eta$ to $\Sa$, we get a minimal projective resolution
\vspace{-4pt}
$$\eta_{_{\hspace{-0.5pt}\it\Sigma}}: \quad
\xymatrixrowsep{15pt}
\xymatrixcolsep{20pt}\xymatrix{0 \ar[r]& \oplus_{j=1}^s
P_{y_j}' \ar[r]^{w_{_{\hspace{-0.5pt}\it\Sigma}}}& \oplus_{i=1}^r P_{x_i}' \ar[r] &N
\ar[r]& 0}$$ of
$N$ in $\rep(\Sa)$. On the other hand, by Lemma \ref{AR-translate}, $L$ has a minimal injective co-resolution
\vspace{-6pt} $$\zeta: \quad \xymatrixrowsep{20pt}
\xymatrixcolsep{20pt}\xymatrix{0 \ar[r] & L \ar[r] &
\oplus_{j=1}^s I_{y_j} \ar[r]^-{\nu(w)} &
\oplus_{i=1}^r I_{x_i} \ar[r]& 0} \vspace{-2pt}$$ in $\rep(Q).$ Restricting $\zeta$ to $\Sa$, we obtain a minimal injective co-resolution
\vspace{-3pt}  $$\zeta_{_{\it\Sigma}}: \quad
\xymatrixcolsep{22pt}\xymatrix{0 \ar[r]& L_{_{\hspace{-0.8pt}\it\Sigma}} \ar[r]&
\oplus_{j=1}^s I_{y_j}' \ar[r]^-{\nu(w)_{_{\hspace{-0.5pt}\it\Sigma}}} &
\oplus_{i=1}^r I_{x_i}' \ar[r] & 0} \vspace{-2pt}$$ of $L_{_{\hspace{-0.8pt}\it\Sigma}}$ in $\rep(\Sa)$. Moreover, it follows from the definition that $\nu(w)_{_{\hspace{-0.5pt}\it\Sigma}}=\nu_{_{\it\Sigma}}(w_{_{\hspace{-0.5pt}\it\Sigma}})\vspace{1pt}$, where $\nu_{_{\it\Sigma}}$ is the Nakayama functor for $\rep(\Sa)$. This implies that $L_{_{\hspace{-0.8pt}\it\Sigma}}=\DTr_{_{\hspace{-1pt}\it\Sigma}} M$. By Lemma \ref{AR-translate}, $L_{_{\hspace{-0.8pt}\it\Sigma}}$ is indecomposable and hence strongly indecomposable. Thus, $\xi_{_{\hspace{-1pt}\it\Sigma}}$ is an almost split sequence in $\rep(\Sa)$; see
\cite[(2.14)]{AuR}.

\vspace{1pt}

(2) Assume that $f: M\to N$ is an irreducible morphism in $\rep(Q)$.
By Lemma \ref{mrasm} and Theorem \ref{AR-sequence}, $\rep(Q)$ has a
minimal right almost split morphism $g: L\to N$. Then $f=gs$, where
$s : M \to L$ is a section. Hence,
$f_{_{\it\Sigma}}=g_{_{\it\Sigma}}s_{_{\it\Sigma}}$, where
$s_{_{\it\Sigma}}$ is clearly a section. If $N$ is not projective,
then $g_{_{\it\Sigma}}$ is minimal right almost split in $\rep(\Sa)$
by Statement (1). Otherwise, $N=P_x$ and $L={\rm rad}\hskip 0.5pt
P_x$ for some $x\in Q_0$. By the hypothesis, both $L$ and $N$ are
supported by $\Sa$. So $g_{_{\it\Sigma}}=g$, which is minimal right
almost split in $\rep(\Sa)$. In any case, $f_{_{\it\Sigma}}$ is
irreducible in $\rep(\Sa)$. The proof of the proposition is
completed.

{\center \section{Auslander-Reiten categories}}

In the following four sections, we shall be mainly concerned with
the study of Auslander-Reiten theory in $\rep^+\hspace{-1pt}(Q)$.
It is left to the reader to formulate the dual results for
$\rep^-(Q)$. In case $Q$ has no left infinite path, Reiten and Van
den Bergh proved that $\rep^+\hspace{-1pt}(Q)$ is right
Auslander-Reiten; see \cite{RVDB}. The main objective of this
section is to find the nece\-ssary and sufficient conditions for
${\rm rep}^+(Q)$ to be left or right Auslander-Reiten.

\medskip

We begin with studying some properties of irreducible morphisms in
$\rep^+\hspace{-1pt}(Q)$.

\medskip

\begin{Lemma}\label{irr-ker} If $f: M\to N$ is an irreducible epimorphism in $\rep^+\hspace{-1pt}(Q)$,  then the kernel of $f$ is finite dimensional.
\end{Lemma}

\noindent{\it Proof.} Let $f: M\to N$ be an irreducible epimorphism
in $\rep^+\hspace{-1pt}(Q)$. Since $\rep^+\hspace{-1pt}(Q)$ is Krull-Schmidt, we may assume
that $N$ is indecomposable; see \cite[(3.1), (3.2)]{Bau}. Since $N$ is not projective, by Theorem \ref{AR-sequence}(1),
there exists in $\rep(Q)$ an almost split sequence  
$\;\eta\hspace{-1pt}:
\xymatrixcolsep{15pt}\xymatrix{0\ar[r]& L \ar[r]& E \ar[r]^g& N
\ar[r]& 0,}$ where $L\in {\rm rep}^-(Q)$. This yields a pushout
diagram
\vspace{-10pt}
$$
\qquad \xymatrixrowsep{15pt}\xymatrixcolsep{20pt}\xymatrix{0 \ar[r] & X \ar[r] \ar[d]^v & \,M \ar[d]^u \ar[r]^f & \,N \ar@{=}[d] \ar[r] & \,0\\
0 \ar[r] & L \ar[r] & \,E \ar[r]^g & \,N \ar[r] & \;\,0.}$$ Set
${\it\Sigma}={\rm supp}\,M$ and $\Oa={\rm supp}\,E$. Let $\Ta$ be
the full subquiver of $\Oa$ generated by the vertices which are
successors in $Q$ of the vertices in $\Sa$. Since $\it\Sigma$ is
top-finite by Lemma \ref{Lemma1.6_shiping} and $Q$ is
interval-finite, $\Ta$ has no left infinite path. Since ${\rm
supp}\hskip 0.5pt N\subseteq \Sa\cap \Oa\subseteq \Ta$, restricting
$\eta$ to $\Ta$ yields a short exact sequence \vspace{-4pt}
$$\eta_{_{\hspace{-0.5pt}\it\Theta}}: \quad \xymatrixcolsep{20pt}\xymatrix{0 \ar[r]& L_{_{\hspace{-1pt}\it\Theta}}\ar[r] &
E_{_{\hspace{-1pt}\it\Theta}} \ar[r]^{g_{_{\hspace{-1pt}\it\Theta}}} & N \ar[r] & 0.}
\vspace{-2pt}$$
Suppose that ${\rm supp}\hspace{0.5pt}L_{_{\hspace{-1pt}\it\Theta}}$, that is $\Ta \cap \,{\rm supp\hskip 0.5pt} L$, contains infinitely
many vertices $x_i$, $i\in \N$. Since ${\rm supp\hskip 0.5pt}L$ is
socle-finite by Lemma \ref{Lemma1.6_shiping}, we may assume that
${\rm supp\hskip 0.5pt}L$ contains a path $p_i: x_i\rightsquigarrow
a$, for each $i\in \N$, where $a$ is some fixed vertex in ${\rm supp\hskip 0.5pt}L$. Since ${\rm supp}\hskip 0.5pt
L\subseteq \Oa$ and $\Ta$ is successor-closed in $\Oa$, the $p_i$ all lie in $\Ta$. Being locally finite,
by K\"onig's lemma, $\Ta\vspace{1pt}$ has a left infinite path ending with $a$, a contradiction. Thus $L_{_{\hspace{-1pt}\it\Theta}}\in \rep^b(Q)$, and consequently, $E_{_{\hspace{-1pt}\it\Theta}}\in \rep^+\hspace{-1pt}(Q)$. Note that $E_{_{\hspace{-1pt}\it\Theta}}$ is a sub-representation of $E$ since $\Ta$ is successor-closed in $\Oa$. Moreover, the support of ${\rm
Im}(u)$ is contained in $\Sa\cap \Oa\subseteq \Ta$. Thus $u=qu'$,
where $u':M \to E_{_{\hspace{-1pt}\it\Theta}}$ is the co-restriction of $u$, and $q:
E_{_{\hspace{-1pt}\it\Theta}}\to E$ is the inclusion. This yields a factorisation $f=(gq) u'$
in $\rep^+\hspace{-1pt}(Q)$. Thus $gq$ is a retraction or $u'$ is a section. Since $g$ is not a
retraction, the first case does not occur. In particular, $u=qu'$ is a
monomorphism, and so is $v$. Since ${\rm supp} \hskip 0.5pt L$ has
no right infinite path by Lemma \ref{Lemma1.6_shiping}, nor does
${\rm supp}\hskip 0.5pt X$. On the other hand, $X\in \rep^+\hspace{-1pt}(Q)$ since it is the kernel of $f$. By Corollary \ref{newcor1.9}, $X\in
\rep^b(Q)$. The proof of the lemma is completed.

\medskip

\begin{Lemma}\label{infdimirr} Let $f: M\to N$ be an irreducible morphism in $\rep^+(Q)$.
If $\,M$ is infinite dimensional, then $N$ is infinite dimensional
while $\DTr N$ is finite dimensional.
\end{Lemma}

\noindent{\it Proof.} Suppose that $M$ is infinite dimensional. If
$N$ is finite dimensional, then $f$ is an epimorphism. By Lemma
\ref{irr-ker}, the kernel of $f$ is finite dimensional, and
consequently, $M$ is finite dimensional. This contradiction shows
that $N$ is infinite dimensional. For proving the second part of the
lemma, we may assume that $N$ is indecomposable and not projective.
Then $\rep(Q)$ has an almost split sequence \vspace{-4pt}
$$\qquad\qquad  \eta: \quad \xymatrixcolsep{20pt}\xymatrix{ 0\ar[r]& L \ar[r]& E \ar[r]^{g}& N \ar[r]& 0,}$$ where $L\in {\rm rep}^-(Q)$. Suppose that $L$ is infinite dimensional. By Corollary \ref{newcor1.9}, ${\rm supp}\hskip 0.5pt L$ has a left infinite path. Since ${\rm supp}\hskip 0.5pt M$ has no left infinite path by Lemma \ref{Lemma1.6_shiping}(2), there exists some $a\in Q_0$ such that $L(a)\ne 0$ but $M(a)=0$. Let $\Sa$ be the successor-closed subquiver of $Q$ generated by $a$ and the vertices in the support of $M\oplus N$. Then $\Sa$ is top-finite.
By Proposition \ref{prop2.9}(1), restricting $\eta$ to $\Sa$, we get an almost split sequence
\vspace{-8pt}
$$\qquad\qquad \eta_{_{\it\Sigma}}: \quad \xymatrixcolsep{20pt}\xymatrix{0\ar[r]& L_{_{\hspace{-0.8pt}\it\Sigma}} \ar[r]& E_{_{\it\Sigma}} \ar[r]^{g_{_{\it\Sigma}}}& N \ar[r]& 0}$$ in $\rep(\Sa)$. By the dual of Lemma \ref{rep+restrict}(1), $L_{_{\hspace{-0.8pt}\it\Sigma}}\in \rep^-(\Sa)$, and hence  ${\rm supp}\hspace{0.4pt}L_{_{\hspace{-0.8pt}\it\Sigma}}$ is socle-finite by Lemma \ref{Lemma1.6_shiping}. On the other hand, since ${\rm supp}\hspace{0.4pt}L_{_{\hspace{-0.8pt}\it\Sigma}}$ is a subquiver of the top-finite quiver $\Sa$, it is finite. As a consequence, $\eta_{_{\it\Sigma}}$ lies in
$\rep^+(\Sa)$ and hence, it is an almost split sequence in $\rep^+(\Sa)$.

Finally, by Lemma \ref{rep+restrict}(2), $f$ lies in $\rep^+(\Sa)$, and hence it is an irreducible morphism in $\rep^+(\Sa)$. Thus we
have an irreducible morphism $h: L_{_{\hspace{-0.8pt}\it\Sigma}}\to M$. Since $L_{_{\hspace{-0.8pt}\it\Sigma}}$ is finite
dimensional while $M$ is infinite dimensional, $h$ is a
monomorphism. Since $L_{_{\hspace{-0.8pt}\it\Sigma}}(a)=L(a)\ne 0$, we have $M(a)\ne 0$, a contradiction. The proof of the lemma is completed.

\medskip

\begin{Cor} \label{Cor3.3}Let $\xymatrixcolsep{15 pt}\xymatrix{0\ar[r] & L\ar[r] & M\ar[r] & N\ar[r] & 0}$ be an almost split sequence in $\rep(Q)$ with $L\in \rep^-(Q)$ and $N\in \rep^+\hspace{-1pt}(Q)$, and let $X$ be an indecomposable direct summand of $M$.

\begin{enumerate}[$(1)$]

\vspace{-1pt}

\item If $L$ is infinite dimensional, then $X$ is finitely presented if and only if $X$ is finite dimensional.

\vspace{1pt}

\item If $N$ is infinite dimensional, then $X$ is finitely co-presented if and only if $X$ is finite dimensional.


\end{enumerate} \end{Cor}

\noindent{\it Proof.} By assumption, there exists an irreducible
morphism $f: X\to N$ in $\rep(Q)$. Suppose that $X\in \rep^+\hspace{-1pt}(Q)$.
Then $f$ is irreducible in $\rep^+\hspace{-1pt}(Q)$. If $X$ is infinite
dimensional, then $L$ is finite dimensional by Lemma
\ref{infdimirr}. This proves Statement (1). Using the dual of Lemma
\ref{infdimirr}, we may prove Statement (2). The proof of the
corollary is completed.

\medskip

For a morphism in $\rep^+\hspace{-1pt}(Q)$,  we shall relate its
irreducibility in $\rep^+\hspace{-1pt}(Q)$ to that in $\rep(Q)$.

\medskip

\begin{Lemma} \label{asm} Let $M$ be an indecomposable representation in $\rep^+\hspace{-1pt}(Q)$.

\vspace{-1pt}

\begin{enumerate}[$(1)$]

\item If $M$ is finite dimensional, then $\rep^+\hspace{-1pt}(Q)$ has a minimal left almost split morphism $f: M\to N$, which is also minimal left almost split in $\rep(Q)$.

\vspace{1pt}

\item If $\,\DTr M$ is finite dimensional, then $\rep^+\hspace{-1pt}(Q)$ has a minimal right almost split morphism $g: L\to M$, which is also minimal right almost split in $\rep(Q)$.

\end{enumerate}\end{Lemma}

\noindent{\it Proof.} (1) Suppose that $M\in \rep^b(Q)$. If $M$ is injective then $M\cong
I_x$ for some $x\in Q_0$ by Lemma \ref{cor1.14shiping}(2), and in this case, the projection $p: M\to M/\,{\rm soc}\,M$ is minimal left almost split in  $\rep^+\hspace{-1pt}(Q)$ and in  $\rep(Q)$. Otherwise, by Corollary
\ref{cor_ARS_exist}(2), $\rep(Q)$ has a minimal left almost split morphism $f: M\to N$ which lies in
$\rep^+\hspace{-1pt}(Q)$. Thus, $f$ is a minimal left almost split morphism in
$\rep^+\hspace{-1pt}(Q)$.

\vspace{1pt}

(2) Suppose that $\DTr M\in \rep^b(Q)$. If $\DTr M=0$, then $M\cong
P_x$ for some $x\in Q_0$, and in this case, the inclusion $q: {\rm rad}\,M\to M$ is minimal right almost split in $\rep(Q)$ and in $\rep^+\hspace{-1pt}(Q)$. Otherwise, by Corollary
\ref{cor_ARS_exist}(1), $\rep(Q)$ has an minimal right almost split morphism $g: L\to M$, which lies in $\rep^+\hspace{-1pt}(Q)$. Hence, $g$ is  a minimal right almost split morphism in $\rep^+\hspace{-1pt}(Q)$. The proof of
the lemma is completed.

\medskip

\begin{Cor} \label{irr+rep} Let $f: M\to N$ be a morphism in $\rep^+\hspace{-1pt}(Q)$. If $M, N$ are indecomposable, then $f$ is irreducible in $\rep^+\hspace{-1pt}(Q)$ if and only if it is irreducible in $\rep(Q)$.
\end{Cor}

\noindent{\it Proof.} Suppose that $f$ is irreducible  in $\rep^+\hspace{-1pt}(Q)$ with $M, N$ indecomposable. Assume
first that $M\in \rep^b(Q)$. By Lemma \ref{asm}(1), $\rep^+\hspace{-1pt}(Q)$ has
a minimal left almost split morphism $g: M\to L$, which is minimal
left almost split in $\rep(Q)$. If $f$ is irreducible in
$\rep^+\hspace{-1pt}(Q)$, then $f=ug$ for some retraction $u: L\to N$, and
hence, $f$ is irreducible in $\rep(Q)$; see \cite[(2.4)]{AuR}.

\vspace{1pt}

Assume now that $M$ is infinite dimensional. By Lemma
\ref{infdimirr}, $\DTr N\in \rep^b(Q)$. By Lemma \ref{asm}(2),
$\rep^+\hspace{-1pt}(Q)$ has a minimal right almost split morphism $h: L\to N$,
which is minimal right almost split in $\rep(Q)$. If $f$ is
irreducible in $\rep^+\hspace{-1pt}(Q)$, then
 $f=hv$ for some section $v: N\to L$, and consequently, $f$ is irreducible in $\rep(Q)$.
The proof of the corollary is completed.

\medskip

The following result is essential in our investigation, since it allows us to apply some well-established results of the representation theory of finite quivers.

\medskip

\begin{Prop}\label{ARrep+} If $\xymatrixcolsep{15 pt}\xymatrix{0\ar[r] &
L\ar[r] & M\ar[r] & N \ar[r] & 0}$ is a short exact sequence in $\rep(Q)$, then it is an almost split sequence in
$\rep^+\hspace{-1pt}(Q)$ if and only if it is an almost split sequence in $\rep(Q)$ with
$L\in \rep^b(Q)$. \end{Prop}

\noindent{\it Proof.} The sufficiency follows from Corollary \ref{cor_ARS_exist}(2) and the uniqueness of an almost split sequence.
For the necessity, assume that $\eta: \; \xymatrixcolsep{15 pt}\xymatrix{0\ar[r] & L\ar[r] &
M\ar[r] & N \ar[r] & 0}$ is an almost split sequence in ${\rm
rep}^+(Q)$. By Lemma \ref{irr-ker}, $L\in \rep^b(Q)$. Since $L$ is not injective, by
Corollary \ref{cor_ARS_exist}(2), $\rep(Q)$ has an almost split sequence $\zeta$ starting with $L$, which is also an almost split sequence in $\rep^+\hspace{-1pt}(Q)$. Then $\zeta$ is isomorphic to $\eta$. In other words, $\eta$ is an almost split sequence in ${\rm rep}(Q)$. The proof of the proposition is completed.

\medskip

We are ready to give conditions for $\rep^+\hspace{-1pt}(Q)$ to be left or right Auslander-Reiten.


\smallskip

\begin{Theo} \label{RARcat} Suppose that $Q$ is a strongly locally finite quiver. 

\begin{enumerate}[$(1)$]

\vspace{-2pt}

\item $\rep^+\hspace{-1pt}(Q)$ is left Auslander-Reiten if and only if $Q$ has no right infinite
path.

\vspace{1pt}

\item $\rep^+\hspace{-1pt}(Q)$ is right Auslander-Reiten if and
only if $Q$ has no left infinite path, or else $Q$ is a left
infinite or double infinite path.

\end{enumerate}

\end{Theo}

\noindent{\it Proof.} (1) Suppose first that $Q$ has a right infinite path $p$ with an initial arrow $x\to y$. In particular, $P_y$ is
infinite dimensional. By Proposition \ref{ARrep+}, $\rep^+\hspace{-1pt}(Q)$
admits no almost split sequence starting with $P_y$. Suppose that
$\rep^+\hspace{-1pt}(Q)$ has a minimal left almost split epimorphism $f: P_y\to L$. By Lemma \ref{pc-ih}(2), $P_y$ is injective in $\rep^+\hspace{-1pt}(Q)$. In particular, the inclusion $q: P_y\to P_x$ is a section, which is absurd. Thus $\rep^+\hspace{-1pt}(Q)$
is not left Auslander-Reiten. Conversely, assume that $Q$ contains no right infinite path. Then
$\rep^+\hspace{-1pt}(Q)=\rep^b(Q)$. Let $M$ be an indecomposable object in
$\rep^b(Q)$. If $M$ is not injective then, by Corollary
\ref{cor_ARS_exist}(2), $\rep^+\hspace{-1pt}(Q)$ admits an almost split sequence
starting with $M$. Otherwise, by Proposition
\ref{cor1.14shiping}(2), $M=I_x$ for some $x\in Q^+$. Thus $M\to
M/\,{\rm soc}\hskip 0.5pt M$ is a minimal left almost split
epimorphism in $\rep^+(Q)$. That is, $\rep^+(Q)$ is left Auslander-Reiten.

\vspace{1pt}

(2) For proving the sufficiency, let $N$ be an indecomposable object
in $\rep^+\hspace{-1pt}(Q)$. If $N$ is projective, then the
inclusion $q: {\rm rad}\hskip 0.5pt N\to N$ is a minimal right
almost split monomorphism in $\rep^+\hspace{-1pt}(Q)$. Otherwise,
$\rep(Q)$ admits an almost split sequence
$$\eta:
\xymatrixcolsep{20pt}\xymatrix{0\ar[r] & L\ar[r] & M\ar[r] &
N\ar[r] & 0},$$ where $L$ is an indecomposable non-injective object
in $\rep^-(Q)$. If $Q$ contains no left infinite path, then
$\rep^-(Q)=\rep^b(Q)$, and hence $L$ is finite dimensional. If $Q$
is a left infinite or double infinite path, then every
indecomposable non-injective object in $\rep^-(Q)$ is finite
dimensional; see (\ref{5.1})(2) below, and hence $L$ is finite
dimensional. In any case, by Proposition \ref{ARrep+}, $\eta$ is
an almost split sequence in $\rep^+\hspace{-1pt}(Q)$. This shows
that $\rep^+\hspace{-1pt}(Q)$ is right Auslander-Reiten.

Conversely, assume that $\rep^+\hspace{-1pt}(Q)\vspace{1pt}$ is
right Auslander-Reiten. By Proposition \ref{ARrep+}, $\DTr X\in
\rep^b(Q)$ for any indecomposable non-projective object $X$ in
$\rep^+\hspace{-1pt}(Q)$. Suppose that $Q$ contains a left
infinite path $p$. Choose arbitrarily a vertex $a$ lying on $p$.
Then $Q$ contains a left infinite path \vspace{-2pt}
$$\xymatrixcolsep{20pt}\xymatrix{\cdots \ar[r] & a_n \ar[r]^-{\alpha_n} & a_{n-1}\ar[r] &  \cdots \ar[r] &  a_1
\ar[r]^-{\alpha_1} & a_0= a.}\vspace{-2pt}$$ We claim that $a$ is the starting point of at most
one arrow, while $\alpha_1$ is the only arrow
ending in $a$. Indeed, assume that $a^+=\{\alpha_i: a\to b_i \mid i=1,
\ldots, r\}$ with $r>1$. Then $S_a$ is not projective with a minimal
projective resolution
$$\xymatrixcolsep{20pt}\xymatrix{0\ar[r]& {\oplus}_{i=1}^r P_{b_i} \ar[r] & P_a \ar[r] & S_a \ar[r] & 0.}$$
By Lemma \ref{AR-translate}, $\DTr S_a$ has a minimal injective co-resolution
\vspace{-2pt}
$$\xymatrixcolsep{20pt}\xymatrix{0\ar[r] & \DTr S_a\ar[r] & \oplus_{i=1}^r I_{b_i} \ar[r]& I_a\ar[r] & 0.}$$
For each $n\ge 0$, since ${\rm dim} I_{b_i}(a_n)\ge {\rm
dim} I_{a}(a_n)>0$ for all $1\le i\le r$, we get
$${\rm dim}\,(\DTr S_a)(a_n) ={\textstyle\sum}_{i=1}^r {\rm
dim}\,I_{b_i}(a_n) - {\rm dim}\,I_a(a_n) \ge
{\textstyle\sum}_{i=2}^r\, {\rm dim}\, I_{b_i}(a_n) > 0.$$
Therefore, $\DTr S_a$ is infinite dimensional, a contradiction.
Next, assume that there exists an arrow $\beta: b\to a$ different
from $\alpha_1$. Similarly, $\DTr S_b$ has a minimal injective
co-resolution \vspace{-2pt}
$$\xymatrixcolsep{20pt}\xymatrix{0\ar[r] & \DTr S_b\ar[r] & I_a
\oplus I \ar[r] & I_b\ar[r] & 0},$$ where $I\in {\rm inj}(Q)$.
Note, for each $n\ge 1$, that there exists a $k$-monomorphism
$\phi_n: I_b(a_n)\to I_a(a_n): \tau\mapsto \beta \tau.$ Since
$\alpha_1\ne \beta$, the path $\alpha_1\cdots \alpha_n$ lies in
$I_a(a_n)$ but not in the image of $\phi_n$. As a consequence,
${\rm dim}\,I_a(a_n)> {\rm dim}\,I_b(a_n)$. Therefore,
$${\rm dim}\,(\DTr S_b)(a_n) = {\rm
dim}\, I_a(a_n)+{\rm dim}\, I(a_n) -{\rm dim}\, I_b(a_n)>0,$$ for
all $n\ge 1$. In particular, $\DTr S_b$ is infinite dimensional, a
contradiction. Our claim is established, from which we infer that
$Q$ is a double infinite path if $p$ is contained in a double
infinite path, and otherwise, $Q$ is a left infinite path. The
proof of the theorem is completed.

\medskip

We conclude this section with an immediate consequence of Theorem \ref{RARcat}.

\medskip

\begin{Cor}\label{AR-category} If $Q$ is a strongly locally finite quiver, then $\rep^+\hspace{-1pt}(Q)$ is
Auslan\-der-Reiten if and only if either $Q$ has no infinite path
or $Q$ is a left infinite path.
\end{Cor}


\vspace{-5pt}

{\center \section{Auslander-Reiten components}}

The objective of this section is to describe the shapes of the
Auslander-Reiten components of $\rep^+\hspace{-1pt}(Q)$, which has
been shown to be a Hom-finite abelian $k$-category. In contrast to
the finite case, we shall see that many new phenomena occur, such
as the number of preinjective components varies from zero to the
infinity, and there exist four types of regular components.

\medskip

First of all, we recall the notion of a section of a translation
quiver since it is essential in the description of the shapes of
Auslander-Reiten components. Let $\Ga$ be a connected valued or
non-valued translation quiver with translation $\tau$; see, for
example, \cite{HPR, R3}. A connected convex subquiver $\Da$ of
$\Ga$ is called a {\it section} if it contains no oriented cycle
and meets each $\tau$-orbit exactly once; see \cite[(2.1)]{L4}.
Now, we say that a section $\Da$ of $\Ga$ is {\it right-most} or
{\it left-most} if the vertices in $\Ga$ are all of the form
$\tau^nx$ or all of the form $\tau^{-n}x$ with $n\in \N$ and $x\in
\Da_0$, respectively.

\medskip

\begin{Lemma}\label{section} Let $(\Ga, \tau)$ be a connected translation quiver with no oriented cycle,
and let $\Da$ be a full subquiver of $\Ga$ meeting any given
$\tau$-orbit at most once. Then

\vspace{-1pt}

\begin{enumerate}[$(1)$]

\item If $\Da$ is successor-closed in $\Ga$ and has the following
property$\,:$ for each arrow $x\to \tau^ny$ in $\Ga$ with $n\ge
0$, $y\in \Da$ implies $x\in \Da$ or $\tau^-x\in \Ga$, then it is
a right-most section of $\Ga$.

\vspace{1pt}

\item  If $\Da$ is predecessor-closed in $\Ga$ and has the
following property$\,:$ for each arrow $\tau^{-n}x\to y$ in $\Ga$
with $n\ge 0$,  $x\in \Da$ implies $y\in \Da$ or $\tau y\in \Ga$,
then it is a left-most section of $\Ga$.

\vspace{1pt}

\end{enumerate}
\end{Lemma}

\noindent{\it Proof.} We shall prove only the first statement.
Assume that $\Da$ satisfies the condition stated in (1). Let $\Sa$
be a connected component of $\Da$. Then $\Sa$ contains no oriented
cycle and meets any $\tau$-orbit in $\Ga$ at most once. Moreover,
since $\Da$ is successor-closed in $\Ga$, so is $\Sa$. In
particular, $\Sa$ is convex in $\Ga$ and the vertices in the
$\tau$-orbit of some vertex $z$ in $\Sa$ are all of the form
$\tau^rz$ with $r\ge 0$. We claim that every vertex $a$ in $\Ga$
lies in the $\tau$-orbit of some vertex in $\Sa$. Since $\Ga$ is
connected, we may assume that $\Ga$ contains an edge $a$ --- $b$,
where $b$ lies in the $\tau$-orbit of some $x\in \Sa_0$. Then
$b=\tau^nx$ for some $n\ge 0$. If $\tau^{-n-1}a\in \Ga$, then
either $x\to \tau^{-n-1}a$ or $x\to \tau^{-n}a$ is an arrow in
$\Ga$. Since $\Sa$ is successor-closed in $\Ga$, we have
$\tau^{-n}a\in \Sa$ or $\tau^{-n-1}a\in \Sa$. Suppose now that
$\tau^{-n-1}a\not\in \Ga$. Then there exists some $0\le m\le n$
such that $\tau^{-m}a\in \Ga$ while $\tau^{-m-1}a \not\in \Ga.$
This yields an arrow $\tau^{n-m}x\to \tau^{-m}a$ or $\tau^{-m}a
\to \tau^{n-m}x$ in $\Ga$. If $n=m$, since
$\tau^{-m-1}a\not\in\Ga$, it follows from the condition stated in
(1) that $\tau^{-m}a\in \Da$, and hence $\tau^{-m}a\in \Sa$. If
$m<n$, then either $\tau^{-m}a \to \tau^{n-m-1}x$ or $\tau^{-m}a
\to \tau^{n-m}x$ is an arrow in $\Ga$. By the property of $\Da$
stated in (1), we have $\tau^{-m}a\in \Da$. Since $\Da$ is
successor-closed, we get $\tau^{n-m-1}x\in \Da$. Since $\Da$ meets
any $\tau$-orbit at most once, we see that $m=n-1$ and
$\tau^{-m}a\to x$ is an arrow in $\Ga$. Thus $\tau^{-m}a\in \Sa$.
This establishes our claim. As a consequence, $\Sa$ is a
right-most section of $\Ga$. Finally, since $\Da$ meets any
$\tau$-orbit at most once, we have $\Da=\Sa$. The proof of the
lemma is completed.

\medskip

Let $(\Ga, \tau)$ be a connected valued or non-valued translation
quiver, and let $x$ be a vertex in $\Ga$. One says that $x$ is
{\it projective} or {\it injective} if $\tau x$ or $\tau^-x$ is
not defined in $\Ga$, respectively. Moreover, we say that $x$ is
{\it left stable} if $\tau^nx$ is defined for all $n\in \N$; {\it
right stable} if $\tau^{-n}x$ is defined for all $n\in \N$; and
{\it stable} if it is both left and right stable. Furthermore,
$\Ga$ is called {\it left stable}, {\it right stable}, or {\it
stable} if its vertices are all left stable, all right stable, or
all stable, respectively.

\medskip

Given a connected quiver $\Da$ with no oriented cycle, one
constructs a stable translation quiver $\Z \Da$; see, for example,
\cite[Section 2]{L4}. We denote by $\N\Da$ the full translation
subquiver of $\Z\Da$ generated by the vertices $(n, x)$ with $n\ge
0$ and $x\in \Da_0$, and by $\N^-{\hskip -3pt}\Da$ the one
generated by the vertices $(n, x)$ with $n\le 0$ and $x\in \Da_0$.
It is evident that $\N\Da$ is right stable with a left-most
section generated by the vertices $(0, x)$ with $x\in \Da_0$,
while $\N^-{\hskip -3pt}\Da$ is left stable with a right-most
section generated by the vertices $(0, x)$ with $x\in \Da_0$.
Assume that $\Da$ is a section of $\Ga$. Then $\Ga$ is isomorphic
to the full translation subquiver of $\Z{\hskip -3pt}\Da$
generated by the vertices $(-n, x)$, where $n\in \Z$ and $x\in
\Da_0$ such that $\tau^nx\in \Ga$; see \cite[(2.3)]{L4}. In
particular, if $\Da$ is left-most or right-most, then $\Ga$ embeds
in $\N\Da$ or $\N^-{\hskip -3pt}\Da$, respectively.

\medskip

Let $\mathcal A\vspace{1pt}$ be a Hom-finite Krull-Schmidt additive
$k$-category. For indecomposable objects $X, Y\in \mathcal A$, write
${\rm irr}(X, Y)={\rm rad}(X, Y)/{\rm rad}^2(X, Y)\vspace{1pt}$, and
its dimensions over ${\rm End}(X)/{\rm rad}(X, X)$ and ${\rm
End}(Y)/{\rm rad}(Y, Y)$ are denoted by $d\hspace{0.4pt}'_{XY}$ and
$d_{XY}\vspace{1.5pt}$, respectively. It is well known that
$\mathcal A$ has an irreducible morphism $f: X\to Y$ if and only if
$d_{XY}>0\,$, and in this case,
$d\hspace{0.4pt}'_{_{XY}}\vspace{1pt}$ and $d_{_{XY}}$ are the
maximal integers for which $\mathcal A$ has irreducible morphisms
$g: X^{d\hspace{0.4pt}'_{XY}}\to Y\vspace{1pt}$ and $h: X\to
Y^{d_{XY}},$ respectively, where $M^n$ denotes the direct sum of $n$
copies of $M\,$; see  \cite[(3.4)]{Bau}. The \emph{Auslander-Reiten
quiver} $\Ga_{_{\hskip -3pt \mathcal A}}$ of $\mathcal A$ is a
valued translation quiver defined as follows: the vertex set is a
complete set of representatives of isomorphism classes of the
indecomposable objects in $\mathcal A$; for vertices $X, Y$, there
exists a unique valued arrow $X\to Y$ with valu\-ation $(d_{_{XY}},
d\hspace{0.4pt}'_{_{XY}})$ if and only if $d_{_{X,Y}}>0\,$; and the
translation $\tau$ is defined so that $\tau Z=X$ if and only if
$\mathcal A$ has an almost split sequence
$\xymatrixcolsep{15pt}\xymatrix{X\ar[r] & Y \ar[r] & Z}$; see
\cite[(2.1)]{L3}. A valuation $(d_{_{XY}},
d\,'_{_{XY}})\vspace{1pt}$ is called {\it symmetric} if
$d_{_{XY}}=d_{_{XY}}'$, and {\it trivial} if
$d_{_{XY}}=d_{_{XY}}'=1$. For some technical reasons; see
(\ref{sqproj}) below, we shall replace each symmetrically valued
arrow $X\to Y$ by $d_{_{XY}}$ unvalued arrows from $X$ to $Y$. In
this way, $\Ga_{_{\hskip -3pt \mathcal A}}$ becomes a partially
valued translation quiver with possible multiple arrows in which all
valuations are non-symmetric. The connected components of
$\Ga_{_{\mathcal{A}}}$ are called the {\it Auslander-Reiten
components} of $\mathcal A$.

\medskip

Now, we specialize to the Auslander-Reiten quiver
$\Ga_{\rep^+\hspace{-0.5pt}(Q)}$ of
$\rep^+\hspace{-1pt}(Q)\vspace{1pt}$. We choose its vertex set in
such a way that it contains the $S_x$ with $x\in Q_0$, the $P_x$
with $x\in Q_0$, and the $I_x$ with $x\in Q^+$. Recall that $Q^+$
denotes the full subquiver of $Q$ generated by the vertices $x$
such that $I_x$ is finite dimensional. We shall see that the
arrows in $\Ga_{\rep^+\hspace{-0.5pt}(Q)}$ are all symmetrically
valued; see (\ref{valuation}) below, and hence
$\Ga_{\rep^+\hspace{-0.5pt}(Q)}\vspace{1pt}$ is an unvalued
translation quiver with multiple arrows in the sense of \cite[page
47]{R3}. For convenience, we say that a representation $M$ in
$\Ga_{\rep^+\hspace{-0.5pt}(Q)}$ is {\it pseudo-projective} if
$\DTr M$ is infinite dimensional, or equivalently, $\DTr M\not\in
\rep^+\hspace{-1pt}(Q)$. The following observation, which follows
from Proposition \ref{ARrep+} and Corollary \ref{cor_ARS_exist},
clarifies the relation between the Auslander-Reiten translations
$\tau$ and $\tau^-$ for $\Ga_{\rep^+(Q)}$ and the Auslander-Reiten
translations $\DTr$ and $\TrD$ for $\rep(Q)$.

\medskip

\begin{Lemma}\label{psproj} If $M$ is a representation lying in $\Ga_{\rep^+\hspace{-0.5pt}(Q)}$, then

\vspace{-1pt}

\begin{enumerate}[$(1)$]

\item $\tau M$ is defined in
$\Ga_{\rep^+\hspace{-0.5pt}(Q)}\vspace{1pt}$ if and only if $\,M$
is neither projective nor pseudo-projective, and in this case,
$\tau M\cong \DTr M$, which is of positive finite dimension$\,;$

\vspace{1pt}

\item  $\tau^-M$ is defined in $\Ga_{\rep^+\hspace{-0.5pt}(Q)}$ if
and only if $M$ is finite dimensional and not injective, and in
this case, $\tau^-\hspace{-1pt}M\cong \TrD M;$

\item $M$ is left stable or right stable in
$\Ga_{\rep^+\hspace{-0.5pt}(Q)}$ if and only if $M$ is
$\DTr$-stable or $\TrD$-stable in $\rep(Q)$, respectively.

\end{enumerate} \end{Lemma}

\medskip

\noindent{\sc Remark.} In other words, $M$ is a projective vertex in
$\Ga_{\rep^+(Q)}$ if and only if $M$ is a projective or
pseudo-projective representation in $\rep^+(Q)$. Moreover, $M$ is an
injective vertex in $\Ga_{\rep^+(Q)}\vspace{1pt}$ if and only if $M$
is an injective or infinite dimensional representation in
$\rep^+(Q)$.

\medskip

\begin{Lemma}\label{sqproj} Let $P_{_{\hspace{-1pt}Q}}$ be the full subquiver of $\Ga_{\rm rep^+(Q)}\vspace{1pt}$
generated by the $P_x$ with $x\in Q_0$, and let $I_{_{Q}}$ be the one generated by the $I_x$ with $x\in Q^+$.

\vspace{0pt}

\begin{enumerate}[$(1)$]

\item The subquiver $P_{_{\hspace{-1pt}Q}}$ is predecessor-closed
in $\Ga_{\rep^+\hspace{-0.5pt}(Q)}$ and isomorphic to
$Q^{\hspace{0.3pt}\rm op}$.

\vspace{0pt}

\item The subquiver $I_{_{\hspace{-1pt}Q}}$ is successor-closed in
$\Ga_{\rm rep^+(Q)}$ and isomorphic to $(Q^+)^{\rm op}$.

\end{enumerate}\end{Lemma}

\noindent{\it Proof.} We prove only the first statement, since the
second one follows dually. For $x, y\in Q_0$, denote by $n_{xy}$
the number of arrows in $Q$ from $x$ to $y$. By Proposition
\ref{prop1.3_shiping}, $\End(P_x)\cong \End(P_y)\cong k$ and ${\rm
irr}(P_y, P_x)$ has $k$-dimension $n_{xy}$. Thus
$\Ga_{\rep^+\hspace{-0.5pt}(Q)}$ contains a valued arrow $P_y\to
P_x$ if and only if $n_{xy}>0$, and in this case, the valuation is
$(n_{xy}, n_{xy})$. By definition, the symmetrically valued arrow
$P_y\to P_x$ is replaced by $n_{xy}$ unvalued arrows from $P_y$ to
$P_x$. Thus, $P_{_{\hspace{-1pt}Q}}\cong Q^{\hspace{0.4pt} \rm
op}$. Moreover, if $M\to P_x$ with $x\in Q$ is an arrow in
$\Ga_{\rep^+\hspace{-0.5pt}(Q)}$, then $M$ is a direct summand of
${\rm rad}\, P_x$, and hence $M=P_y$ for some $y\in Q$. Thus
$P_{_{\hspace{-1pt}Q}}$ is predecessor-closed in
$\Ga_{\rep^+\hspace{-0.5pt}(Q)}$. The proof of the lemma is
completed.

\medskip

The following result is well
known in the finite case.

\medskip

\begin{Lemma}\label{nocycle} If $Q$ is connected, then $\Ga_{\rep^+\hspace{-0.5pt}(Q)}$ contains
an oriented cycle if and only if $Q$ is finite of Euclidean type.
\end{Lemma}

\noindent{\it Proof.} We only need to consider the case where $Q$
is connected and infinite. Suppose that
$\Ga_{\rep^+\hspace{-0.5pt}(Q)}$ contains an oriented cycle
$$\;\eta: \xymatrixcolsep{18pt}\xymatrix{M_1 \ar[r]& M_2 \ar[r]&
\cdots \ar[r]& M_n = M_1.}$$ Since $Q$ has no oriented cycle, it
follows from Lemma \ref{sqproj} that none of the $M_i$ is
projective. If some of the $M_i$ is infinite dimensional, then the
$M_i$ are all infinite dimensional and the $\DTr M_i$ are of
finite positive dimension by Lemma \ref{infdimirr}. Thus
$\Ga_{\rep^+\hspace{-0.5pt}(Q)}$ has an oriented cycle
$\xymatrixcolsep{18pt}\xymatrix{\tau M_1\ar[r] & \tau M_2\ar[r] &
\cdots \ar[r] & \tau M_n=\tau M_1}$, which contains only finite
dimensional representations. Thus, we may assume that the $M_i$
with $1\le i\le n$ are all finite dimensional. In particular,
${\rm supp}(M_1\oplus \cdots \oplus M_n)$ is contained in a finite
connected full subquiver $\Sa$ of $Q$. Observing that $\eta$ is
also an oriented cycle in $\Ga_{\rep({\it\Sigma})}$, we see that
$\Sa$ is of Euclidean type. Since $Q$ is connected and infinite,
$\Sa$ is contained in a connected finite full subquiver $\Da$ of
$Q$ which is of wild type. Again, $\eta$ is an oriented cycle in
$\Ga_{\rep({\it\Delta})}$, a contradiction. The proof of the lemma
is completed.

\medskip

The preceding lemma yields the following important consequence.

\medskip

\begin{Lemma}\label{infrepsection} Let $\Ga$ be a connected component of $\Ga_{\rep^+\hspace{-0.5pt}(Q)}$.

\begin{enumerate}[$(1)$]

\item If $\Ga$ contains infinite dimensional representations, then such representations generate a right-most section of $\Ga$.

\vspace{1pt}

\item If $\Ga$ contains pseudo-projective representations, then such representations gene\-rate a left-most section of $\Ga$.

\vspace{1pt}

\end{enumerate} \end{Lemma}

\noindent{\it Proof.} (1) By Lemma \ref{nocycle}, $\Ga$
has no oriented cycle. Assume that the full subquiver $\Da$ of
$\Ga$ generated by the infinite dimensional representations is
non-empty. By Proposition \ref{ARrep+}, $\Da$ meets any $\tau$-orbit
in $\Ga$ at most once, and by  Lemma \ref{infdimirr}, $\Da$ is
successor-closed in $\Ga$. Let $M\to \tau^nN$ be an arrow in $\Ga$, where $n\ge 0$, $N\in \Da$,
and $M\not \in \Da$. Then $M$ is finite dimensional. If $M$ is injective then, by Proposition
\ref{cor1.14shiping}, $M=I_x$ for some $x\in Q^+$. By Lemma
\ref{sqproj}, $n=0$ and $N=I_y$ for some $y\in Q^+$, contrary to that $N\in \Da$. Thus $M$ is not injective, and by Lemma
\ref{psproj}(2), $\tau^- M\in \Ga$. It follows then from Lemma \ref{section}(1) that
$\Da$ is a right-most section of $\Ga$.

(2) Assume that the full subquiver $\Sa$ of $\Ga$ generated by the
pseudo-projective representations is non-empty. Clearly, $\Sa$ meets
any $\tau$-orbit in $\Ga$ at most once. Fix an arrow $M\to N$ in
$\Ga$. Suppose first that $N\in \Sa$.  Then $\rep(Q)$ admits an
almost split sequence $\xymatrixcolsep{15pt}\xymatrix{0\ar[r] & \DTr
N\ar[r] & E\ar[r] & N\ar[r] & 0}$, where $\DTr N\not\in \rep^+\hspace{-1pt}(Q)$.
By Corollary \ref{irr+rep}, an irreducible morphism $f: M\to N$ in
$\rep^+\hspace{-1pt}(Q)$ is irreducible in $\rep(Q)$. Thus there exists an
irreducible morphism $g: \DTr N\to M$ in $\rep(Q)$. If $M\not\in
\Sa$, by Lemma \ref{asm}(2), $\rep^+\hspace{-1pt}(Q)$ has a minimal right almost
split morphism $h: L\to M$ which is minimal right almost split in
$\rep(Q)$. Therefore, $\DTr N$ is a direct summand of $L$, and hence
$\DTr N\in \rep^+\hspace{-1pt}(Q)$, a contradiction. Therefore, $M\in \Sa$.
In particular, $\Sa$ is predecessor-closed in $\Ga$. Suppose next
that $M=\tau^{-n}\hspace{-1.5pt}X$ with $n\ge 0$ and $X\in \Sa$. If $N$
is projective, then $M$ is projective, and hence $n=0$ and $X$ is projective, which contradicts that $X$ is
pseudo-projective. Thus $N$ is not projective, and hence either $N$
is pseudo-projective or $\tau N$ is defined in $\Ga$, that is, either $N\in \Sa$
or $\tau N\in \Ga$. By Lemma \ref{section}(2), $\Sa$ is a left-most section of $\Ga$. The proof of the lemma is completed.

\medskip

We are ready to describe the connected components of
$\Ga_{\rep^+\hspace{-0.5pt}(Q)}$. Such a connected component is
called {\it preprojective} if it contains some of the $P_x$ with
$x\in Q_0$. In case $Q$ is connected, by Lemma \ref{sqproj}(1),
$\Ga_{\rep^+\hspace{-0.5pt}(Q)}$ has a unique preprojective
component which we denote by ${\mathcal P}_{_{\hspace{-1pt}Q}}$.

\medskip

\begin{Theo}\label{prep-component} Let $Q$ be an infinite
connected strongly locally finite quiver. Then the preprojective
component ${\mathcal P}_{_{\hspace{-1pt}Q}}$ of
$\Ga_{\rep^+\hspace{-0.5pt}(Q)}\vspace{1pt}$ has a left-most
section genera\-ted by the $P_x$ with $x \in Q_0$, and
consequently, $\mathcal{P}_Q$ embeds in $\mathbb{N}Q^{\hskip 0.3pt
\rm op}$. Furthermore,

\vspace{-2pt}

\begin{enumerate}[$(1)$]

\item if $Q$ has no right infinite path, then ${\mathcal
P}_{_{\hspace{-1pt}Q}}$ is right stable of shape
$\mathbb{N}Q^{\hskip 0.3pt \rm op};$

\vspace{1pt}

\item if $Q$ has right infinite paths, then ${\mathcal
P}_{_{\hspace{-1pt}Q}}$ has a right-most section, and
consequently,  ${\mathcal P}_{_{\hspace{-1pt}Q}}$ contains only
finite $\tau$-orbits.

\vspace{0pt}

\end{enumerate} \end{Theo}

\noindent{\it Proof.} By Lemma \ref{sqproj}(1), the full subquiver
$P_{_{\hspace{-1pt}Q}}$ of ${\mathcal P}_{_{\hspace{-1pt}Q}}$
generated by the $P_x$ with $x \in Q_0$ is predecessor-closed in
${\mathcal P}_{_{\hspace{-1pt}Q}}$ and is isomorphic to
$Q^{\rm\,op}$. Clearly, $P_{_{\hspace{-1pt}Q}}$ meets at most once
any given $\tau$-orbit in $\Ga$.  Let $\tau^{-n}\hspace{-1pt}M \to
N$ be an arrow in ${\mathcal P}_{_{\hspace{-1pt}Q}}$ with $n \ge
0$ and $M \in P_{_{\hspace{-1pt}Q}}$. Then
$\tau^{-n}\hspace{-1pt}M$ is not pseudo-projective. By Lemma
\ref{asm}(2), $\rep^+\hspace{-1pt}(Q)$ has a minimal right almost
split morphism $f: L\to \tau^{-n}\hspace{-1pt}M$, which is minimal
right almost split in $\rep(Q)$. If $N \not\in
P_{_{\hspace{-1pt}Q}}$, then $\rep(Q)$ has an almost split
sequence $$\xymatrixcolsep{15pt}\xymatrix{0\ar[r] & \DTr N \ar[r]
& E\ar[r] & N\ar[r] & 0}.$$ In view of Corollary \ref{irr+rep}, we
see that $\rep(Q)$ admits an irreducible morphism $g:
\tau^{-n}\hspace{-1pt}M\to N$, and hence an irreducible morphism
$h: \DTr N \to \tau^{-n}\hspace{-1pt}M$. As a consequence, $\DTr
N$ is a direct summand of $L$. In particular, $\DTr N \in
\rep^+\hspace{-1pt}(Q)$, and therefore, $\tau N\in {\mathcal
P}_{_{\hspace{-1pt}Q}}$. By Lemma \ref{section}(2),
$P_{_{\hspace{-1pt}Q}}$ is a left-most section of ${\mathcal
P}_{_{\hspace{-1pt}Q}}$. In particular, ${\mathcal
P}_{_{\hspace{-1pt}Q}}$ embeds in $\mathbb{N}Q^{\hspace{0.4pt} \rm
op}$; see \cite[(2.3)]{L4}.

\vspace{1pt}

Furthermore, if $Q$ has no right infinite path, then
$\rep^+\hspace{-1pt}(Q)=\rep^b(Q)$. In parti\-cular, ${\mathcal
P}_{_{\hspace{-1pt}Q}}$ contains only finite dimensional
representation. Containing no injective representation by
Proposition \ref{orbits}(1) and Corollary \ref{cor1.14shiping}(2),
${\mathcal P}_{_{\hspace{-1pt}Q}}$ is right stable by Lemma
\ref{psproj}(2). As a consequence, ${\mathcal
P}_{_{\hspace{-1pt}Q}}\cong \mathbb{N}Q^{\hspace{0.4pt} \rm op}$.
Otherwise, some of the $P_x$ are infinite dimensional, and hence
${\mathcal P}_{_{\hspace{-1pt}Q}}$ has a right-most section by
Lemma \ref{infrepsection}. Now, since ${\mathcal
P}_{_{\hspace{-1pt}Q}}$ has a left-most section and a right-most
one, every $\tau$-orbit in ${\mathcal P}_{_{\hspace{-1pt}Q}}$ is
finite. The proof of the theorem is completed.

\medskip

\noindent{\sc Remark.} In case $Q$ has right infinite paths, we
can describe ${\mathcal P}_{_{\hspace{-1pt}Q}}$ more explicitly in
the following way. Consider the right stable translation quiver
$\mathbb{N}Q^{\hspace{0.4pt} \rm op}$. We first define $f(0,
x)={\rm dim}_k\, P_x\in \mathbb{N}\cup \{\infty\}$ for $x\in Q_0$,
and then extend this in a unique way to an additive function
\vspace{-3pt}$$f:\mathbb{N}Q^{\hspace{0.4pt} \rm op} \to
\mathbb{N}\cup \{\infty\}$$ such that $f(v)=\infty$ whenever $v$
is a successor of some $u$ for which $f(u)=\infty$. Then
${\mathcal P}_{_{\hspace{-1pt}Q}}$ is isomorphic to the full
translation subquiver of $\mathbb{N}Q^{\hspace{0.4pt} \rm op}$
generated by the vertices $(n, x)$ with $n\in \N$ and $x\in Q_0$
such that $n=0$, or otherwise, $f(n-1, x)<\infty$.

\medskip

\noindent{\sc Example.} If $Q$ is the infinite quiver

\vspace{-6pt}
$$\xymatrixrowsep{15pt}
\xymatrixcolsep{18pt}\xymatrix{\circ & \ar[l] \circ & \ar[l] \circ
\ar[r] & \circ \ar[r]& \circ \ar[r]&  \cdots,} \vspace{1pt}$$ then
the preprojective component of $\Ga_{\rep^+(Q)}$ has of the
following shape$\,:$

\vspace{-6pt}
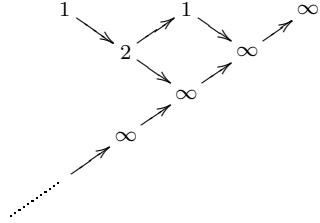
$$\xymatrixrowsep{14pt}
\xymatrixcolsep{21pt} \footnotesize \xymatrix@!C=1pt@!R=2pt{
&&& 1 \ar[rd] & & 1 \ar[rd] & & \infty\\
&&&&2 \ar[ru] \ar[rd] && {\footnotesize\infty} \ar[ru]\\
&&&&& \infty \ar[ru]\\
&&&& \tiny\infty \ar[ru]\\
&&& \ar[ru]  \\
&&\ar@{.}[ru]&& } \vspace{2pt}$$
 where each vertex is labeled with the dimension of
the corresponding representation.

\medskip

Next, we shall describe the connected components of
$\Ga_{\rep^+\hspace{-0.5pt}(Q)}$ containing some of the $I_x$ with
$x\in Q^+$, called the {\it preinjective components}. To do so,
for each $x\in Q^+$, we denote by $Q_x^+$ the connected component
of $Q^+$ containing $x$.

\medskip


\begin{Theo} \label{preinjcomp} Let $Q$ be an infinite connected strongly locally finite
quiver. If $\Ga$ is a preinjective component of
$\Ga_{\rep^+\hspace{-0.5pt}(Q)}\vspace{1pt}$ containing $I_x$ for
some $x\in Q^+$, then it has a right-most section generated by the
$I_y$ with $y\in Q_x^+$, and consequently, it contains only finite
dimensional representations and embeds in $\N^-{\hskip
-1pt}(Q_x^+)^{\rm op}$.  Furthermore,

\vspace{-2pt}

\begin{enumerate}[$(1)$]

\item if $Q$ has no left infinite path, then
$\Ga_{\rep^+\hspace{-0.5pt}(Q)}$ has a unique preinjective
component of shape $\N^-{\hskip -1pt}Q^{\hspace{0.4pt} \rm op};$

\vspace{2pt}

\item if $Q$ has left infinite paths, then every preinjective
component of $\Ga_{\rep^+\hspace{-0.5pt}(Q)}$ has a left-most
section generated by its pseudo-projective representations, and
consequently, it contains only finite $\tau$-orbits.

\end{enumerate} \end{Theo}

\noindent{\it Proof.} We fix a preinjective component $\Ga$ of
$\Ga_{\rep^+\hspace{-0.5pt}(Q)}$, which contain some $I_x$ with
$x\in Q_0$. Let ${\it \Delta}$ be the full subquiver of $\Ga$
generated by the $I_y\vspace{1pt}$ with $y \in Q_x^+$. From Lemma
\ref{sqproj}(2), we deduce that ${\it \Delta}$ is successor-closed
in $\Ga$ and isomorphic to $(Q_x^+)^{\rm\,op}$. Moreover, $\Da$
clearly meets at most once any $\tau$-orbit in $\Ga$. Let $M \to
\tau^n N$ be an arrow in $\Ga$, where $n \ge 0$ and $N \in {\it
\Delta}$. Then $N = I_y$ for some $y \in Q_x^+$. Since $N$ is finite
dimensional, it follows from Lemma \ref{infdimirr} that $M$ is
finite dimensional. If $M=I_z$ for some $z \in Q^+$, then
$\tau^n\hspace{-1pt}I_y$ is injective by Lemma \ref{sqproj}(2).
Therefore, $n=0$, and hence $z \in Q_x^+$. That is, $M\in \Da$.
Otherwise, by Lemma \ref{psproj}(2), $\tau^-\hspace{-1pt}M$ is
defined in $\Ga$. By Lemma \ref{section}, ${\it \Delta}$ is a
right-most section of $\Ga$. Since ${\it \Delta}$ contains only
finite dimensional representations, by Lemma \ref{infdimirr}, every
representation in $\Ga$ is finite dimensional.

If $Q$ contains no left infinite path, then $Q=Q^+$ and ${\rm
rep}^-(Q)={\rm rep}^b(Q)$. Since $\Ga$ contains no projective
representation by Proposition \ref{orbits}(1), we see from Lemma
\ref{psproj}(1) that $\tau$ is defined everywhere in $\Ga$, that
is, $\Ga$ is left stable. As a consequence, $\Ga\cong
\N^-(Q_x^+)^{\rm op}$. On the other hand, since $Q$ is connected,
$Q^+_x=Q$. Thus, $\Ga$ contains all the $I_y$ with $y\in Q_0$. In
particular, $\Ga$ is the unique preinjective component, which is
of shape $\N^- Q^{\hspace{0.4pt} \rm op}$.

Finally, suppose that $Q$ contains left infinite paths. Since
$Q^+$ is predecessor-closed in $Q$ by definition, $Q$ has some
arrow $y \to z$ with $y \in Q_x^+$ and $z\not\in Q^+$. Then
$\rep(Q)$ has an irreducible morphism $f: I_z \to I_y$ with $I_z$
infinite dimensional. Since $I_y\not\in {\rm proj}(Q)$ by
Proposition \ref{orbits}(1), $\rep(Q)$ has an almost split
sequence $\xymatrixcolsep{15pt}\xymatrix{0\ar[r] & \DTr I_y\ar[r]
& E\ar[r] & I_y\ar[r] & 0.}$ Thus $I_z$ is an infinite dimensional
direct summand of $E$. Since $I_y$ is finite dimensional, $\DTr
I_y$ is infinite dimensional, that is, $I_y$ is a
pseudo-projective representation. By Lemma \ref{infrepsection},
the pseudo-projective representations in $\Ga$ generate a
left-most section.  The proof of the theorem is completed.

\medskip

\noindent {\sc Remark.} (1) Theorem \ref{preinjcomp} says that the
preinjective components of $\Ga_{\rep^+\hspace{-0.5pt}(Q)}$
correspond bijectively to the connected components of $Q^+$. In
particular, $\Ga_{\rep^+\hspace{-0.5pt}(Q)}$  has no preinjective
component if $Q^+$ is empty.

\vspace{1pt}

(2) In case $Q$ has left infinite paths, the preinjective
components can be found in the following way. Consider the left
stable translation quiver
$\mathbb{N}^-\hspace{-1.5pt}Q^{\hspace{0.4pt} \rm op}$. We define
$f(0, x)={\rm dim}_k \,I_x\in \mathbb{N} \cup \{\infty\}$ for
$x\in Q_0$, and extend this in a unique way to an additive
function \vspace{-4pt}
$$f:
\mathbb{N}^-\hspace{-1.5pt}Q^{\hspace{0.4pt} \rm op} \to
\mathbb{N} \cup \{\infty\}$$ such that $f(v)=\infty$ if $v$ is a
predecessor of some vertex $u$ with $f(u)=\infty$. Then the
preinjective components of $\Ga_{\rep^+(Q)}$ correspond
bijectively to the connected components of the full translation
subquiver of $\mathbb{N}^-\hspace{-1.5pt}Q^{\hspace{0.4pt} \rm
op}$ generated by the vertices $(n, x)$ with $f(n, x) < \infty$.

\medskip

\noindent{\sc Example.} If $Q$ is the infinite quiver
$$\footnotesize
\xymatrixrowsep{18pt} \xymatrixcolsep{25pt}\xymatrix@!=1pt{&&&& \;0 \ar[d]&&&&\\
&\ar@{.}[l] \ar[r] & -2  \ar[r] & -1\ar[r]&1&2\ar[l]& 3\ar[l]\ar[r] &4 \ar[r] &5&6\ar[l]& \ar@{.}[r]\ar[l]&} \vspace{6pt}$$

\noindent then $\Ga_{\rep^+(Q)}$ has a trivial preinjective
component $\{I_0\}$ and another preinjective component of the
following shape$\,:$
$$\footnotesize\xymatrixrowsep{11pt}
\xymatrixcolsep{25pt} \xymatrix@!C=1pt@!R=2pt{& I_2\ar[dr] &\\ \tau
I_3 \ar[ur]\ar[dr] & & I_3\\ & I_4 \ar[ur]& } \vspace{1pt}$$

\medskip

A representation lying in $\Ga_{\rep^+(Q)}$ is called {\it
preprojective} or {\it preinjective} if it lies in a preprojective
component or in a preinjective component, respectively. Before
going further, we shall study some properties of these
representations.

\medskip

\begin{Lemma}\label{fin-prec-suc} Let $M$ be a representation in $\Ga_{\rep^+\hspace{-0.5pt}(Q)}\vspace{1pt}$.

\vspace{-2pt}

\begin{enumerate}[$(1)$]

\item If $M$ is preprojective, then it has only finitely many non-projective predecessors in the
preprojective component.

\vspace{1pt}

\item If $M$ is preinjective, then it has only finitely many successors in its preinjective component.

\end{enumerate} \end{Lemma}

\noindent{\it Proof.} We may assume that $Q$ is connected. Suppose
first that $M$ lies in the preprojective component
$\mathcal{P}_{_{\hspace{-1pt}Q}}$. By Theorem
\ref{prep-component} and Lemma \ref{sqproj},
$\mathcal{P}_{_{\hspace{-1pt}Q}}$ has a left-most section
$P_{_{\hspace{-1pt}Q}}$, which is generated by the $P_x$ with
$x\in Q_0$ and isomorphic to $Q^{\hspace{0.4pt} \rm op}$. Suppose
that $M$ has infinitely many non-projective predecessors in
$\mathcal{P}_{_{\hspace{-1pt}Q}}$. By K\"{o}nig's lemma,
$\mathcal{P}_{_{\hspace{-1pt}Q}}$ has a left infinite path of
non-projective representations as follows$\,:$ \vspace{-1pt}
$$\xymatrixcolsep{15pt}\xymatrix{\cdots \ar[r] & M_i\ar[r] & M_{i-1}\ar[r] & \cdots \ar[r] & M_1\ar[r] &
M_0=M}.$$ \vspace{-2 pt}Since $P_{_{\hspace{-1pt}Q}}$ is a left-most section, we
can write $M_i=\tau^{-r_i}\hspace{-1.5pt}P_{x_i}$, where $x_i\in
Q_0$ and $r_i\in \N$ such that $r_i\ge r_{i+1}>0$ for all $i\ge 0$.
Thus, we may assume that $r_i=r_0$ for all $i\ge 1$. This implies
that $P_{_{\hspace{-1pt}Q}}$ has a left infinite path ending in $P_{x_0}$, which
in turn implies that $Q$ has a right infinite path starting in
$x_0$. In particular, $P_{x_0}$ is infinite dimensional. By Lemma
\ref{psproj}(2), $\tau^-P_{x_0}$ is not defined, which is absurd
since $r_0>0$. This proves Statement (1). Since the preinjective
components of $\Ga_{\rep^+\hspace{-0.5pt}(Q)}$ contain only finite dimensional
representations, we can prove Statement (2) in a dual manner. The
proof of the lemma is completed.

\medskip

\begin{Lemma}\label{prec-suc} Let $f: M\to N$ be a non-zero non-invertible morphism in $\rep^+\hspace{-1pt}(Q)$ with
$M, N\in \Ga_{\rep^+\hspace{-0.5pt}(Q)}\vspace{1pt}$. If $M$ is preinjective or $N$ is
preprojective, then $\Ga_{\rep^+\hspace{-0.5pt}(Q)}$ contains a non-trivial path
from $M$ to $N$.
\end{Lemma}

\noindent{\it Proof.} We consider only the case where $N$ lies in
the preprojective component $\mathcal{P}_{_{\hspace{-1pt}Q}}$ of
$\Ga_{\rep^+\hspace{-0.5pt}(Q)}$, since the other case can be treated
in a dual manner. Let $P_{_{\hspace{-1pt}Q}}$ be the left-most section in
$\mathcal{P}_{_{\hspace{-1pt}Q}}$. Suppose first that $N=P_y$ for some $y\in
Q_0$. Since $f$ is non-zero and $\rep^+\hspace{-1pt}(Q)$ is
hereditary, $M=P_x$ for some $x\in Q_0$. Since $f$ is
non-invertible, we deduce from Proposition \ref{prop1.3_shiping}
that $x$ is a proper successor of $y$ in $Q$. Hence $P_y$ is a
proper successor of $P_x$ in $P_{_{\hspace{-1pt}Q}}$. Suppose now that $N$ is not
projective while $M$ is not a proper predecessor of $N$ in
${\mathcal P}_{_{\hspace{-1pt}Q}}$. Since every representation in ${\mathcal P}_{_{\hspace{-1pt}Q}}$ is the co-domain of a minimal right almost split morphism
in $\rep^+\hspace{-1pt}(Q)$, using induction, we get an infinite
path \vspace{-2pt}
$$\xymatrixcolsep{15pt}\xymatrix{\cdots \ar[r]& N_i \ar[r]& N_{i-1} \ar[r]& \cdots \ar[r]& N_1 \ar[r]& N_0=N}$$ in $\mathcal{P}_{_{\hspace{-1pt}Q}}$
and non-zero non-invertible morphisms $f_i: M\to N_i$ for $i\ge 0$.
By lemma \ref{fin-prec-suc}, there exists a positive integer $n$
such that $N_n$ is projective. Hence, $M$ is a proper predecessor of
$N_n$ by our previous consideration, and hence a proper predecessor
of $N$ in ${\mathcal P}_{_{\hspace{-1pt}Q}}$, a contradiction. The proof of the
lemma is completed.

\medskip

\begin{Prop}\label{no_ext} Let $M$ be an indecomposable representation lying in $\Ga_{\rep^+\hspace{-0.5pt}(Q)}$.

\vspace{-1pt}

\begin{enumerate}[$(1)$]

\item If $M$ is preprojective, then $\Hom(L, M)=0$ for all but finitely many non-projective representations $L$ in $\Ga_{\rep^+\hspace{-0.5pt}(Q)}$.

\vspace{1pt}

\item If $M$ is preinjective, then $\Hom(M, L)=0$ for all but finitely many representations $L$ in $\Ga_{\rep^+\hspace{-0.5pt}(Q)}$.

\vspace{1pt}

\item If $M$ is preprojective or preinjective, then $\Ext^1(M,M) = 0$.

\end{enumerate} \end{Prop}

\noindent{\it Proof.} The first two statements follow immediately
from Lemmas \ref{fin-prec-suc} \vspace{0pt} and \ref{prec-suc}.
Suppose that $\rep^+(Q)$ admits a non-trivial extension
$\xymatrixcolsep{18pt}\xymatrix{0\ar[r] & M\ar[r] &  E \ar[r] &
M\ar[r] & 0.}$ Then, $M$ is neither projective nor injective. Assume
that $M$ lies in the preprojective component
$\mathcal{P}_{_{\hspace{-1pt}Q}}$. If $M$ is pseudo-projective then,
by Lemma \ref{infrepsection}(2), $\mathcal{P}_{_{\hspace{-1pt}Q}}$
has a left-most section generated by its pseudo-projective
representations, which coincides with the left-most section
generated by the $P_x$ with $x\in Q_0$, a contradiction. Hence $\tau
M$ is defined in $\mathcal{P}_{_{\hspace{-1pt}Q}}$. This yields a
commutative diagram \vspace{-2pt}
$$\xymatrixrowsep{15pt}\xymatrixcolsep{18pt}\xymatrix{
0\ar[r] & M\ar[r]\ar[d]^{f}& E \ar[r]\ar[d] & M\ar[r] \ar@{=}[d]& 0\\
0\ar[r] & \tau M\ar[r] & N \ar[r] & M\ar[r]& 0}\vspace{-1pt}$$ in $\rep^+(Q)$, where the lower row is an almost split sequence. By Lemma \ref{prec-suc},
$f=0,$ and hence the lower row splits, a contradiction.
Suppose next that $M$ lies in a preinjective component $\mathcal I$. Observing that $\tau^-\hspace{-1pt}M$ is defined in $\mathcal{I}$, we get a contradiction by  a dual argument. The proof of the proposition is completed.

\medskip

The rest of this section is  mainly devoted to describing the {\it regular components} of $\Ga_{\rep^+\hspace{-0.5pt}(Q)}\vspace{1pt}$, that is, the connected components which contain none of the $P_x$ and the $I_x$
with $x\in Q_0$. For this purpose, we shall need a proposition to dualize results on $\rep^+(Q)$ to those on $\rep^-(Q)$. For its proof, the following easy result is useful.

\medskip

\begin{Lemma} \label{Irr_b}
Let $f: M \to N$ be a morphism in $\rep^b(Q)$. Then $f$ is
irreducible in $\rep^b(Q)$ if and only if it is irreducible in
$\rep(Q)$.
\end{Lemma}

\noindent {\it Proof.} We only need to prove the necessity. Suppose that
$f=hg$, where $g:M \to L$ and $h: L \to N$ are morphisms in $\rep(Q)$. Let
$\Sa$ be a finite full subquiver of $Q$, containing the vertices in ${\rm supp}(M \oplus N)$ as well as their
immediate predecessors and immediate successors in $Q$. This yields a factorization $f = h_{_{\it \Sigma}}g_{_{\it \Sigma}}$ in
$\rep^b(Q)$. Therefore, $g_{_{\it \Sigma}}$ is a section or
$h_{_{\it \Sigma}}$ is a retraction. In the first case, $g$ is a section since $\Sa$ contains the immediate predecessors of the vertices in ${\rm supp}\hspace{0.4pt}M$. In the second case, $h$ is a retraction  since $\Sa$ contains the immediate successors of the vertices in ${\rm supp}\hspace{0.4pt}N$. The proof of the lemma is completed.

\medskip

\begin{Prop} \label{dualityregularcomp} Let $\Ga$ be a regular component of $\Ga_{\rep^+\hspace{-0.5pt}(Q)}$.

\vspace{-1pt}

\begin{enumerate}[$(1)$]

\item If $\Ga$ has infinite dimensional but no pseudo-projective representations, then the full translation subquiver of $\Ga$ obtained by deleting the infinite dimensional representations is a left stable regular component of $\Ga_{\rep^-(Q)}$.

\vspace{1pt}

\item If $\Ga$ has pseudo-projective but no infinite dimensional representations, then $\Ga$ is the full translation subquiver of
a right stable regular component of $\Ga_{\rep^-(Q)}$ obtained by deleting the non-empty set of infinite dimensional representations.

\end{enumerate}\end{Prop}

\noindent{\it Proof.} (1) Suppose that $\Ga$ contains infinite
dimensional but no pseudo-projective representations. By Proposition
\ref{infrepsection}(1), the infinite dimensional representations in
$\Ga$ generate a right-most section $\Da$. Let $\Ga\hspace{0.2pt}'$ be the full
translation subquiver of $\Ga$ gene\-rated by the representations
not in $\Da$. Containing no projective or pseudo-projective representation, $\Ga\hspace{0.2pt}'$ is left stable.
Being finite dimensional, the representations in $\Ga\hspace{0.2pt}'$ may be assumed to all lie in
$\Ga_{\rep^-(Q)}$. Fix representations $M, N$ in $\Ga\hspace{0.2pt}'$.
It follows from Lemma \ref{Irr_b} that a morphism $f: M \to N^r$ with $r>0$ is irreducible in $\rep^+(Q)$ if and only if it is irreducible in $\rep^-(Q)$ and a morphism $g: M^s \to N$ with $s>0$ is irreducible in $\rep^+(Q)$ if and only if it is irreducible in $\rep^-(Q)$. Therefore,
$M\to N$ is a valued arrow with valuation $(d, d\,')$ in $\Ga\hspace{0.2pt}'$ if and only if $M\to N$ is a valued arrow with valuation $(d, d\,')$ in $\Ga_{\rep^-(Q)}.\vspace{1pt}$ In particular, $\Ga\hspace{0.2pt}'$ is a full valued subquiver of some connected
component $\mathcal C$ of $\Ga_{\rep^-(Q)}$. Next, since $M$ is neither projective nor
pseudo-projective, $\rep^+\hspace{-1pt}(Q)$ has an almost split
sequence $\eta: \xymatrixcolsep{15pt}\xymatrix{0\ar[r] & \tau M\ar[r] & E\ar[r] & M\ar[r] & 0,}$
where $\tau M$ is finite dimensional. By Lemma \ref{Irr_b}, $\eta$ is also an
almost split sequence in $\rep^-(Q)$. This shows that $\Ga\hspace{0.2pt}'$ is a
predecessor-closed valued translation subquiver of $\mathcal C$. Next,
let $M\to X$ be an arrow in $\mathcal C$ with an irreducible
morphism $u: M \to X$ in $\rep^-\hspace{-1pt}(Q)$. By the dual of
Corollary \ref{irr+rep}, $h$ is irreducible in $\rep(Q)$. On the other hand, since $M$
is finite dimensional and not injective, by Corollary
\ref{cor_ARS_exist}(2), $\rep^+\hspace{-1pt}(Q)$ has an almost split
sequence $\xymatrixcolsep{15pt}\xymatrix{0 \ar[r] & M \ar[r] & E
\ar[r] & \tau^-\hspace{-1pt}M \ar[r] & 0}$, which is also almost split in
$\rep(Q)$. Then $X$ is a direct summand of $E$. If $\tau^-\hspace{-1pt} M\in \Ga\hspace{0.2pt}'$, then
$E$ is finite dimensional and so is $X$. If $\tau^-\hspace{-1pt} M\in\Da$ then, by Corollary
\ref{Cor3.3}(2), $X$ is finite dimensional. In any case, $M\to X$ is an arrow in $\Ga\hspace{0.2pt}'$. This shows that
$\Ga\hspace{0.2pt}'$ is successor-closed in $\mathcal C$, and hence,
$\Ga\hspace{0.2pt}'={\mathcal C}$. That is,  $\Ga\hspace{0.2pt}'$ is a regular component of
$\Ga_{\rep^-(Q)}$ which is left stable.

\vspace{1pt}

(2) Suppose that $\Ga$ contains pseudo-projective but no
infinite dimensional representations. Since
the representations in $\Ga$ are finite dimensional and non-injective, $\Ga$
is right stable. Using an argument dual to the above one, we see
that $\Ga$ is a successor-closed valued translation subquiver of a connected
component $\mathcal C$ of $\Ga_{\rep^-(Q)}$. Since $\Ga$
is right stable, $\mathcal C$ has no right-most section. In particular, by the dual of Theorem
\ref{prep-component} and the dual of Lemma \ref{infrepsection}(2), $\mathcal C$ contains no representation $M$ which is injective or
pseudo-injective, where $M$ is {\it pseudo-injective} if $\TrD M$ is infinite dimensional.
Now, fix a pseudo-projective representation $N$ in $\Ga$. Since it is finite dimensional and non-projective,  by
Corollary \ref{cor_ARS_exist}(1), $\rep(Q)$ has an almost split sequence $\xymatrixcolsep{15pt}\xymatrix{0\ar[r] &
\DTr N\ar[r] & E\ar[r] & N\ar[r] & 0}$ with $\DTr N$ of infinite
dimension, which is also an almost split sequence in $\rep^-(Q)$. In particular, $\DTr N$ is an
infinite dimensional representation in $\mathcal C$. By the dual
of Theorem \ref{preinjcomp}, $\mathcal C$ is not a preprojective
component, and hence a regular component of $\Ga_{\rep^-(Q)}$. By
the dual of Statement (1), the full translation subquiver $\mathcal C'$ of $\mathcal C$ obtained by deleting the infinite
dimensional representations is a connected component of
$\Ga_{\rep^+\hspace{-0.5pt}(Q)}$. Since $\Ga$ is a connected component of $\Ga_{\rep^+\hspace{-0.5pt}(Q)}$ which is contained in $\mathcal C'$,
we see that $\Ga$ and $\mathcal C'$ coincide. The proof of the proposition is
completed.

\medskip

We shall also need the following result to deal with the regular components containing infinite dimensional representations.

\medskip

\begin{Lemma}\label{infdimrep} Let $\Ga$ be a regular component of $\Ga_{\rep^+\hspace{-0.5pt}(Q)}$, and let $M$ be an infinite dimensional representation lying in $\Ga$.

\vspace{-1pt}

\begin{enumerate}[$(1)$]

\item If $N, L$ are representations in $\Ga$, then $\rep(Q) \vspace{-3pt}$ admits no chain of irreducible monomorphisms $\xymatrixcolsep{15pt}\xymatrix{\DTr L \ar[r]^-f & N\ar[r]^g & L}$.

\vspace{1pt}

\item If $M\to N$ is an arrow in $\Ga$, then $\tau N\in \Ga$ with ${\rm dim}_k \,\DTr M > {\rm dim}_k \,\tau N$,
 and 
 any irreducible morphism $f: M \to N$ in $\rep^+\hspace{-1pt}(Q)$ is an epimorphism.

\vspace{1pt}

\item If $\xymatrixcolsep{15pt}\xymatrix{M \ar[r] & M_1 \ar[r] & \cdots \ar[r] & M_{n-1} \ar[r] & M_n}$ is a path in $\Ga$, then
$\tau^jM_i\in \Ga$ for all $i=1, \ldots, n; \, j=0,\ldots, i$.

\vspace{1pt}

\item If $\xymatrixcolsep{15pt}\xymatrix{M\ar[r] & M_1 \ar[r] & M_2 \ar[r] & M_3}$ is a path in $\Ga$ and $f:  M_3\to N$ is
an irreducible morphism in $\rep^+\hspace{-1pt}(Q)$, then $N$ is
indecomposable.

\vspace{1pt}

\item If $M$ is not pseudo-projective, then $\rep^+(Q)$ has a minimal right almost split morphism $f: N_1\oplus N_2 \to M$, where $N_1$ is
indecomposable of infinite dimension and $N_2$ is of finite
dimension.

\end{enumerate} \end{Lemma}

\noindent{\it Proof.} Let $X$ be a representation lying in $\Ga$.
Write $d(X)={\rm dim}_k\,X\in \N\cup\{\infty\}$. Since $\Ga$ is regular, $\DTr
X$ is an indecomposable representation in $\rep^-(Q)$. Thus, $\tau X$
is defined in $\Ga$ if and only if $X$ is not pseudo-projective.

\vspace{-1pt}

(1) Suppose that $\rep(Q)$ admits irreducible monomorphisms
$\xymatrixcolsep{15pt}\xymatrix{\DTr L \ar[r]^-f & N \ar[r]^g & L}$
where $N, L$ in $\Ga$. Since $\DTr L$ is a finitely co-presented
sub-representation of $N$, making use of Lemma
\ref{Lemma1.6_shiping}(2) and Corollary \ref{newcor1.9}, we see that
$\DTr L$ is finite dimensional, and hence $\tau L \vspace{-1pt}$ is
defined in $\Ga$. Now $\rep(Q)$ has an irreducible morphism $g_1:
\DTr N \to \DTr L,$ which is a monomorphism by Corollary
\ref{preservemono}. In particular, $\DTr N$ is finite dimensional,
and hence $\tau N\vspace{-1pt}$ is defined in $\Ga$. Applying the
same argument to $\xymatrixcolsep{15pt}\xymatrix{\tau N\ar[r]^{g_1}
& \tau L \ar[r]^f & N,}$ we get an irreducible monomorphism $f_1:
\DTr^2\hspace{-1pt}L\to \tau N$ in $\rep(Q)$. Repeating this
process, we see that $\tau^i\hspace{-1pt}N$ and
$\tau^i\hspace{-1pt}L$ are defined in $\Ga$ for all $i\ge 0$, and
$\rep(Q)$ admits an infinite chain of irreducible monomorphisms
\vspace{-5pt}
$$\xymatrixcolsep{18pt}\xymatrix{\cdots \ar[r] & \tau^3L \ar[r]^{f_2}
& \tau^2N \ar[r]^{g_2} & \tau^2L \ar[r]^{f_1} & \tau N \ar[r]^{g_1}
& \tau L \ar[r]^{f} & N \ar[r]^{g} & L,}\vspace{-2pt}$$ which is
absurd since $d(\tau L)<\infty.$

\vspace{1pt}

(2) Let $M\to N$ be an arrow in $\Ga$. By Lemma \ref{infdimirr},
$\DTr N$ is finite dimensional, and hence $\tau N\in \Ga$, and
$\rep(Q)$ has an irreducible monomorphism $g: \tau N\to M$. Since
$M$ is not projective, $\rep(Q)$ has an irreducible morphism $h:
\DTr M\to \tau N$, which is an epimorphism by Statement (1). Hence
${\rm dim}_k \,\DTr M > d(\tau N).$ By Corollary \ref{preservemono},
every irreducible morphism $f: M \to N$ in $\rep^+\hspace{-1pt}(Q)$
is an epimorphism.

\vspace{1pt}

(3) Let $\xymatrixcolsep{15pt}\xymatrix{M \ar[r] & M_1 \ar[r] &
\cdots \ar[r] & M_{n-1} \ar[r] & M_n}$ be a path in $\Ga$. By Lemma
\ref{infdimirr}, $M_i$ is infinite dimensional while $\tau M_i$ is defined of
finite dimensional, for every $1\le i\le n$. Thus, $\rep^+(Q)$ has irreducible
monomorphisms $f_{i,i-1}: \tau M_i \to M_{i-1}$, $i=1, \ldots, n$,
where $M_0=M$. Let $i$ with $0<i\le n$ be such that
$\rep^+\hspace{-1pt}(Q)$ has a chain of monomorphisms \vspace{-4pt}
$$\xymatrixrowsep{14pt}
\xymatrixcolsep{22pt}\xymatrix{ \tau^{n-i}M_n\ar[r]^-{f_{n,
i}}&\tau^{n-i-1}M_{n-1}\ar[r] & \cdots \ar[r]&  \tau
M_{i+1}\ar[r]^-{f_{i+1,i}}&M_i.} \vspace{-2pt}$$ By Corollary \ref{preservemono},
$\rep^-(Q)$ has a chain of monomorphisms
\vspace{-4pt}$$\xymatrixrowsep{18pt} \xymatrixcolsep{21pt}\xymatrix{
\DTr^{n-i+1}M_n\ar[r]^-{f_{n, i-1}}&\DTr^{n-i}M_{n-1}\ar[r] & \cdots
\ar[r]&  \DTr^2 M_{i+1}\ar[r]^-{f_{i+1,i-1}}&\tau M_i \ar[r]^-{f_{i,
i-1}} & M_{i-1}. } \vspace{-3pt}$$ Since $\tau M_i$ is finite
dimensional, the $\DTr^{j-i+1}M_j$ with $i\le j\le n$ are all finite
dimensional. That is, $\tau^{j-i+1}M_j$ is defined in $\Ga,$ for all $i\le j\le
n$.

\vspace{1pt}

(4) Let $\xymatrixcolsep{15pt}\xymatrix{M \ar[r] & M_1 \ar[r] & M_2
\ar[r] & M_3}$ be a path in $\Ga$, and let $f: M_3 \to N$ be an
irreducible morphism in $\rep^+\hspace{-1pt}(Q)$. Write
$N=\oplus_{i=1}^n\, N_i$ with the $N_i$ indecomposable. By
Statement (3), $\tau^jM_i$ is defined of
finite dimension for $1\le i\le 3$ and $1\le j\le i$, and  $\tau^j N_i$ is defined of finite dimension for
$1\le i \le n$ and $1\le j\le 4$. Since  $M_2$ is infinite dimensional by Lemma \ref{infdimirr}, $\rep^+\hspace{-1pt}(Q)$ has an irreducible
monomorphism $g: \tau M_3 \to M_2$, and by Statement (2), $d(\tau
M_3)<d(\tau M_2)$. Moreover, it follows from Corollary \ref{preservemono} that
$d(\tau^3M_3) < d(\tau^2M_2)$. This yields

\vspace{-6pt}$$  d(\tau^3M_3) + d(\tau M_3)<
d(\tau^2M_2) + d(\tau M_3)< d(\tau^2M_2) + d(\tau M_2).$$
Suppose that
$n\ge 2$. Since $d$ is additive,  $\,d(\tau^3 N_i) + d(\tau^2 N_i)
\ge d(\tau^2 M_3)$ for $i=1,2,$ and $d(\tau^2 M_3) + d(\tau M_3) \ge
d(\tau^2 N_1) + d(\tau^2 N_2) + d(\tau M_2).$ Furthermore,
\vspace{-1.5pt}
$$\begin{array}{rcl}
d(\tau^3M_3) + d(\tau^2M_3) & \ge& d(\tau^3N_1) + d(\tau^3N_2) + d(\tau^2M_2)\\
&\ge &d(\tau^2 M_3)-d(\tau^2 N_1) + d(\tau^2 M_3)-d(\tau^2 N_2) + d(\tau^2M_2)\\
&\ge & d(\tau^2 M_3) + d(\tau M_2)- d(\tau M_3) + d(\tau^2M_2). \vspace{-1.5pt}
\end{array} $$
As a consequence, we get $d(\tau^3M_3) + d(\tau M_3) \ge d(\tau^2M_2)+d(\tau
M_2)$, a contradiction. Thus $N$ is indecomposable.

\vspace{1pt}

(5) Let $M$ be not pseudo-projective. Then $\rep^+(Q)$ has an almost
split sequence \vspace{-5pt}
$$\xymatrixrowsep{15pt} \xymatrixcolsep{18pt}\xymatrix{0\ar[r] &
\tau M \; \ar[rr]^-{(g_1, \ldots, g_n)^T} && \;
N_1\oplus\cdots\oplus N_n \;\;
\ar[rr]^-{(f_1, \ldots, f_n)} & & \;\; M \ar[r] & 0,}$$ where $\tau
M$ is finite dimensional and the $N_i$ are indecomposable. Since
$d(M)=\infty$, we may assume that $d(N_1)=\infty$. Then $g_1: \tau
M\to N_1$ is a monomorphism. By Statement (1), $f_1$ is an
epimorphism. Hence, $g_i: \tau M\to N_i$ is an epimorphism, and
hence $d(N_i)<\infty$, for $1<i\le n$. The proof of the lemma is
completed.

\medskip

Finally, we recall that a valued translation quiver is called a {\it wing} if it is isomorphic to the following trivially valued translation quiver$\,:$
\vspace{-1pt} $$\xymatrixrowsep{16pt}
\xymatrixcolsep{20pt}
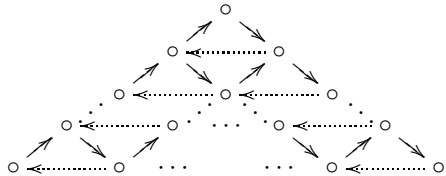
\xymatrix@!=0.1pt{&&&& \circ \ar[dr]&&&&\\
&&&\circ \ar[ur]\ar[dr] \ar@{<.}[rr]&& \circ \ar[dr]&&&\\
&& \circ\ar[ur] \ar@{<.}[rr] && \circ\ar[ur] \ar@{<.}[rr] && \circ
&&\\
}$$ \vspace{-16 pt}
$$\iddots \quad\;\; \quad \; \iddots \;\;\;  \ddots \;\quad \;\;\;\;\; \ddots$$
\vspace{-24 pt}
$$\xymatrixrowsep{16pt}
\xymatrixcolsep{20pt}\xymatrix@!=0.1pt{& \circ\ar[dr] \ar@{<.}[rr] && \circ   &\cdots& \circ\ar[dr]  \ar@{<.}[rr]  && \circ \ar[dr] & \\
\circ \ar[ur]\ar@{<.}[rr] && \circ \ar[ur]  &\cdots&&\cdots&
\circ\ar[ur] \ar@{<.}[rr]
 && \circ}\vspace{6pt}
$$
where the dotted arrows indicate the translation; see \cite[(3.3)]{R3}.

\medskip

Now we have the promised description of the regular components of $\Ga_{\rep^+\hspace{-0.5pt}(Q)}$.

\medskip

\begin{Theo}\label{reg-cpt} Let $Q$ be an infinite connected strongly locally finite quiver, and let $\Ga$ be a regular component of $\Ga_{\rep^+\hspace{-0.5pt}(Q)}$.

\vspace{-1pt}

\begin{enumerate}[$(1)$]

\item If $\Ga$ has no infinite dimensional or pseudo-projective representation, then it is of shape $\mathbb{Z} \mathbb{A}_\infty$.

\vspace{1pt}

\item If $\Ga$ has infinite dimensional but no pseudo-projective representations, then it is of shape $\mathbb{N}^-\hskip -0.5pt \mathbb{A}_\infty$ and its right-most section is a left infinite path.

\vspace{1pt}

\item If $\Ga$ has pseudo-projective but no infinite dimensional representations, then it is of shape $\mathbb{N}\mathbb{A}_\infty$ and its left-most section is a right infinite path.

\vspace{1pt}

\item  If $\Ga$ has both pseudo-projective representations and infinite dimensional representations, then $\Ga$ is a finite wing.

\end{enumerate}
\end{Theo}

\noindent{\it Proof.} (1) Write $d(M)={\rm dim}_k\,M\in \N\cup\{\infty\}$, for $M\in
\Ga$. Suppose that $\Ga$ contains no infinite dimensional or
pseudo-projective representation. By Lemma \ref{psproj}, $\Ga$ is
stable. Having no oriented cycle by Lemma \ref{nocycle}, $\Ga$ is
isomorphic to $\mathbb{Z}\Da$, where $\Da$ is a section of $\Ga$;
see \cite[(2.3)]{L4}. Consider the additive function
$$d: \Ga_0\to \mathbb{N}: M\mapsto d(M),$$ which is strictly monotone by Lemma \ref{preservemono},
and consequently, $\Da$ is either finite or of type
$\mathbb{A}_\infty$; see \cite{R}. Suppose that $\Da$ is finite. Let
$\it\Theta$ be the full translation subquiver of $\Ga$ generated by
the representations lying in ${\it \Delta} \cup \tau\Da \cup
\tau^-\hskip -2pt \Da$. Being connected and infinite, $Q$ has a finite connected
full subquiver $\Sa$, which contains the supports of the
representations lying in $\it \Theta$ and which has more vertices than
$\it\Delta$ does. Then $\it \Theta$ is a full translation subquiver
of some connected component $\Ga\hspace{0.2pt}'$ of $\Ga_{\rep(\it\Sigma)}$. Since
$\it \Theta$ is finite and satisfies the condition $S4$ stated in
\cite[(3.1)]{L4}, $\Ga\hspace{0.2pt}'$ is the preprojective or preinjective
component of $\Ga_{\rep(\it\Sigma)}$ having $\Da$ as a section. In
particular, $\Sa$ and $\Da$ have the same number of vertices, a
contradiction. Thus, $\Da$ is of type $\mathbb{A}_\infty$.

\vspace{1pt}

(2) Assume that $\Ga$ has infinite dimensional but no pseudo-projective representations. By Lemma \ref{psproj}(1), $\Ga$ is left stable, and by Lemma \ref{infrepsection}(1), the infinite dimensional representations
lying in $\Ga$ generate a right-most section $\Da$. Therefore,
$\Ga\cong \mathbb{N}^-\hspace{-2pt}\Da$; see \cite[(2.3)]{L4}.
By Lemma \ref{infdimrep}(2), $\Da$ contains no right infinite path, and hence it has
a sink-vertex $M_0$. Having no pseudo-projective representation, by Lemma \ref{infdimrep}(5), $\Da$ contains a left infinite
path \vspace{-0pt}$$(*) \qquad \xymatrixrowsep{18pt}
\xymatrixcolsep{15pt}\xymatrix{\cdots \ar[r] &  M_n \ar[r]  &
M_{n-1} \ar[r] & \cdots \ar[r] & M_1 \ar[r] & M_0.}\vspace{-2pt}$$
For each $n>0$, denote by $(d_n, d_n')$ the valuation of the arrow
$M_n\to M_{n-1}$. By Lemma \ref{infdimrep}(5), $d\,_n'=1$ and $M_n$ is
the only immediate predecessor of $M_{n-1}$ in $\Da$. By Lemma
\ref{infdimrep}(4), $d_n=1$ and $M_{n-1}$ is the only immediate
successor of $M_n$ in $\Da$. This shows that the path $(*)$ is trivially valued and coincides with $\Da$. In particular, $\Ga$ is of
shape $\mathbb{N}^-\hskip -2pt \A_\infty.$

\vspace{1pt}

(3) Assume that $\Ga$ has pseudo-projective but no
infinite dimensional representations. By Proposition
\ref{dualityregularcomp}(2), $\Ga$ is the full translation subquiver
of a right stable regular component $\Ga\hspace{0.2pt}'$ of
$\Ga_{\rep^-(Q)}$ obtained by deleting the infinite dimensional
representations. By the dual of Statement (2), the infinite
dimensional representations in $\Ga\hspace{0.2pt}'$ generate a left-most
section which is a right infinite path. As a consequence, the
pseudo-projective representations in $\Ga$ generate a left-most
section which is a right infinite path. In particular, $\Ga$ is of
shape $\N\A_\infty$.

\vspace{1pt}

(4) Suppose that $\Ga$ has both infinite dimensional
representations and pseudo-projective ones. By Lemma
\ref{infrepsection}, $\Ga$ has a right-most section $\Da$ generated by the infinite dimensional
representations, and a left-most section generated by the pseudo-projective ones. In particular, $\Ga$ contains no left or right stable representation.

Next, we show that $\Da$ is a finite trivially valued path.
If $\Da$ has no pseudo-projective representation, then it follows from Lemma
\ref{infdimrep}(3),(5) that $\Da$ has a left infinite path in which every representation is left stable, a contradiction.
Thus $\Da$ contains a pseudo-projective representation $M_0$. By Lemma \ref{infdimrep}(2),(5), $M_0$ is a unique source vertex in $\Da$. If $N\in \Da$ with $N\ne M_0$, then $N$ is a proper successor of $M_0$ in $\Da$, and by
Lemma \ref{infdimrep}(2), $\rep^+\hspace{-1pt}(Q)$ has a chain of irreducible epimorphisms from $M_0$ to $N$. In particular, we have shown that ${\rm supp}\hspace{0.4pt}X\subseteq {\rm supp}\hspace{0.4pt}M_0$, for all $X\in \Da$.
Consider now an almost split sequence
\vspace{-3pt}
$$
\xymatrixcolsep{20pt}\xymatrix{\eta: \quad 0 \ar[r]& L \ar[r]& E \ar[r]& M_0
\ar[r]& 0} \vspace{-2pt}$$ in $\rep(Q)$, where $L\in \rep^-(Q)$ is infinite
dimensional. By Corollary
\ref{newcor1.9}, ${\rm supp}\hskip 0.5pt L$ contains a left infinite
path $p$ with $e(p)$ lying in the socle-support of $L$. By Lemma \ref{DTr-supp},
$e(p)$ is a successor in $Q$ of some vertex in ${\rm
supp}\hspace{0.4pt}M_0$. Not being
contained in ${\rm supp}\hspace{0.4pt}M_0$ by Lemma
\ref{Lemma1.6_shiping}(2), $p$ passes through a vertex $x$ which is
not a successor of any vertex in ${\rm supp}\hspace{0.4pt}M_0$. Let
$\Ta$ be the successor-closed subquiver of $Q$ generated by $x$ and
${\rm supp}\hspace{0.4pt}M_0$. By Proposition \ref{prop2.9}(1), $\rep(\Ta)$ has
an almost split sequence$\,:$
\vspace{-2pt}
$$
\xymatrixcolsep{20pt}\xymatrix{\eta_{_{\hspace{-0.5pt}\it\Theta}}: \quad  0
\ar[r] & L_{_{\it \Theta}} \ar[r] & E_{_{\it \Theta}} \ar[r] & M_0
\ar[r] & 0.}\vspace{-2pt}$$

We define a new quiver $Q\hspace{0.2pt}'$ by attaching to $\Ta$ a right infinite path
\vspace{-2pt}
$$u: \;
\xymatrixcolsep{20pt}\xymatrix{ x\ar[r]^-{\alpha_1}& a_1
\ar[r]^-{\alpha_2}& a_2 \ar[r]^-{\alpha_3}& \cdots,} \vspace{-2pt}$$
where $a_i\not\in Q$ for every $i\ge 1$. Let $X\in \Da$. Since ${\rm supp}\hspace{0.4pt}X\subseteq {\rm supp}\hspace{0.4pt}M_0$, we have $X\in \rep(\Ta)$ and $X\in \rep(Q\hspace{0.2pt}')$. Since $x$ is not a successor in $Q$ of any vertex in ${\rm supp}\hspace{0.4pt}M_0$, the vertices in
${\rm supp}\hspace{0.4pt}M_0$ have the same successors in $Q$ and in $Q\hspace{0.2pt}'$. This implies that $X\in \rep^+(Q\hspace{0.2pt}')$. Now, let $q: X\to Y$ with $X, Y\in \Da$ be an irreducible morphism in $\rep^+\hspace{-1pt}(Q)$. By Corollary \ref{irr+rep}, $q$ is irreducible in $\rep(Q)$, and hence irreducible in $\rep(\Ta)$. Since $\Ta$ contains the successors in $Q\hspace{0.2pt}'$ of the vertices in ${\rm supp}\hspace{0.4pt}Y$, by Proposition \ref{ARS-extension}(2), $q$ is irreducible in $\rep(Q\hspace{0.2pt}')$, and hence irreducible in
$\rep^+(Q\hspace{0.2pt}')$. This shows that $\Da$ is a connected subquiver of $\Ga_{\rep^+(Q\hspace{0.2pt}')}$. In particular, $\Da$ is a subquiver of a connected component $\Ga\hspace{0.2pt}'$ of  $\Ga_{\rep^+(Q\hspace{0.2pt}')}\vspace{1pt}$.

Observe that $x$ is a source vertex in $Q\hspace{0.2pt}'$ which is not a successor of any vertex in
${\rm supp}\hspace{0.4pt}M_0$. Thus $\Ta$ contains all the
predecessors in $Q\hspace{0.2pt}'$ of the successors of the vertices in ${\rm
supp}\hskip 0.5pt M_0$. By Proposition \ref{ARS-extension}(1), $\eta_{_{\hspace{-0.5pt}\it\Theta}}$
is an almost split sequence in $\rep(Q\hspace{0.2pt}')$. Since $L_{_{\hspace{-0.5pt}\it \Theta}}$ is finite
dimensional, $\eta_{_{\hspace{-0.5pt}\it\Theta}}$ is also an almost
split sequence in $\rep^+(Q\hspace{0.2pt}')$. That is, $L_{_{\hspace{-0.5pt}\it \Theta}}=\tau_{_{\hspace{-1pt}\it\Theta}}M_0.$
Furthermore, Since $\Ta$ is top-finite, $Q\hspace{0.2pt}'$ has no left infinite path, and hence, $\rep^-(Q\hspace{0.2pt}')={\rm rep}^b(Q\hspace{0.2pt}')$. In particular, $\Ga_{{\rm
rep}^+(Q\hspace{0.2pt}')}$ contains no pseudo-projective representation. Since
$x\in {\rm supp}\hspace{0.4pt} L_{_{\hspace{-0.5pt}\it\Theta}}$, applying Lemma
\ref{stableorbit2} to the infinite acyclic walk $u$, we
see that $L_{_{\hspace{-0.5pt}\it\Theta}}$ is left stable in
$\Ga\hspace{0.2pt}'$, and consequently, $\Ga\hspace{0.2pt}'$ is not preprojective.
Moreover, since $M_0$ is infinite dimensional, $\Ga\hspace{0.2pt}'$ is not
preinjective by Theorem \ref{preinjcomp}. That is, $\Ga\hspace{0.2pt}'$ is a
regular component of $\Ga_{\rep^+(Q\hspace{0.2pt}')}$. Since $\Ga\hspace{0.2pt}'$ contains infinite dimensional but no
pseudo-projective representations, by Statement (2), its infinite
dimensional representations generate a right-most section $\Da'$, which is a trivially valued left infinite path. As a consequence, $\Ga\hspace{0.2pt}'$ is trivially valued.

Let $X\to Y$ be an arrow in $\Da$ whose valuation in $\Ga$ is $(d, d')$. Then $Y\ne M_0$, and by
Lemma \ref{infdimrep}(2), $Y$ is not pseudo-projective. In view of Lemma \ref{infdimrep}(5), we see that $d\hspace{0.2pt}'=1$.
Suppose that $d>1$. Then $\rep^+\hspace{-1pt}(Q)$ has an
irreducible morphism $f: X\to N$ with $N=Y\oplus Y$, which is also an irreducible morphism in $\rep^+(\Ta)$.
We claim that $f$ is irreducible in $\rep^+(Q\hspace{0.2pt}')$.
Indeed, let $f=hg$, where $g: X\to M$ and $h: M\to N$ are morphisms
in $\rep^+\hspace{-1pt}(Q\hspace{0.2pt}')$. Since
$\Ta$ is predecessor-closed in $Q\hspace{0.2pt}'$, by Lemma \ref{rep+restrict}(1),
$M_{_{\it\Theta}}\in \rep^+(\it\Theta)$. This yields a factorization
$f=h_{_{\it\Theta}} \hspace{0.3pt} g_{_{\it\Theta}}$ in $\rep^+(\it\Theta)$.
Therefore, $g_{_{\it\Theta}}$ is a section or $h_{_{\it\Theta}}$ is a
retraction. In the first case, $v' g_{_{\it\Theta}}=\id_X$ for some
morphism $v': M_{_{\it\Theta}}\to X$. Since $\Ta$ contains the
predecessors in $Q\hspace{0.2pt}'$ of the vertices in ${\rm supp}\hspace{0.4pt}X$,
we can extend $v'$ to a morphism $v: M\to X$ in $\rep(Q\hspace{0.2pt}')$ such that
$vg=\id_X$. That is, $g$ is a section. Dually, if $h_{_{\it\Theta}}$
is a retraction, then $h$ is a retraction, since $\Ta$ contains the
successors of the vertices in ${\rm supp}\hspace{0.4pt}N$. This
establishes the claim. As a consequence, the arrow $X\to Y$ in
$\Ga\hspace{0.2pt}'$ has a non-trivial valuation, a contradiction. This proves that
$\Da$ is trivially valued as a valued subquiver of $\Ga$.
Furthermore,
since the representations in $\Da$ are all infinite dimensional,
$\Da$ is a full subquiver of $\Da'$. Having a source vertex, $\Da$
is of the form \vspace{-0pt}
$$\xymatrixcolsep{20pt}\xymatrix{M_0\ar[r]& M_1\ar[r]& \cdots \ar[r] &
M_{n-1}\ar[r]& M_n.}
\vspace{-0pt}$$

If $n=0$, then $\Da=\{M_0\}$, and consequently $\Ga=\{M_0\}$, and we
are done. Suppose now that $n>0$. By Lemma \ref{infdimrep}(2), $\Ga$
contains a path \vspace{-0pt}$$\xymatrixrowsep{15pt}
\xymatrixcolsep{20pt}\xymatrix{p: \qquad \tau^nM_n\ar[r]&
\tau^{n-1}M_{n-1}\ar[r]& \cdots \ar[r] & \tau M_1\ar[r]&
M_0.}\vspace{-2pt}$$
Since $M_0$ is pseudo-projective, by Lemma \ref{infrepsection}(2),
$p$ contains only pseudo-projective representations. Moreover, $p$
meets every $\tau$-orbit in $\Ga$ since so does $\Da$. Thus, $p$ is
the left-most section of $\Ga$ generated by its pseudo-projective
representations. This shows that $\Ga$ is a finite wing. The
proof of the theorem is completed.

\medskip

\vspace{5pt}

\noindent{\sc Remark.} (1) By Theorem \ref{reg-cpt}(4), we
have a one-one correspondence between the infinite
dimensional pseudo-projective representations in $\Ga_{\rep^+\hspace{-0.5pt}(Q)}$
and the finite regular components of $\Ga_{\rep^+\hspace{-0.5pt}(Q)}$.

\smallskip

\noindent (2) Let $\Ga$ be a non-trivial regular component of
$\Ga_{\rep^+\hspace{-0.5pt}(Q)}$. By Theorem \ref{reg-cpt}, $\Ga$ contains a
unique non-trivial $\tau$-orbit $\mathcal O$ in which every
representation which is not pseudo-projective is the ending term of an almost
split sequence with an indecomposable middle term. The
representations in $\mathcal O$ are called {\it quasi-simple}.
Moreover, each representation $M$ in $\Ga$ has a unique sectional path
$\xymatrixcolsep{15pt}\xymatrix{M=M_n \ar[r] & \cdots \ar[r] &
M_1,}$ with $n\ge 1$ and $M_1$ quasi-simple. One calls $n$ the {\it quasi-length} of $M$. For convenience, the only representation in
any trivial regular component is also called {\it quasi-simple}.

\medskip

Applying Theorems \ref{preinjcomp},
\ref{prep-component} and \ref{reg-cpt}, we get immediately the following
result.

\medskip

\begin{Cor}\label{valuation} If $Q$ is an infinite connected strongly locally finite quiver,
then $\Ga_{\rep^+\hspace{-0.5pt}(Q)}$ has a symmetric valuation.
\end{Cor}

\medskip

In the next two sections, we shall see that each of the four types
of regular components does occur. To conclude this section, we
give some conditions on $Q$ such that
$\Ga_{\rep^+\hspace{-0.5pt}(Q)}$ has at most one type of regular
components. We start with the case where $Q$ has no infinite path.

\medskip

\begin{Cor}\label{no-infpath-cpt} Let $Q$ be an infinite connected strongly locally finite quiver.
If $Q$ has no infinite path, then $\Ga_{\rep^+\hspace{-0.5pt}(Q)}$
consists of a preprojective component of shape $\N Q^{\,\rm op}$,
a preinjective component of shape $\N^-\hspace{-1pt}Q^{\,\rm op}$,
and possibly some regular components of shape
$\Z\hspace{-0pt}\A_{\infty}$.

\end{Cor}

\noindent{\it Proof.} Assume that $Q$ has no infinite path. Then
${\rm rep}^+(Q)=\rep^b(Q)={\rm rep}^-(Q)$. In parti\-cular,
$\Ga_{\rep^+\hspace{-0.5pt}(Q)}$ has no infinite dimensional or
pseudo-projective representation. Now the result follow
immediately from Theorems \ref{prep-component}(1),
\ref{preinjcomp}(1), and \ref{reg-cpt}(1). The proof of the
corollary is completed.

\medskip

Finally, for convenience, we shall call a non-trivial walk in $Q$
an {\it almost-path} if all but finitely many of its edges are
arrows.

\medskip

\begin{Theo}\label{res-reg-cpt} Let $Q$ be an infinite connected strongly locally finite quiver.

\vspace{-1pt}

\begin{enumerate}[$(1)$]

\item If every right infinite acyclic walk in $Q$ is an
almost-path, then every regular component of
    $\Ga_{\rep^+\hspace{-0.5pt}(Q)}$ is of shape $\N^-\hspace{-1.5pt}\A_{\infty}$.

\vspace{2pt}

\item If every left infinite acyclic walk in $Q$ is an almost-path, then every regular component of
$\Ga_{\rep^+\hspace{-0.5pt}(Q)}$ is of shape $\N\A_{\infty}$.

\end{enumerate}\end{Theo}

\noindent{\it Proof.} (1) Suppose that the right infinite acyclic
walks in $Q$ are all almost-paths. Since the inverse of a left
infinite path is a right infinite acyclic walk which is not an
almost-path, $Q$ contains no left infinite path. Hence,
$\rep^-(Q)=\rep^b(Q)$, and in particular,
$\Ga_{\rep^+\hspace{-0.5pt}(Q)}\vspace{1pt}$ has no
pseudo-projective representation. Let $\Ga$ be a regular component
of $\Ga_{\rep^+\hspace{-0.5pt}(Q)}\vspace{1pt}$. By Theorem
\ref{reg-cpt}, $\Ga$ is of shape $\Z\A_{\infty}$ or
$\N^-\hspace{-1.5pt}\A_{\infty}$.

Suppose that $\Ga$ is of shape
$\Z\A_{\infty}$. Then, by Theorem \ref{reg-cpt} again,
the representations in $\Ga$ are all finite dimensional. Fix
arbitrarily a representation $M$ in $\Ga\vspace{1pt}$. Observe that ${\rm
supp}(M\oplus \tau^-\hspace{-1pt}M)$ is connected since ${\rm
Ext}^1(\tau^-\hspace{-1pt}M, M)\ne 0$. As a consequence, $\Sa={\rm
supp}(\oplus_{i \ge 0}\,\tau^{-i}\hspace{-1.5pt}M)$ is connected. We
claim that $\Sa$ is finite. Indeed, if $\Sa$ is infinite, then it
contains a right infinite acyclic walk $w$, which is an almost-path
by hypothesis. Write $w = vu$, where $u$ is a finite walk and $v$ is
a right infinite path$\,:$

\vspace{-10pt}
$$\xymatrixcolsep{15pt}\xymatrix{a_0 \ar[r] & a_1 \ar[r] & a_2 \ar[r] & \cdots }\vspace{-0pt}$$

Observe that $a_0\in {\rm supp}\,\tau^{-r}\hspace{-2pt}M$ for some
$r\ge 0$. Since $M\oplus \cdots \oplus \tau^{-r}\hspace{-1.5pt}M$
is finite dimensional, there exists a maximal integer $s$ such
that $a_s\in {\rm supp}(M\oplus \cdots\oplus
\tau^{-r}\hspace{-1.5pt}M)$. Then $a_{s+1}$ lies in the support of
some $\tau^{-j}M$ with $j>r$. Let $t$ be minimal such that the
support of $\tau^{-t}\hspace{-1.5pt}M$ contains some of the $a_i$
with $i > s$. Then $t>r$. Since $\tau^{-t}\hspace{-1.5pt}M$ is
finite dimensional, there exists a maximal integer $l$ such that
$a_l\in {\rm supp}\,\tau^{-t}\hspace{-1.5pt}M$. Then $l>s$. Since
$a_{l+1}\not\in {\rm supp}\,\tau^{-t}\hspace{-1pt}M$ and $a_l\to
a_{l+1}$ is an arrow, by Corollary \ref{DTr-supp}, $a_{l+1} \in
{\rm supp}\,\tau^{-(t-1)}\hspace{-1pt}M$, which is contrary to the
minimality of $t$. Our claim is established. As a consequence,
there exist two distinct integers $m, n\ge 0$ such that
$\tau^{-m}\hspace{-1.5pt}M$ and $\tau^{-n}\hspace{-1.5pt}M$ have
the same support, a contradiction to Proposition \ref{orbits}(2).
Therefore, $\Ga$ is of shape
$\N^-\hspace{-1.5pt}\A_{\infty}\vspace{1pt}$.

(2) Suppose that the left infinite acyclic walks in $Q$ are all
almost-paths. Then $Q$ has no right infinite path, and hence
$\rep^+\hspace{-1pt}(Q)=\rep^b(Q)$. Hence
$\Ga_{\rep^+\hspace{-0.5pt}(Q)}$ has no infinite dimensional
representation. Let $\Ga$ be a regular component of
$\Ga_{\rep^+\hspace{-0.5pt}(Q)}$. By Theorem \ref{reg-cpt}, $\Ga$
is of shape $\Z\A_{\infty}$ or $\N \A_{\infty}$. If $\Ga$ is of
shape $\Z\A_{\infty}$, then every representation in $\Ga$ is
finite dimensional. Hence $\Ga$ is a regular component of
$\Ga_{\rep^-(Q)}$. On the other hand, by the dual of Statement
(1), all the regular components of $\Ga_{\rep^-(Q)}$ are of shape
$\N \A_{\infty}$, a contradiction. Thus $\Ga$ is of shape $\N
\A_{\infty}$. The proof of the theorem is completed.

\medskip

\noindent {\sc Remark.} If $Q$ is constructed from a finite quiver
by attaching finitely many disjoint right infinite paths, then it
clearly satisfies the condition stated in Theorem
\ref{res-reg-cpt}(1). Note, however, that the quiver
\vspace{-2pt}
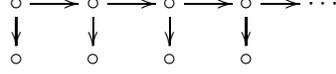
$$
\xymatrixrowsep{20pt}\xymatrixcolsep{28pt}
\xymatrix@!=1pt
{\circ \ar[r] \ar[d] & \circ
\ar[r]\ar[d] & \circ \ar[r]\ar[d] & \circ \ar[r] \ar[d] &  \cdots\\
\circ & \circ &\circ &\circ &}\vspace{0pt}$$ satisfies the stated
condition, but it can not be constructed in this way.

\vspace{-5pt}

{\center \section{Representations of infinite Dynkin quivers}}

Throughout this section, let $Q$ stand for an infinite Dynkin quiver, that
is, its underlying graph is one of the following three infinite diagrams$\,:$

\vspace{10pt}

\hspace{30pt} $\begin{array}{lllllllllllllllllll}
\vspace{40pt} \mathbb{A}_\infty^\infty : && \cdots \; \mbox{---}\; \circ \;\mbox{---}\; \circ \; \mbox{---} \; \circ \; \mbox{---} \;  \circ \;\mbox{---}\;\cdots\\ [2.5mm]
\end{array}$

\vspace{-40pt}

\hspace{30pt} $\begin{array}{lllllllllllllllllll}
\mathbb{A}_\infty :        && \circ  \; \mbox{---}\; \circ \; \mbox{---} \; \circ\; \mbox{---} \; \circ \; \mbox{---} \; \cdots\\ [2.5mm]
\end{array} $

\hspace{30pt} $\begin{array}{lllllllllllllllllll}
\mathbb{D}_\infty :        && \circ  \; \mbox{---}\;\circ \; \mbox{---} \; \circ \; \mbox{---} \; \circ \; \mbox{---} \;  \cdots \\
                           && \hskip 24.5pt | \\
                           &&  \hskip 23pt \circ \\
 \end{array} \vspace{8pt}$

\noindent

As main results, we shall give a complete list of the non-isomorphic
indecomposable representations in $\rep^+\hspace{-1pt}(Q)$ and describe
explicitly its Auslander-Reiten components. Note that
Reiten and Van den Bergh have done so (with no proof) for each type of infinite Dynkin quivers with
the alternating orientation; see
\cite[(III.3)]{RVDB}.

\medskip

As usual, some combinatorial consideration is needed. Let $w$ be a
reduced walk in $Q$. Denote by $Q(w)$ the full subquiver of $Q$
generated by the vertices appearing in $w$. We say that $w$ has no
left infinite path or no right infinite path if so does $Q(w)$. Now,
$w$ is called a {\it string} if the quiver $Q(w)$ contains at least
one, and at most finitely many, sink or source vertices. If $w$ is a
non-trivial string, then neither $w$ nor $w^{-1}$ is a double
infinite path, and we can write uniquely $w$ as $w=w_n\cdots w_1$,
where $w_1, \ldots, w_n$ are non-trivial paths or inverses of
non-trivial paths such that $w_{i+1}w_i$ is neither a path nor the
inverse of a path, for $1\le i<n$. In this case, we call $w_1$ the
{\it initial walk}, and $w_n$ the {\it terminal walk}, of $w$. Let
$v, w$ be strings. In case $e(v)=s(w)$, we define the {\it
composite} of $v, w$ to be $wv$ if it is a non-trivial reduced walk,
or $w$ if $v$ is trivial, or $v$ if $w$ is trivial. For instance, if
$\alpha:x\to y$ is an arrow, then
$\alpha^{-1}\varepsilon_y=\alpha^{-1}=\varepsilon_x\alpha^{-1}$, but
$\varepsilon_x\ne \alpha^{-1}\alpha$, since $\alpha, \alpha^{-1}$
are not composable as strings. Now, $v$ is called a {\it substring}
of $w$ if $w=s v r$, where $r, s$ are strings.

\medskip

\begin{Defn}\label{spec-string} Let $Q$ be a quiver of type $\A_\infty^\infty$ or $\A_\infty$. Each arrow $\alpha: y\leftarrow x$ determines a unique triple $(q, \alpha, p)$, called a {\it double-hook}, where $q$ is the longest path ending in $y$ but not with $\alpha$, and $p$ is the longest path stating in $x$ but not with $\alpha$.

\end{Defn}

\medskip

\noindent{\sc Remark.} A double-hook $(q, \alpha, p)$ has no left infinite path if and only if $q$ is finite. In this case, $p\,\alpha^{-1}\hspace{-1pt}q$ is a string with no left infinite path.

\medskip

\noindent {\sc Example.} Let $Q$ be a quiver of type $\A_\infty$ as follows$\,:$

\vspace{-6pt}
$$\footnotesize
\xymatrixcolsep{22pt}\xymatrix@!=1pt{0 \ar[r]^\alpha &1 \ar[r]^\beta & 2 &3\ar[r]\ar[l]_\gamma & 4\ar[r]\ar[r] &5\ar[r]&6\ar[r]& \cdots &}
\vspace{2pt}$$
Then $(\varepsilon_1, \alpha, \varepsilon_0)$ and $(\beta\alpha, \gamma, p_{_{\infty}})$ are double-hooks,  where $p_{_{\infty}}$ denotes the infinite path starting in $3$.

\medskip

For convenience, we shall say that $Q$ is {\it canonical} of type $\A_\infty$ if
$Q_0=\mathbb{N}$ and the edges are
of the form $x$ --- $x+1$; {\it canonical} of type $\A_\infty^\infty$ if
$Q_0=\mathbb{Z}$ and the edges are of the form $x$ ---
$x+1$; and {\it canonical} of type $\mathbb{D}_\infty$ if $Q_0=\N$ such
that $0, 1$ are of weight one and $2$ is of weight three and the
edges not attached to $0$ are of the form $x$ --- $x+1$.

\medskip

\begin{Defn} \label{order}
Let $Q$ be a canonical infinite Dynkin quiver, and let $\mathscr{S}$
be a set of paths having pairwise distinct starting points.
\vspace{-1pt}

\begin{enumerate}[$(1)$]

\item For $p, q\in \mathscr{S}$, we define $p\preceq q$ if and only if $s(p)\le s(q)$.

\vspace{1pt}

\item For $p, q\in \mathscr{S}$, we define
$p=\sigma_{_{\hspace{-1.5pt}\mathscr{S}}}(q)$ and $q=\sigma_{_{\hspace{-1.5pt}\mathscr{S}}}^-(p)$ if $p\prec q$ and there exists no $u$ in $\mathscr{S}$ such that $p\prec u\prec q$.

\vspace{1pt}

\item We call $\sigma_{_{\hspace{-1.5pt}\mathscr{S}}}$ the
\emph{source-translation} in $\mathscr{S}$.

\end{enumerate} \end{Defn}

\medskip

\noindent{\sc Remark.} (1) The relation $\preceq$ is a well order in $\mathscr{S}.$

\vspace{1pt}

\noindent  (2) If $p, q\in \mathscr{S}$, then $p\prec q$ if
and only if $p=\sigma^i(q)$ for some $i>0$.

\medskip

Let $Q$ be canonical of type $\A_\infty$ or $\A_\infty^\infty$. A
non-trivial path in $Q$ is called {\it right-oriented} or {\it
left-oriented} if its arrows are all of the form $x\to x+1$ or all
of the form $x-1\leftarrow x$, respectively. Moreover, a string $w$
is called {\it normalized} if $s(u)\le e(u)$ for any finite
substring $u$ of $w$. It is evident that $w$ or $w^{-1}$ is
normalized.

\medskip

\begin{Notn}\label{RL-sets}
Suppose that $Q$ is a canonical quiver of type $\A_\infty$ or
$\A_\infty^\infty$.

\vspace{-2pt}

\begin{enumerate}[$(1)$]

\item Let $Q_R$ denote the set of right-oriented
maximal paths having a starting point and the trivial paths
$\varepsilon_x$, where $x$ is either a middle point of a left-oriented path or a sink vertex of weight one.

\vspace{1pt}

\item Let $Q_L$ denote the set of left-oriented maximal paths having a starting
point and the trivial paths $\varepsilon_x$, where $x$ is either a
middle point of a right-oriented path or a source vertex of weight
one.

\end{enumerate}\end{Notn}

\medskip

We state some alternative defining properties of the paths in $Q_R$.

\medskip

\begin{Lemma} \label{right-maximal} Let $Q$ be a canonical quiver of type $\mathbb{A}_{\infty}$ or
$\mathbb{A}_{\infty}^{\infty}$.

\vspace{-1pt}

\begin{enumerate}[$(1)$]

\item If $x\in Q_0$ and $p$ is a path, then $p\in Q_R$ with $s(p)=x+1$ if and only if $x\leftarrow x+1$ is an arrow, and $p$ is the longest path starting in $x+1$ but not with $x\leftarrow x+1$.

\item If $q$ is a finite path, then $q\in Q_R$ with $e(q)=x$ if and only if $x\leftarrow x+1$ is an arrow, and $q$ is the longest path ending in $x$ but not with $x\leftarrow x+1$.

\end{enumerate} \end{Lemma}

\noindent{\it Proof.} (1) Let $x\in Q_0$ and $p$ be a path. Suppose
that $\alpha: x\leftarrow x+1$ is an arrow, and $p$ is the longest
path starting in $x+1$ but not with $\alpha$. If $Q$ has an arrow
$x+1\leftarrow x+2$, then $p=\varepsilon_{x+1}\in Q_R$. If $Q$ has
an arrow $\gamma: x+1\to x+2$, then $p$ is the maximal path stating
with $\gamma$. Being right-oriented, $p$ lies in $Q_R$.
Conversely, suppose that $p\in Q_R$ with $s(p)=x+1$. Assume that $p$ is non-trivial. Then $p$ is a right-oriented maximal path starting in $x+1$. Hence, $Q$ has an arrow $\alpha: x\leftarrow x+1$. Being maximal and right-oriented, $p$ is the longest path starting in $x+1$ but not with $\alpha$. Otherwise, $p=\varepsilon_{x+1}$ with $x+1$ the middle point of a left-oriented path or a sink vertex of weight one. Since $x\in Q_0$, the second case does not occur. That is, $Q$ has a path $x\leftarrow x+1\leftarrow x+2$. In this situation, $\varepsilon_{x+1}$ is the longest path starting in $x+1$ but not with $\alpha$.

(2) Let $q$ be a finite path in $Q$. Suppose that $\alpha:
x\leftarrow x+1$ is an arrow, and $q$ is the longest path ending in
$x$ but not with $\alpha$. If $x$ is a sink vertex of weight one or
$x-1\leftarrow x$ is an arrow, then $q=\varepsilon_x\in Q_R$.
Otherwise, $\delta: x-1\to x$ is an arrow, and $q$ is the maximal
path ending with $\gamma$. Being right-oriented and a maximal path,
$q\in Q_R$.
Conversely, suppose that $q\in Q_R$ with $e(q)=x$. Assume that $q$
is non-trivial. Then it is a right-oriented maximal path ending in
$x$. Hence, $\alpha: x\leftarrow x+1$ is an arrow. Being maximal and
right-oriented, $q$ is the longest path ending in $x$ but not with
$\alpha$. Otherwise, $q=\varepsilon_x$ with $x$ the middle point of
a left-oriented path or a sink vertex of weight one. In either case,
$\alpha: x\leftarrow x+1$ is an arrow and $\varepsilon_x$ is the
longest path ending in $x$ but not with $\alpha$. The proof of the
lemma is completed.

\medskip

Similarly, we have some alternative defining properties of the paths in $Q_L$.

\medskip

\begin{Lemma} \label{left-maximal} Let $Q$ be a canonical quiver of type $\mathbb{A}_{\infty}$ or
$\mathbb{A}_{\infty}^{\infty}$.

\vspace{-1pt}

\begin{enumerate}[$(1)$]

\item If $x\in Q_0$ and $p$ is a path, then $p\in Q_L$ with $s(p)=x$ if and only if $x\to x+1$ is an arrow, and $p$ is the longest path starting in $x$ but not with $x\to x+1$.

\vspace{1pt}

\item If $x\in Q_0$ and $q$ is a finite path, then $q\in Q_L$ with $e(q)=x+1$ if and only if $x\to x+1$ is an arrow, and $q$ is the longest path ending in $x+1$ but not with $x\to x+1$.

\end{enumerate} \end{Lemma}

\medskip

Let $Q$ be canonical of type $\mathbb{A}_{\infty}$ or
$\mathbb{A}_{\infty}^{\infty}$. Since $Q_R$ and
$Q_L$ are sets of paths having pairwise distinct starting points, they are equipped with the well order
$\preceq$ and the source-translation which is denoted by $\sigma_{_{\hspace{-1.5pt} R}}$ for $Q_R$ and by $\sigma_{_{\hspace{-1.5pt} L}}$ for $Q_L$. If no risk of confusion is possible, the subscripts in $\sigma_{_{\hspace{-1.5pt} R}}$ and $\sigma_{_{\hspace{-1.5pt} L}}$ will be dropped. The following result reveals the link between the source-translates and the double-hooks.

\medskip

\begin{Lemma} \label{sigma-ass} Let $Q$ be a canonical quiver of type $\mathbb{A}_{\infty}$ or
$\mathbb{A}_{\infty}^{\infty}$.

\begin{enumerate}[$(1)$]

\item If $p, q$ are paths in $Q$, then $p, q\in Q_R$ with $q=\sigma_{_{\hspace{-1pt}R}}(p)$ if and only if $q$ is a finite path in $Q_R$ and $\alpha:  e(q)\leftarrow s(p)$ is an arrow such that $(q, \alpha, p)$ is a double-hook.

\vspace{1pt}

\item If $p, q$ are paths in $Q$, then $p, q\in Q_L$ with $q=\sigma_{_{\hspace{-1pt}L}}^-(p)$ if and only if $q$ is a finite path in $Q_L$ and $\beta: e(q)\leftarrow s(p)$ is an arrow such that $(q, \beta, p )$ is a double-hook.

\vspace{1pt}

\item If $(q, \alpha, p)$ is a double-hook with $q$ finite, then $p, q\in Q_R$ or $p, q\in Q_L$.

\end{enumerate}\end{Lemma}

\noindent{\it Proof.} Let $p, q$ be paths in $Q$. Suppose that $q$
is a finite path in $Q_R$ with $e(q)=x$ and $\alpha: x\leftarrow
s(p)$ is an arrow such that $(q, \alpha, p)$ is a double-hook. By
Lemma \ref{right-maximal}(2), $s(p)=x+1$, and by Lemma
\ref{right-maximal}(1), $p\in Q_R$. Since $s(q)\le x<x+1=s(p)$, we
have $q\prec p$. If $q$ is trivial, then $q=\varepsilon_x$ with
$s(q)=s(p)-1$, and hence $q=\sigma_{_{\hspace{-1pt}R}}(p)$.
Otherwise, $q$ is a right-oriented maximal path ending in $x$. Thus,
for any $y\in Q_0$ with $s(q)<y\le x$, we have an arrow $y-1\to y$,
and by Lemma \ref{right-maximal}(1), $y\ne s(v)$, for any $v\in
Q_R$. Thus, $\sigma_{_{\hspace{-1pt}R}}(p)=q$.

Conversely, suppose that $p, q\in Q_R$ with $q=\sigma_{_{\hspace{-1pt}R}}(p)$. Write $b=s(q)<s(p)=a$. In particular, $a-1\in Q_0$. By Lemma \ref{right-maximal}(1), $\alpha: a-1\leftarrow a$ is an arrow, and $p$ is the longest path starting in $a$ but not with $\alpha$.
Let $(u, \alpha, p)$ be the double-hook determined by $\alpha$, that is, $u$ is the longest path ending in $a-1$ but not with $\alpha$.
If $u$ is infinite, then it is right-oriented and contains the arrow $b-1\to b$, which contradicts Lemma \ref{right-maximal}(1) since $q\in Q_R$.
Thus $u$ is finite. By the sufficiency we have proved, $u=\sigma_{_{\hspace{-1pt}R}}(p)=q$. This establishes Statement (1). Similarly, we may prove Statement (2).

\vspace{1pt}

Finally, let $(q, \alpha, p)$ be a double-hook with $q$ finite. If
$\alpha$ is a left-oriented arrow $x\leftarrow x+1$, then $p, q\in
Q_R$ by Lemma \ref{right-maximal}. If $\alpha$ is a right-oriented
arrow $x\to x+1$, then $p, q\in Q_L$ by Lemma \ref{left-maximal}.
The proof of the lemma is completed.

\medskip

\begin{Cor} \label{rl-sigma} Let $Q$ be a canonical quiver of type $\mathbb{A}_{\infty}$ or
$\mathbb{A}_{\infty}^{\infty}$.

\vspace{-2pt}

\begin{enumerate}[$(1)$]

\item If $q\in Q_R$, then $\sigma_{_{\hspace{-1pt}R}}^-(q)$ is defined in $Q_R$ if and only if $q$ is finite.

\vspace{1pt}

\item If $p\in Q_L$, then $\sigma_{_{\hspace{-1pt}L}}(p)$ is defined in $Q_L$ if and only if $p$ is finite with $e(p)-1\in Q_0$.

\end{enumerate}\end{Cor}

\noindent{\it Proof.} Let $q\in Q_R$. If $\sigma_{_{\hspace{-1pt}R}}^-(q)=p\in Q_R$, then $q=\sigma_{_{\hspace{-1pt}R}}(p)$, which is finite by Lemma \ref{sigma-ass}(1). Conversely, suppose that $q$ is finite with $e(q)=x$. By Lemma \ref{right-maximal}(2), $\alpha: x\leftarrow x+1$ is an arrow, and $q$ is the longest path ending in $x$ but not with $\alpha$. Let $(q, \alpha, p)$ be the double-hook determined by $\alpha$.
By Lemma \ref{sigma-ass}(1), $p\in Q_R$ and $q=\sigma_{_{\hspace{-1pt}R}}(p)$. That is, $p=\sigma^-_{_{\hspace{-1pt}R}}(q)$. This establishes Statement (1). Similarly, we can prove Statement (2). The proof of the corollary is completed.

\medskip

We are now ready to study the representation theory of $Q$. Firstly,
we show that, as in the finite Dynkin case, the indecomposable
representations in $\rep^+\hspace{-1pt}(Q)$ are uniquely determined
by their dimension vector.

\medskip

\begin{Prop}\label{5.0}

Let $Q$ be an infinite Dynkin quiver. If $M, \,N$ are indecomposable
objects in $\rep^+\hspace{-1pt}(Q)$, then $M\cong N$ if and only if ${\rm dim}
M(x)={\rm dim}N(x),$ for all $x\in Q_0$.

\end{Prop}

\noindent{\it Proof.} Let $M, \,N$ be representations in $\rep^+\hspace{-1pt}(Q)$
such that ${\rm dim}M(x)={\rm dim}N(x),$ for all $x\in Q_0$. In
particular, $M, N$ have the same support which
is denoted as $Q\hspace{0.2pt}'$. By Theorem \ref{theo1.10_shiping}(1), $M\oplus
N$ is projective restricted to a co-finite successor-closed
subquiver $\Sa$ of $Q\hspace{0.2pt}'$. Let $\Oa$ be the predecessor-closed
subquiver of $Q\hspace{0.2pt}'$ generated by the augmented complement of $\Sa$ in
$Q\hspace{0.2pt}'$ and the top-support of $(M\oplus N)_{_{\hspace{-1pt}\it\Sigma}}$.
By the dual of Theorem \ref{restriction}(2),
$M_{_{\hspace{-1pt}\it\Omega}}$ and $N_{_{\hspace{-1.5pt}\it\Omega}}$ are indecomposable. Since
$Q\hspace{0.2pt}'$ is top-finite and $\Sa$ is co-finite in $Q\hspace{0.2pt}'$, we see that
$\Oa$ is finite. Since $Q$ is of infinite Dynkin type, $\Oa$ is a finite
Dynkin quiver. Therefore, $M_{_{\hspace{-0.8pt}\it\Omega}}\cong N_{_{\hspace{-1.5pt}\it\Omega}}$,
and hence $M\cong N$ by the dual of Theorem \ref{restriction}(1).
The proof of the proposition is completed.

\medskip

In order to classify the indecomposable representations in $\rep^+(Q)$,
we need to define the {\it string} representation $M(w)$, associated to a string $w$ in $Q$, as follows: if $x\in Q_0$, then $M(w)(x)=k$ in
case $x$ appears in $w$ and $M(w)(x)=0$ otherwise; and if $\alpha\in
Q_1$, then $M(w)(\alpha)=\id$ in case $\alpha$ or $\alpha^{-1}$
appears in $w$, and $M(w)(\alpha)=0$ otherwise; compare \cite[page
158]{BR}. By definition, $M(w)=M(w^{-1})$. In case $Q$ is of type $\A_\infty$ or
$\A_\infty^\infty$, we shall prove that the indecomposable representations in $\rep^+(Q)$
are parameterized by the strings having no left infinite path.
Note that this does not follow directly from the result of Butler
and Ringel stated in \cite{BR}. Indeed, if $Q$ contains infinite paths, then
the path algebra $kQ$ is not a string algebra as defined in \cite[Section 3]{BR}. Nevertheless,
Theorem \ref{restriction} allows us to apply their results.

\medskip

\begin{Prop}\label{5.1} Let $Q$ be a quiver of type $\mathbb{A}_\infty^\infty$ or $\mathbb{A}_\infty$. If $M$ is an indecomposable object in $\rep(Q)$, then $M$ is finitely presented if and only if $M\cong M(w)$ with $w$ a string having no left infinite path.

\end{Prop}

\noindent{\it Proof.} We may assume that $Q$ is canonical. For
proving the sufficiency, let $w$ be a string in $Q$ having no left
infinite path, which we may assume to be infinite and normalized.
Let $w_1$ be the initial walk and $w_n$ the terminal walk of $w$.
Suppose that $w_1, w_n$ are distinct and both infinite. Then
$w_1^{-1}$ and $w_n$ are right infinite paths. Observe that ${\rm
supp}{\hskip 0.5pt}M(w)=Q$. Let $\alpha: x\to y$ be the initial
arrow of $w_1^{-1}$ and $\beta: a\to b$ the initial arrow of $w_n$.
Consider the full subquiver $\Sa$ of $Q$ generated by the successors
of $y$ and those of $b$. Then $\Sa$ is successor-closed and
co-finite in $Q$ such that $M(w)_{_{\it\Sigma}}=P_y\oplus
P_b\vspace{1pt}$. By Theorem \ref{theo1.10_shiping}(1), $M(w)\in
\rep^+\hspace{-1pt}(Q)$. In case either $w_1$ or $w_n$ is infinite,
we can prove in a similar manner that $M(w)\in
\rep^+\hspace{-1pt}(Q)$.

Conversely, suppose that $M$ is an indecomposable representation in $\rep^+\hspace{-1pt}(Q)$. Then
${\rm supp}(M)$ is connected, and hence, ${\rm supp}(M)=Q(\omega)$, for some reduced walk  $\omega$.
Being top-finite, ${\rm supp}(M)$ has at most finitely many source vertices and no left infinite path.
That is, $\omega$ is a string having no left infinite path. By Theorem \ref{theo1.10_shiping}(1),
$M$ is projective restricted to a co-finite successor-closed
subquiver $\Sa$ of ${\rm supp}(M)$. Fix a vertex $x$ in ${\rm supp}(M)$. Let $\Oa$
be the predecessor-closed subquiver of ${\rm supp}(M)$ generated by $x$,
the augmented complement of $\Sa$ in ${\rm supp}(M)$ and the top-support of
$M_{_{\hspace{-1pt}\it \Sigma}}$, which is finite.
By Theorem \ref{restriction}(2), $M_{_{\hspace{-0.8pt}\it\Omega}}$ is
indecomposable. Then, by the theorem stated in \cite[Section 3]{BR},
$M_{_{\hspace{-0.8pt}\it\Omega}}$ is a string representation of $\Oa$. In
particular, ${\rm dim}_k M(x)=1$. Since ${\rm supp}(M)=Q(\omega)={\rm supp}\hspace{0.4pt}M(\omega),$
by Lemma \ref{5.0}, $M\cong M(\omega).$ The proof of the proposition is
completed.

\medskip

We describe the irreducible morphisms in $\rep^+\hspace{-1pt}(Q)$ in the following proposition; compare \cite[page 166]{BR}.

\medskip

\begin{Prop} \label{irr-mor} Let $Q$ be a quiver of type $\A_\infty^\infty$ or $\A_\infty$, and let $\alpha$ be an arrow  in $Q$ and $w$ be a string having no left infinite path.

\vspace{-2pt}

\begin{enumerate}[$(1)$]

\item If $p$ is a path of maximal length such that $w\alpha
p^{-1}$or $p{\hskip 1.3pt}\alpha^{-1}w$ is a string, then the canoni\-cal embedding
$M(w)\to M(w\alpha p^{-1})$ or
$M(w)\to M(p{\hskip 1.3pt}\alpha^{-1}w)$ is irreducible in $\rep^+\hspace{-1pt}(Q)$, respectively.

\vspace{1.5pt}

\item If $q$ is a finite path of maximal length such that $w\alpha^{-1} q$ or $q^{-1} \alpha w$ is a string, then the canonical projection $M(w\alpha^{-1}q)\to M(w)$ or $M(q^{-1} \alpha w)\to M(w)$ is
irreducible in $\rep^+\hspace{-1pt}(Q)$, respectively,

\end{enumerate}\end{Prop}

\noindent{\it Proof.} We shall prove the proposition only for one case,
since the other cases can be treated similarly. Suppose
that $Q$ is canonical and $w$ is normalized such that $w\alpha$ is a string. Let $p$ be the longest path
such that $v=w\alpha p^{-1}$ is a string. Since $w$ has no
left infinite path, nor does $v$. The canonical embedding $f: M(w)\to M(v)$ is
defined so that $f(x)=\id$ if $x$ is a vertex lying in $w$; and otherwise, $f(x)=0$.

Suppose that $f=hg$, where $g: M(w)\to N$ and $h: N\to M(v)$ are
morphisms in $\rep^+\hspace{-1pt}(Q)$.
It suffices to prove that $g$ is a
section or $h$ is a retraction. For this purpose, we may assume that
$\Hom(M, M(v))\ne 0$ for any indecomposable summand $M$ of $N$.
Write $L=M(w)\oplus N \oplus M(v)$, and consider $Q\hspace{0.2pt}'={\rm supp}{\hskip 1pt}L$, which is connected by the assumption.
By Theorem \ref{theo1.10_shiping}(1), $L$ is projective restricted to a
co-finite successor-closed subquiver $\Sa$ of $Q\hspace{0.2pt}'$. Being connected,
$Q\hspace{0.2pt}'$ has a finite connected subquiver $\Da$ which contains the arrow
$\alpha$, the top-support of $L_{_{\hspace{-0.8pt}\it\Sigma}}$ and the augmented
complement of $\Sa$ in $Q\hspace{0.2pt}'$. Since $Q\hspace{0.2pt}'$ is top-finite, the
predecessor-closed subquiver $\Oa$ of $Q\hspace{0.2pt}'$ generated by $\Da$ is
finite and connected. Then $M(w)_{_{\hspace{-1pt}\it\Omega}}=M(u)$ and
$M(v)_{_{\hspace{-1pt}\it\Omega}}=M(u\alpha q^{-1})$, where $u=w\,\cap\,\Oa$ and
$q=p\,\cap\,\Oa$. Observe that $q$ is the longest path in $\Oa$ such that
$u\alpha q^{-1}$ is a string in $\Oa$, and $f_{_{\it\Omega}}:
M(u)\to M(u\alpha q^{-1})$ is the canonical embedding. By the lemma
stated in \cite[page 166]{BR}, $f_{_{\hspace{-1pt}\it\Omega}}$ is irreducible in
${\rm rep}(\Oa)$. Thus $g_{_{\hspace{-1pt}\it\Omega}}$ is a section or
$h_{_{\hspace{-1pt}\it\Omega}}$ is a retraction. By the dual of Theorem \ref{restriction}(3), $g$ is a section or $h$ is a retraction. The proof of the
proposition is completed.

\medskip

\begin{Lemma}\label{qw-Lemma} Let $Q$ be a quiver of type $\A_\infty^\infty$ or $\A_\infty$, and let $w$ be a string such that $M(w)$ and $\DTr M(w)$ are indecomposable objects in $\rep^+\hspace{-1pt}(Q)$. If $\alpha$ is an arrow and $q$ is a path such that $w\alpha^{-1}q$ is a string, then $q$ is finite.
\end{Lemma}

\noindent{\it Proof.} Let $\alpha$ be an arrow and $q$ be a path such that $w\alpha^{-1}q$ is a string. Then $M(w)$ admits
a minimal projective resolution \vspace{-2pt}$$
\xymatrixcolsep{20pt}\xymatrix{0\ar[r]&
\oplus_{i=1}^{\,r}\,P_{x_i}\ar[r] &
\oplus_{j=1}^{\,s}\,P_{y_j}\ar[r]& M(w) \ar[r] & 0,} \vspace{-2pt}$$
where $y_1, \ldots, y_s$ are the source vertices in $Q(w)$, and by Lemma \ref{resol-supp}, we may assume that $x_1=e(q)$. By Lemma \ref{AR-translate}, $\DTr
M(w)$ admits a minimal injective co-resolution \vspace{-4pt}
$$\xymatrixcolsep{20pt}\xymatrix{0\ar[r]&\DTr M(w)\ar[r]& \oplus_{i=1}^{\,r}\,I_{x_i}\ar[r] & \oplus_{j=1}^{\,s}\,I_{y_j}\ar[r]& 0.}
\vspace{-2pt}$$
Since $x_1\not\in Q(w)$, none of the vertices in ${\rm supp}(\oplus_{j=1}^{\,s}\,I_{y_j})$ appears in $q$. Hence,
$q$ is contained in the support of $\DTr M(w)$. By Lemma \ref{Lemma1.6_shiping}(2), $q$
is not a left infinite path, and hence it is finite. The proof of the
lemma is completed.

\medskip

As shown below, the almost split sequences with an indecomposable middle term in $\rep^+(Q)$ are parameterized by the double-hooks with no left infinite path; compare \cite[page 174]{BR}.

\medskip

\begin{Prop}\label{5.3} Suppose that $Q$ is a quiver of type $\A_\infty^\infty$ or $\A_\infty$.

\begin{enumerate}[$(1)$]

\item If $(q, \alpha, p)$ is double-hook with $q$ finite, then $\rep^+\hspace{-1pt}(Q)$ has an almost split sequence
\vspace{-7pt}
$$\xymatrixcolsep{20pt}
\xymatrix{0\ar[r]& M(q) \ar[r] & M(p\,\alpha^{-1}{\hskip -1pt}q) \ar[r] & M(p) \ar[r] & 0.}$$

\vspace{0pt}

\item Every almost split sequence with an indecomposable middle term in $\rep^+\hspace{-1pt}(Q)$ is of the form stated above.

\end{enumerate}\end{Prop}

\noindent{\it Proof.} Let $(q, \alpha, p)$ be double-hook with $q$ finite. Consider the short exact sequence
\vspace{-6pt}
$$(1) \qquad
\xymatrixcolsep{20pt}\xymatrix{0\ar[r]& M(q) \ar[r]^-f & M(p\,\alpha^{-1}{\hskip -1pt}q) \ar[r]^-g & M(p) \ar[r] & 0}$$ in $\rep^+\hspace{-1pt}(Q)$, where $f$ is the canonical embedding and $g$ is the canonical projection. By Proposition \ref{irr-mor}, $f$ and $g$ are irreducible. Since $\rep^+\hspace{-1pt}(Q)$ is a Krull-Schmidt category, the sequence (1) is almost split; see \cite[(2.15)]{AuR}. Conversely, assume that
\vspace{-3pt}
$$\xymatrixcolsep{22pt}\xymatrix{(2) \qquad 0\ar[r] & L \ar[r]^f & M \ar[r]^g & N \ar[r] & 0\hskip 66pt}$$
is an almost split sequence in $\rep^+\hspace{-1pt}(Q)$ with $M$
indecomposable. By Proposition \ref{5.1}, $N=M(w)$, where $w$ is a
string with no left infinite path. We shall show that the sequence (2) is as stated in Statement (1) by considering all possible
cases.

Firstly, assume that there exists an arrow $\beta$ such that $w\beta^{-1}$ is a string. Let $q$ be the longest
path such that $w\beta^{-1}q$ is a string. By Lemma \ref{qw-Lemma}, $q$ is finite. Now  $\rep^+\hspace{-1pt}(Q)$ has a short exact sequence
\vspace{-4pt} $$\xymatrixcolsep{20pt}\xymatrix{(3)\qquad 0\ar[r]& M(q)\ar[r]^-{f'} & M(w\beta^{-1}q) \ar[r]^-{g'} & M(w) \ar[r] & 0,} \vspace{-3pt}$$ where $f'$ is the canonical embedding and $g'$ is the canonical projection. By
Lemma \ref{5.1}(2), $g'$ is irreducible. Since $M$ is
indecomposable, the sequence $(3)$ is isomorphic to the sequence $(2)$, and hence it is almost split. Let $p$ be the longest path
such that $p\beta^{-1}q$ is a string. Then $(q, \alpha, p)$ is a double-hook with $q$ finite. By the sufficiency we have proved,
$\xymatrixcolsep{15pt}\xymatrix{ 0\ar[r] & M(q) \ar[r] &
M(p\,\beta^{-1} q) \ar[r] & M(p)\ar[r] & 0}$ is an almost split sequence in $\rep^+(Q)$,
which is isomorphic to the sequence (3). That is, the sequence (2) is of the desired form. Similarly, we can treat the case where there
exists an arrow $\beta$ such that $\beta w$ is a string.

Next, suppose that there exists an arrow $\gamma$ such that
$\gamma^{-1}w$ is a string. Let $p$ be the longest path such
that $p\,\gamma^{-1}w$ is a string. Consider the short exact sequence \vspace{-2pt}  $$
\xymatrixcolsep{20pt}\xymatrix{(4)\qquad 0\ar[r]& M(w) \ar[r]^-{u} &
M(p\,\gamma^{-1}w) \ar[r]^-{v} & M(p) \ar[r] & 0}\vspace{-3pt}$$ in $\rep^+\hspace{-1pt}(Q)$, where $u$ is
the canonical embedding and $v$ is the canonical projection. By
Lemma \ref{irr-mor}(1), $u$ is irreducible. In particular, $M(w)$ is
not injective. Since $M(w)$ is not projective, neither is $M(p
\gamma^{-1}w).$  Since $M$ is indecomposable, $\tau
M(p\,\gamma^{-1}w)=M$. By Proposition \ref{ARrep+}, $M$ is of finite dimension, and so is $M(w)$. Therefore,
$\rep^+\hspace{-1pt}(Q)$ has an almost split sequence
\vspace{-2pt}$$(5) \qquad  \xymatrixcolsep{22pt}\xymatrix{ 0\ar[r]&
M(w) \ar[r] & E \ar[r] & \TrD M(w) \ar[r] & 0.}\vspace{-2pt}$$ Since $M(w)$ is
not projective, $E$ has no projective direct summand, and since $M$ is
indecomposable, so is $E$. As a consequence, the sequence $(4)$ is isomorphic to the sequence $(5)$, and hence it is
almost split. Let $q$ be the longest path such that $p\,\gamma^{-1}q$ is a string. By Lemma \ref{qw-Lemma}, $q$ is finite, and consequently, $\rep^+(Q)$ has an almost split sequence
$$\xymatrixcolsep{15pt}\xymatrix{ 0\ar[r] & M(q) \ar[r] &
M(p\,\beta^{-1} q) \ar[r] & M(p)\ar[r] & 0,}$$ which is isomorphic to the sequence (4).
In particular, $w=q$. Since $M(q)=M(w)$ is
not projective, there exists some arrow $\beta$ such that
$q\beta^{-1}$ is a string. This turns out to be the first case we
have treated. Thus, the sequence (2) is of the desired form.
Similarly, we can deal with the case where there exists an arrow
$\gamma$ such that $w\gamma$ is a string.

Now, we consider the case where $e(w)$ is not defined. Since $w$ is a string with
no left infinite path, we may write $w=p\,v$, where $v$ is a
string and $p$ is a right infinite path with $s(p)$ a source vertex in $Q(w)$.
If $v$ is not trivial, then we can write $w=p\,\beta^{-1}u$,
where $u$ is a string and $\beta$ is an arrow. By Lemma
\ref{irr-mor}(1), there exists an irreducible monomorphism $j: M(u) \to
M(w)$, which is impossible since $M$ is indecomposable. Thus $w=p$. Since $M(p)$ is not
projective, there exists an arrow $\beta$ such that $p\,\beta^{-1}$ is
a string. This is again the first case we have treated.
Finally, we can similarly treat the case where $s(w)$ is not defined.
The proof of the proposition is completed.

\medskip

We shall now describe the Auslander-Reiten components of $\rep^+(Q)$ in case $Q$ is of type $\mathbb{A}_{\infty}$ or
$\mathbb{A}_{\infty}^{\infty}$. By Lemma \ref{5.1}, we may choose the vertex set of $\Ga_{\rep^+(Q)}$ to be the set of the finitely presented string representations.

\medskip

\begin{Lemma} \label{right-path-lemma} Let $Q$ be a canonical quiver of type $\mathbb{A}_{\infty}$ or
$\mathbb{A}_{\infty}^{\infty}$. If $p\in Q_R$, then

\begin{enumerate}[$(1)$]

\item $\tau\hspace{-0.5pt}M(p)$ is defined in $\Ga_{\rep^+(Q)}\vspace{1pt}$ if and only if $\sigma(p)$
 is defined in $Q_R$, and in this case, $\tau M(p)=M(\sigma(p))\,;$

\vspace{1pt}

\item $\tau^-\hspace{-1.5pt}M(p)$ is defined  in $\Ga_{\rep^+(Q)}\vspace{1pt}$ if and only if $\sigma^-\hspace{-1.5pt}(p)$ is defined  in $Q_R$, and in this case, $\tau^-\hspace{-1pt}M(p)=M(\sigma^-\hspace{-1pt}(p)).$

\end{enumerate}\end{Lemma}

\noindent{\it Proof.} (1) Let $p\in Q_R$ with $s(p)=x$. Suppose that
$\tau M(p)$ is defined in $\Ga_{\rep^+(Q)}\vspace{1pt}$. Since $M(p)$ is not
projective, $Q$ has an arrow $\alpha: x-1\leftarrow x$. By Lemma
\ref{right-maximal}(1), $p$ is the longest path starting in $x$ but
not with $\alpha$. Let $q$ be the longest path ending in $x-1$ but
not with $\alpha$. By Lemmas \ref{qw-Lemma} and
\ref{right-maximal}(2), $q$ is a finite path in $Q_R$. By Lemma
\ref{sigma-ass}(1), $q=\sigma(p)$. Conversely, suppose that
$\sigma(p)=q\in Q_R$. By Lemma \ref{sigma-ass}(1), $q$ is finite and
$\alpha: e(q)\leftarrow s(p)$ is an arrow such that $(q, \alpha, p)$
is a double-hook. By Proposition \ref{5.3}(1), $\tau M(p)=M(q)$.

(2) If $\sigma^-(p)=q\in Q_R$, then $p=\sigma(q)$, and hence $\tau^-\hspace{-1pt}M(p)=M(q)$ by Statement (1). Conversely, suppose that $\tau^-\hspace{-1.5pt}M(p)\in \Ga_{\rep^+(Q)}\vspace{1pt}$. By Lemma \ref{psproj}(2), $M(p)$ is finite dimensional, that is, $p$ is finite. By Corollary \ref{rl-sigma}(1), $\sigma^-(p)$ is defined in $Q_R$. The proof of the lemma is completed.

\medskip

Similarly, we have the following statement.

\medskip

\begin{Lemma} \label{left-path-lemma} Let $Q$ be a canonical quiver of type $\mathbb{A}_{\infty}$ or
$\mathbb{A}_{\infty}^{\infty}$. If $q\in Q_L$, then

\vspace{0pt}

\begin{enumerate}[$(1)$]

\item $\tau M(q)$ is defined $\Ga_{\rep^+(Q)}\vspace{1pt}$ if and only if $\sigma^-(q)$ is defined in $Q_L$, and in this case, $\tau M(q)=M(\sigma^-(q))\,;$

\vspace{1pt}

\item $\tau^-\hspace{-1pt}M(q)$ is defined $\Ga_{\rep^+(Q)}\vspace{1pt}$ if and only if $\sigma(q)$ is defined in $Q_L$, and in this case, $\tau^-\hspace{-1pt}M(q)= M(\sigma(q)).$

\end{enumerate}\end{Lemma}

\medskip

Let $\mathcal{O}$ be a $\tau$-orbit in
$\Ga_{\rep^+\hspace{-0.5pt}(Q)}$. We shall say that $\mathcal{O}$
is {\it preprojective, preinjective} or {\it regular} if it
contains preprojective, preinjective or regular representations,
respectively. Furthermore, $\mathcal{O}$ is called {\it
quasi-simple} if it consists of quasi-simple regular
representations.

\medskip

\begin{Prop} \label{distinguishedorbits} Let $Q$ be a canonical quiver of type $\A_\infty^\infty$ or $\A_\infty$, and write $\mathcal{O}_R=\{M(p)\mid p \in Q_R\}$ and $\mathcal{O}_L=\{M(p)\mid p \in Q_L\}$.

\begin{enumerate}[$(1)$]

\item If $Q_R$ is non-empty, then $\mathcal{O}_R$ is a preprojective or quasi-simple $\tau$-orbit in $\Ga_{\rep^+\hspace{-0.5pt}(Q)}$, where the second case occurs if and only if $Q$ is of type $\mathbb{A}_\infty^\infty$.

\vspace{1pt}

\item If $Q_L$ is non-empty, then $\mathcal{O}_L$ is a preinjective or regular $\tau$-orbit in $\Ga_{\rep^+\hspace{-0.5pt}(Q)}$, where the second case occurs if and only if $Q$ is of type $\mathbb{A}_\infty^\infty$.

\vspace{1pt}

\item If $\Ga$ is a regular component of $\Ga_{\rep^+\hspace{-0.5pt}(Q)}\vspace{1pt}$, then it contains either $\mathcal{O}_R$ or $\mathcal{O}_L$ but not both, and consequently, $Q$ is of type $\A_\infty^\infty$.

\end{enumerate}
\end{Prop}

\noindent{\it Proof.} It is evident that for any $p, q\in Q_R$, we have $q=\sigma^i(p)$ for some $i\in \Z$. Making use of Lemma \ref{right-path-lemma}, we deduce that $\mathcal{O}_R$ is a $\tau$-orbit in $\Ga_{\rep^+\hspace{-0.5pt}(Q)}$. Suppose that $\mathcal{O}_R$ is preinjective. Then  $M(p)=I_x$ for some $p\in Q_R$ and $x\in Q^+$. In particular, $p$ is finite. By Corollary \ref{rl-sigma}(1) and Lemma \ref{right-path-lemma}(2), $\tau^-\hspace{-1pt} M(p)$ is defined in $\Ga_{\rep^+(Q)}$, a contradiction.
Thus, $\mathcal{O}_R$ is preprojective or regular. Assume that the second case occurs. Suppose that $\mathcal{O}_R$ is non-trivial.
Let $M(p)$ with $p\in Q_R$ be not pseudo-projective. By Lemma \ref{right-path-lemma}(1), $q=\sigma(p)\in Q_R$. By Lemma \ref{sigma-ass}(1), we have a double hook $(q, \alpha, p)$ with $q$ finite. By Proposition \ref{5.3}, $M(p)$ is quasi-simple. That is, $\mathcal{O}_R$ is quasi-simple. Suppose that now $\mathcal{O}_R$. Being regular, $\mathcal{O}_R=\{M(p)\}$, where $p\in Q_R$ and $M(p)$ is infinite-dimensional and pseudo-projective. In particular, $p$ is infinite. Since $M(p)$ is not projective,  $Q$ contains an arrow $\alpha: s(p)-1\leftarrow s(p)$. Let $(q, \alpha, p)$ be a double hook in $Q$. By Lemma \ref{sigma-ass}(1), $q$ is infinite. As a consequence, $Q=Q(w)$ with $w=p\alpha^{-1}q$. In this case, there exists a canonical projection $f: M(w)\to M(p)$. Let $g: M(u) \to M(p)$ be a non-zero non-isomorphism in $\rep^+(Q)$, where $u$ is a string with no left infinite path. If $M(u)$ is projective, then $u$ is a proper subpath of $p$, and in this case, $g$ factors through $f$. Otherwise,
$u=p\alpha^{-1}v$, where $v$ is a finite subpath of $q$ ending in $s(p)-1$. Thus $g$ factors through $f$. That is, $f$ is not irreducible in $\rep(Q)$, and by Corollary \ref{irr+rep}, it is not irreducible in $\rep^+\hspace{-1pt}(Q)$. Thus, $\Ga_{\rep^+(Q)}$ has no arrow ending in $M(p)$. By Theorem \ref{reg-cpt}, $\{M(p)\}$ is a trivial component of $\Ga_{\rep^+(Q)}\vspace{1pt}$. In particular, $\mathcal{O}_R$ is quasi-simple.

Let $Q$ be of type $\mathbb{A}_{\infty}$. Then $0$ is a sink or source vertex of weight one. In the first
case, $\varepsilon_0\in Q_R$ with $M(\varepsilon_0)=P_0$, and in
the second case, the maximal path $p\,_0$ starting in $0$ lies
in $Q_R$ and $M(p\,_0)=P_0$. In either case, $\mathcal{O}_R$ is preprojective. Conversely, suppose that
$\mathcal{O}_R$ is preprojective. Then, $M(p)=P_x$ for some $p\in Q_R$ with $x=s(p)$. If $x-1\in Q_0$, then $Q$ has an arrow $x-1\leftarrow x$, which is absurd since $M(p)$ is projective. Thus $x-1\not\in Q_0$, that is, $Q$ is of type $\A_\infty$. This establishes Statement (1). In a similar way, we can prove Statement (2).

To prove Statement (3), let $\Ga$ be a regular component of
$\Ga_{\rep^+\hspace{-0.5pt}(Q)}$. Suppose first that $\Ga=\{M(w)\}$,
where $w$ is a normalized string with no left infinite path. By
Theorem \ref{reg-cpt}(4), $M(w)$ is infinite dimensional, and hence
$w$ is infinite. Then either the initial walk $w_1$ of $w$ or the
terminal walk $w_n$ is infinite. In the second case,
being normalized and not the inverse of a left infinite path, $w_n$ is a right-oriented maximal path with a starting point. In particular, $w_n\in Q_R$. If $w\ne w_n$, then $w=w_n\alpha^{-1}v$, where $\alpha$
is an arrow and $v$ is a string. Since $w_n$ is a maximal path, by Proposition \ref{irr-mor}(1), $\rep^+\hspace{-1pt}(Q)$ has an irreducible morphism $M(v)\to M(w)$, a contradiction. Therefore, $w=w_n\in Q_R$, and hence $\Ga=\mathcal{O}_R$. Similarly, if $w_1$ is infinite, then $\Ga=\mathcal{O}_L$. Suppose now that $\Ga$ is non-trivial. By Theorem \ref{reg-cpt}, $\rep^+(Q)$ has an almost split
sequence $\xymatrixcolsep{20pt}\xymatrix{0\ar[r] & L \ar[r] & M
\ar[r] & N\ar[r] & 0},$ where $L, M, N\in \Ga$.
By Proposition \ref{5.3}, there exists a double-hook $(q, \alpha, p)$ such that $L=M(q)\vspace{1pt}$ and $N=M(p)$. By Lemma \ref{sigma-ass}(3), $L, N\in \mathcal{O}_R$ or $L, N\in \mathcal{O}_L$. That is, $\Ga$ contains $\mathcal{O}_R$ or $\mathcal{O}_L$. Finally, observe that
$Q_R\cap Q_L=\emptyset$. Having only one quasi-simple $\tau$-orbit, $\Ga$ does not contain $\mathcal{O}_R$ and $\mathcal{O}_L$. The proof of the proposition is completed.

\medskip

\noindent {\sc Example.} Let $Q$ be a canonical quiver of type
$\A_\infty^\infty$ as follows$\,:$
\vspace{2pt}
$$\footnotesize
\xymatrixcolsep{28pt}\xymatrix@!=1pt{\cdots &0\ar[l] &1\ar[l]& 2\ar[l]\ar[r]
&3 & 4\ar[l]\ar[r] &5\ar[r]&6& 7\ar[l]& \cdots \ar[l]}
\vspace{1pt}$$

\noindent We denote by $p_{\,i,j}$ the path from $i$ to $j$, and by
$p\,_\infty$ the infinite path starting in $2$. Then
$Q_L=\{p\hspace{0.4pt}_\infty, p_{_{\hskip 0.5pt 4,3}}, \varepsilon_5\}$ in which
the action of $\sigma$ is indicated as follows:

\vspace{-10pt}

$$\xymatrixcolsep{15pt}\xymatrix{\varepsilon_5\ar@{.>}[r] & p_{_{\hskip 0.5pt 4,3}} \ar@{.>}[r] & p_\infty.}\vspace{-2pt}$$

\noindent By Lemma \ref{left-path-lemma}(1), the action of $\tau$ in $\mathcal{O}_L$ is indicated as follows:
$$\xymatrixcolsep{15pt}\xymatrix{S_5 & M(p_{_{\hskip 0.5pt 4,3}}) \ar@{.>}[l] & M(p\,_\infty)\ar@{.>}[l]}\hspace{-2pt}.\vspace{-2pt}$$

\noindent On the other hand,
$Q_R=\{\varepsilon_i \mid i \le 1\}\cup \{p_{_{\hskip 0.3pt 2,3}}, p_{_{\hskip 0.3pt 4,6}}\} \cup \{\varepsilon_i \mid i \ge 7\},$
in which the action of $\sigma$ is indicated as follows$\,:$

\vspace{-10pt}

$$\xymatrixcolsep{15pt}\xymatrix{\cdots & \ar@{.>}[l] \varepsilon_{-1} & \varepsilon_0\ar@{.>}[l] & \varepsilon_1 \ar@{.>}[l] & p_{\,2,3}\ar@{.>}[l] &p_{\hskip 0.3pt 4,6}\ar@{.>}[l] & \varepsilon_7 \ar@{.>}[l] &
\varepsilon_8 \ar@{.>}[l] & \cdots \ar@{.>}[l]}$$

\noindent By Lemma \ref{right-path-lemma}(1), the action of $\tau$ in $\mathcal{O}_R$ is indicated as follows:

\vspace{-11pt}

$$
\xymatrixcolsep{15pt}\xymatrix{\cdots & \ar@{.>}[l] S_{-1} & S_0\ar@{.>}[l] & S_1 \ar@{.>}[l] & M(p_{_{\hskip 0.3pt 2,3}})\ar@{.>}[l] &
M(p_{_{\hskip 0.3pt 4,6}}) \ar@{.>}[l] & S_7 \ar@{.>}[l] & S_8 \ar@{.>}[l] & \cdots \ar@{.>}[l]}$$

\medskip

The following result describes the Auslander-Reiten components in the $\A_\infty$-case.

\medskip

\begin{Theo} \label{A-infinity} Suppose that $Q$ is an infinite Dynkin quiver of type $\mathbb{A}_\infty$.

\begin{enumerate}[$(1)$]

\vspace{0pt}

\item If $Q$ has no left infinite path, then $\Ga_{\rep^+\hspace{-0.5pt}(Q)}\vspace{1pt}$ consists of
the preprojective component and a preinjective component of shape
$\mathbb{N}^-\hspace{-1.5pt}\mathbb{A}_{\infty}$.

\vspace{2pt}

\item If $Q$ is a left infinite path, then $\Ga_{\rep^+\hspace{-0.5pt}(Q)}$ consists of the preprojective component
of shape $\mathbb{N}\mathbb{A}_{\infty}$.

\vspace{2pt}

\item If $Q$ is not a left infinite path but has left infinite paths, then $\Ga_{\rep^+\hspace{-0.5pt}(Q)}$ consists of the preprojective component of shape $\mathbb{N}\mathbb{A}_{\infty}$ and a finite preinjective component.

\end{enumerate} \end{Theo}

\noindent{\it Proof.} We may assume that $Q$ is canonical. By
Theorem \ref{prep-component} and Proposition
\ref{distinguishedorbits}, $\Ga_{\rep^+(Q)}$ has a unique
preprojective component but no regular component. Moreover, by
Theorem \ref{preinjcomp}, the preinjective components correspond
bijectively to the connected component of $Q^+$. If $Q$ has no left
infinite path, then $Q^+=Q$, which is connected. Thus
$\Ga_{\rep^+\hspace{-0.5pt}(Q)}$ has a unique preinjective component
which is of shape $\mathbb{N}^-\hspace{-1.5pt}\mathbb{A}_{\infty}$
by Theorem \ref{preinjcomp}(1).

Suppose now that $Q$ contains left infinite paths. Then it has no
right infinite path. By Theorem \ref{prep-component}(1), the
preprojective component is of shape $\mathbb{N}\mathbb{A}_{\infty}$.
If $Q$ is a left infinite path, then $Q^+=\emptyset$, and hence
$\Ga_{\rep^+\hspace{-0.5pt}(Q)}$ has no preinjective component.
Otherwise, $Q$ has a left infinite maximal path with an ending point
$x>0$. Then $x$ is a sink vertex, and $Q^+$ is generated by the
vertices $y$ with $0\le y<x$. In particular, $Q^+$ is finite and
connected. By Theorem \ref{preinjcomp},
$\Ga_{\rep^+\hspace{-0.5pt}(Q)}$ has a unique preinjective component
$\mathcal I$, which has a right-most section of shape $(Q^+)^{\rm \,
op}$ and contains only finite $\tau$-orbits. In particular,
$\mathcal I$ is finite. The proof of the theorem is completed.

\medskip

Next, we shall describe the Auslander-Reiten components of
$\rep^+\hspace{-1pt}(Q)$ in the $\A_\infty^\infty$-case. Since the
preprojective component and the possible preinjective components
have been described in Theorems \ref{prep-component} and
\ref{preinjcomp}, we shall concentrate on the regular components.

\medskip

\begin{Theo}\label{A-double-infity}  Let $Q$ be a quiver of type
$\mathbb{A}_{\infty}^{\infty}\vspace{1pt}$, having $r$ right infinite maximal paths and $l$ left infinite maximal paths with $0\le l, r\le 2$.
Then $\Ga_{\rep^+\hspace{-0.5pt}(Q)}$ consists of the preprojective component, at most one preinjective component,
and at most two regular components which are described as follows.

\vspace{-2pt}

\begin{enumerate}[$(1)$]

\item If $Q$ is a double infinite path, then $\Ga_{\rep^+\hspace{-0.5pt}(Q)}$ has a unique regular component of shape $\mathbb{Z}\mathbb{A}_{\infty}$.

\vspace{1pt}

\item If $Q$ has no left infinite path, then $\Ga_{\rep^+\hspace{-0.5pt}(Q)}$ has two regular
components  of which $r$ are of shape $\mathbb{N}^-\hspace{-1.5pt}\mathbb{A}_{\infty}$ and $(2-r)$ are of shape $\mathbb{Z}\mathbb{A}_{\infty}$.

\vspace{1pt}

\item If $Q$ has no right infinite path, then $\Ga_{\rep^+\hspace{-0.5pt}(Q)}$ has two regular
components  of which $l$ are of shape $\mathbb{N} \mathbb{A}_{\infty}$ and $(2-l)$ are of shape $\mathbb{Z}\mathbb{A}_{\infty}$.

\vspace{1pt}

\item If $Q$ has a left infinite maximal path and a right infinite maximal path,
then $\Ga_{\rep^+\hspace{-0.5pt}(Q)}$ has two regular components of which one is of shape
$\mathbb{Z}\mathbb{A}_{\infty}$ and the other one is a finite wing.

\end{enumerate}\end{Theo}

\noindent{\it Proof.} We may assume that $Q$ is canonical. It is easy to see that
$Q^+$ is either empty or connected. Hence $\Ga_{\rep^+\hspace{-0.5pt}(Q)}$ has at most one
preinjective component. Moreover, by Proposition \ref{distinguishedorbits}(3), $\Ga_{\rep^+\hspace{-0.5pt}(Q)}\vspace{1pt}$ has
at most two regular components $\mathcal R$ and $\mathcal L$, such that $M(p)\in \mathcal{R}$ for $p\in Q_R$,
and $M(q)\in \mathcal{L}$ for $q\in Q_L$.

\vspace{1pt}

(1) Assume that $Q$ is a double infinite path in which the arrows are all right-oriented. Then $Q_R =
\emptyset$ and $Q_L=\{\varepsilon_x \mid x\in Q_0\}$. By Proposition \ref{distinguishedorbits}, $\mathcal L$ is the only
regular component. Since $\varepsilon_i=\sigma^{-i}_{_{\hspace{-1pt}L}}(\varepsilon_0)\vspace{1pt}$, by Lemma \ref{left-path-lemma}, $\tau^iS_0=S_i$, for $i\in \Z$. Hence,
$\mathcal L$ is stable. By Theorem \ref{reg-cpt}(1), $\mathcal L$ is of shape $\mathbb{Z}\mathbb{A}_{\infty}$.

\vspace{1pt}

(2) Suppose that $Q$ has no left infinite path. Since $Q_0=\Z$,
there exist at most two right infinite maximal paths. Assume first
that $Q$ has no right infinite path, that is, it has no infinite
path. Then $Q$ can be viewed of the following form$\,:$
$$ \cdots p_nq_n^{-1}\cdots  q_1^{-1}p_{\,0}
q_0^{-1}p_{_{-1}}\cdots p_{_{-m}}q_{-m}^{-1} \cdots$$
where the $p_n$ are the right-oriented maximal paths and the $q_n$ are the
left-oriented maximal paths. Thus $Q_R$ contains a double infinite chain
$$\cdots \prec p_{-m}\prec \cdots \prec p_0 \prec \cdots \prec p_n\prec \cdots,$$
and hence $\sigma^i_{_{\hspace{-1pt}R}}\hspace{-1pt}(p_{\hspace{0.5pt}0})$ is defined for all $i\in \Z$. By Lemma \ref{right-path-lemma},
$\tau^i\hspace{-1.5pt}M(p_{\hspace{0.5pt}0})$ is defined for all $i\in \Z$, that is, $M(p_{\hspace{0.5pt}0})$ is stable. Moreover,
$Q_L$ contains an infinite chain
$$\cdots \prec q_{-m}\prec \cdots \prec q_0 \prec \cdots \prec q_n\prec \cdots,$$
and thus $\sigma^{-i}_{_{\hspace{-1pt}L}}\hspace{-1pt}(q_0)$ is defined for all $i\in \Z$. By Lemma \ref{left-path-lemma}, $\tau^i\hspace{-1.5pt}M(q_0)$ is defined for all $i\in \Z$. That is, $M(q_0)$ is stable. Therefore, $\Ga_{\rep^+(Q)}$ has two stable regular components $\mathcal R$ and $\mathcal L$. By Theorem \ref{reg-cpt}(1), they both are of shape $\Z\A_\infty$.

Assume next that $Q$ has exactly one right infinite maximal path
$q_0$. We may assume that $q_{\,0}$ is left-oriented starting in
$x$. Since $q_0$ is the unique infinite maximal path, $Q$ can be
viewed of the following form$\,:$
$$\cdots q_n^{-1}p_n\cdots q_1^{-1}p_1q_{\hskip 0.5pt 0}^{-1},$$
where the $q_n$ are the left-oriented maximal paths and the $p_n$ are the right-oriented maximal paths. Then
$Q_L$ contains a right infinite chain
$$q_0 \prec q_1 \prec \cdots \prec q_n\prec \cdots,$$ and thus,
$\sigma^{-i}_{_{\hspace{-1pt}L}}\hspace{-1pt}(q_0)$ is defined for
all $i\ge 0$. By Lemma \ref{left-path-lemma}(1),
$\tau^i\hspace{-1pt}M(q_0)$ is defined for all $i\ge 0$. That is,
$M(q_0)$ is left stable. Hence $\mathcal L$ is a left stable regular
component. Since $M(q_0)$ is infinite dimensional, by Theorem
\ref{reg-cpt}(2), $\mathcal L$ is of shape
$\mathbb{N}^-\hspace{-1.5pt}\mathbb{A}_\infty$. On the other hand,
$Q_R$ has a double infinite chain $$\cdots \prec \varepsilon_{x-i}
\prec \cdots \prec \varepsilon_{x-1}\prec p_1 \prec \cdots \prec
p_n\prec \cdots.$$ Making use of Lemma \ref{right-path-lemma} again,
we see that $\mathcal R$ is stable, which is of shape $\Z\A_\infty$
by Theorem \ref{reg-cpt}.

Assume finally that $Q$ has two right infinite maximal paths $p$ and $q$. We may assume that $p$ is right-oriented starting in $x$, and $q$ is
left-oriented starting in $y$. Then $Q_R$ has a left infinite chain
$$\;\cdots \prec \varepsilon_{y-i} \prec \cdots \prec \varepsilon_{y-1}\prec p.$$ In view of Lemma \ref{right-path-lemma}(1), we deduce that
$M(p)$ is left stable. Hence $\mathcal R$ is a left stable regular component. Since $M(p)$ is infinite dimensional, $\mathcal R$ is of shape $\N^-\hspace{-1.5pt}\A_\infty$. Moreover, since $Q_L$ contains a right infinite chain
$$q \prec\varepsilon_{x+1}  \prec \cdots \prec \varepsilon_{x+i} \prec  \cdots,\,$$ by Lemma \ref{left-path-lemma}(1), $M(q)$ is left stable. Hence, $\mathcal L$ is a left stable regular component. Since $M(q)$ is infinite dimensional, $\mathcal L$ is of shape  $\N^-\hspace{-1.5pt}\A_\infty$.

(3) Suppose that $Q$ has no right infinite path. Using an argument dual to that for proving Statement (2), we may show that if
$Q$ has exactly one left infinite maximal path which is assumed to be
right-oriented, then $\mathcal L$ is a regular component
of shape $\Z\A_\infty$ and $\mathcal R$ is a regular component of shape
$\N\A_\infty;$ and  if $Q$ has two left infinite maximal paths,
then $\mathcal R$ and $\mathcal L$ are two regular components of shape $\N\A_\infty$.

(4) Suppose that $Q$ has a left infinite maximal path $p$ and a right infinite maximal path $q$. Then $p, q$ are either both right-oriented or both left-oriented. We need only to consider the case where $p$ is right-oriented with $e(p)=x$, while $q$ is right-oriented with $s(q)=y$.
Then $Q_L$ has a double infinite chain
\vspace{-3pt}
$$\cdots \prec \varepsilon_{x-i}\prec\cdots \prec \varepsilon_{x-1}\prec\varepsilon_{y+1}\prec\cdots\prec\varepsilon_{y+j}\prec\cdots. \vspace{-3pt}$$
Therefore, $\mathcal L$ is stable of shape $\Z\A_\infty$. On the other hand, $q\in Q_R$. If $z$ is a vertex with $z<x$ or $z>y$, then  $z$ is a middle point of a right-oriented path. By Lemma \ref{right-maximal}(1), $z\ne s(v)$ for any $v\in Q_R$. Hence, $Q_R$ is finite. By Proposition \ref{distinguishedorbits}(1), $\mathcal R$ has a finite $\tau$-orbit. By Theorem \ref{reg-cpt}(4), $\mathcal R$ is a finite wing. The proof of the theorem is completed.

\medskip

\noindent {\sc Example.} Reconsider the canonical quiver $Q$  of type
$\A_\infty^\infty$ as follows$\,:$
\vspace{2pt}
$$\footnotesize
\xymatrixcolsep{28pt}\xymatrix@!=1pt{\cdots &0\ar[l] &1\ar[l]& 2\ar[l]\ar[r]
&3 & 4\ar[l]\ar[r] &5\ar[r]&6& 7\ar[l]& \cdots \ar[l]}
\vspace{1pt}$$
As seen before,
$Q_L=\{p\hspace{0.3pt}_\infty, p_{\hskip
0.5pt 4,3}, \varepsilon_5\}$. Therefore, the regular component containing the $\tau$-orbit
$\mathcal{O}_L$ is a wing as follows:
\vspace{-2pt}
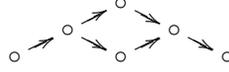
$$\xymatrixrowsep{10pt} \xymatrixcolsep{20pt}
\xymatrix@!=0.1pt{&&&& \circ \ar[dr]&&&&\\
&&&\circ \ar[ur]\ar[dr] 
&& \circ \ar[dr]&&&\\
&& \circ\ar[ur] 
&& \circ\ar[ur] 
&& \circ &&}\vspace{0pt}$$ On the other hand, since both $\sigma$ and $\sigma^-$ are defined everywhere in $Q_R$, the
regular component containing $\mathcal{O}_R$ is of shape
$\Z\A_\infty$.

\medskip

Finally, it comes to the point for us to study the $\D_\infty$-case.
Recall from Lemma \ref{5.0} that the indecomposable representations
in $\rep^+\hspace{-1pt}(Q)$ are uniquely determined by their
dimension vector.

\medskip

\begin{Notn}\label{notation} Let $Q$ be an infinite Dynkin quiver of type $\D_\infty$.

\vspace{-1pt}

\begin{enumerate}[$(1)$]

\item For integers $i\ge 0$ and $j\ge 1$, denote
by $N_{i, j}$ the finite dimensional indecomposable representation in $\rep^+(Q)$
with dimension vector indicated as follows$\,:$
$$\vspace{8pt}
 \begin{array}{llllllllllllllll}  1 \hspace{-5pt} &  \hspace{-5pt}  \mbox{ ---} \hspace{-5pt}
 &  \overbrace{  2  \;\,  \mbox{---}  \;\,
 2  \;\,  \mbox{---}  \;\, \cdots \;\, \mbox{---} \;\,  2 }^{i \;\text{times}} & \hspace{-5pt} \mbox{---} \hspace{-5pt}
 & \overbrace{ 1 \;\, \mbox{---} \;\, \cdots  \;\, \mbox{---} \;\, 1}^{j \;\text{times}} & \hspace{-5pt} \mbox{---} \hspace{-5pt}
 & \hspace{-1pt} 0 \hspace{-1pt} & \hspace{-5pt} \mbox{---} \hspace{-5pt} &
\hspace{-1pt} 0 \hspace{-1pt} & \hspace{-5pt} \mbox{---} \; \cdots \\
&& \,| &&&&&&&&\\
&& 1&&&&&&&&\end{array} \vspace{-8pt}$$

\item In case $Q$ has infinite paths, for each $i\ge 0$, denote by $N_{i, \infty}$ the sincere
indecom\-posable representation in $\rep(Q)$ with dimension vector indicated as
follows$\,:$
$$\vspace{8pt}
 \begin{array}{llllllllllllllllll}  1 \hspace{-5pt} &  \hspace{-5pt}  \mbox{ ---} \hspace{-5pt}
 &  \overbrace{  2  \;\,  \mbox{---}  \;\, 2  \;\,  \mbox{---}  \;\, \cdots \;\, \mbox{---} \;\,  2 }^{i \;\text{times}} & \hspace{-5pt} \mbox{---} \hspace{-5pt}
& 1 \;\, \mbox{---} \;\, \cdots  \;\, \mbox{---} \;\, 1  &  \hspace{-5pt} \mbox{---} \; \cdots \\
&& \,| &&&&&&&&\\
&& 1&&&&&&&&\end{array} \vspace{-8pt}$$

\end{enumerate} \end{Notn}

\medskip

The following result describes the indecomposable representations in $\rep^+\hspace{-1pt}(Q)$.

\medskip

\begin{Prop} \label{indec-D} Let $Q$ be a quiver of type $\D_\infty$. If $M$ is an indecomposable representation in $\rep^+\hspace{-1pt}(Q)$, then

\vspace{-2pt}

\begin{enumerate}[$(1)$]

\item $M\cong M(w)$ with $w$ a string having no left infinite path, or

\item $M\cong N_{i, j}$ with  $i\ge 0$ and $j\ge 1$, or

\item $M\cong N_{i, \infty}$ with $i\ge 0$, and this occurs if and only if $Q$ has right infinite paths.

\end{enumerate} \end{Prop}

\noindent{\it Proof.} We may assume that $Q$ is canonical. Let $M$ be an
indecomposable representation in $\rep^+\hspace{-1pt}(Q)$. If one of the vertices $0, 1$
is not in ${\rm supp}{\hskip 0.5pt}M$, then $M$ is an indecomposable
representation of a quiver of type $\A_\infty$. By Proposition
\ref{5.1}, $M=M(w)$, where $w$ is a string without left infinite
paths.

Suppose that $0, 1\in {\rm supp}{\hskip 0.5pt}M$. If $M$ is finite
dimensional, then it is an indecomposable representation of ${\rm
supp}\hspace{0.5pt}M$. Since ${\rm supp}\hspace{0.5pt}M$ is of type
$\D_n$, it is well known that $M\cong N_{i, j}$ for some $i\ge 0$
and $j\ge 1\,;$ see, for example, \cite[p. 299]{ASS}. Assume that
$M$ is infinite dimensional. Since $M$ is indecomposable, ${\rm
supp}\hspace{0.5pt}M=Q.$ By Corollary \ref{newcor1.9}, $Q$ contains
a right infinite path which we assume starts in a vertex $a\ge 3$.
By Theorem \ref{theo1.10_shiping}(1), $M$ is projective restricted
to a co-finite successor-closed subquiver $\Sa$ of $Q$. By the dual
of Lemma \ref{lemma1.8_shiping}, we may assume that the vertices $x$
with $x\le a$ are not in $\Sa$. Let $b$ be maximal such that $b$
lies in the top-support of $M_{_{\hspace{-1pt}\it \Sigma}}$, and let
$\Oa$ be the full subquiver of $Q$ generated by the vertices $x$
with $0\le x\le b$. Then $\Oa$ contains the top-support of
$M_{_{\hspace{-1pt}\it \Sigma}}$ and is predecessor-closed in $Q$
since $b>a$. By Theorem \ref{restriction}(1),
$M_{_{\hspace{-0.8pt}\it\Omega}}$ is an indecomposable sincere
representation of $\it\Omega$. Since $\Omega$ is of type $\D_{b+1}$,
it is well known that there exists some $r$ with $0\le r\le b-2$
such that ${\rm dim}\,M(x)=2$ for any $2\le x< 2+r$ and ${\rm
dim}\,M(y)=1$ if $y<2$ or $2+r\le y\le b$. In particular, ${\rm
dim}\,M(b)=1$. Let $c$ be an arbitrary vertex with $c>b$. Applying
the same argument to the full subquiver of $Q$ generated by the
vertices $x$ with $0\le x\le c$, we see that ${\rm dim}\,M(c)=1$.
This shows that $M\cong N_{r, \infty}.$ The proof of the proposition
is completed.

\medskip

\noindent {\sc Remark.} The preceding result says particularly that if $M$ is an indecomposable representation in $\rep^+\hspace{-1pt}(Q)$, then ${\rm dim}\, M(x)\le 2$ for all $x\in Q_0$.

\medskip

As for the two other types, the study of quasi-simple
representations is essential in the description of the
Auslander-Reiten components of $\rep^+\hspace{-1pt}(Q)$.

\medskip

\begin{Lemma} \label{D-qs-orbit}
Suppose that $Q$ is quiver of type $\D_{\infty}$. If $\Ga$ is a regular component
of $\Ga_{\rep^+\hspace{-0.5pt}(Q)}$, then every vertex in $Q$ lies in the support of at most two quasi-simple representations in $\Ga$.

\end{Lemma}

\noindent {\it Proof.} Fix $x\in Q_0$. Let $\Ga$ be a regular
component of $\Ga_{\rep^+\hspace{-0.5pt}(Q)}$ with a sectional path
$L_n\to L_{n-1}\to \cdots \to L_1,$ where $L_1$ is quasi-simple. In
view of the shape of $\Ga$ described in Theorem \ref{reg-cpt}, we
see that $\tau^iL_1\in \Ga\vspace{1pt}$, for $i=1, \ldots, n-1$. Now
an easy induction on $n$ shows that ${\rm dim}\,L_n(x) =
{\sum}_{i=0}^{n-1} \,{\rm dim}\, \tau^iL_1(x).\vspace{1pt}$

Suppose that $\Ga$ contains some distinct quasi-simple
representations $N_1$, $N_2,$ and $N_3$ such that $N_i(x)\ne 0$, for
$i=1, 2, 3$. With no loss of generality, we may assume that
$N_3=\tau^sN_1$ and $N_2=\tau^rN_1$ with $0<r<s$. Then $\Ga$
contains a sectional path $M_{s+1}\to M_{s}\to \cdots \to M_1=N_1$.
This yields \vspace{-1pt}
$${\rm dim}\hspace{0.3pt}
M_{s+1}(x) ={\textstyle \sum}_{i=0}^{s} \,{\rm dim}\,\tau^iN_1(x)
\ge {\rm dim}N_1(x) + {\rm dim}N_2(x) + {\rm dim}N_3(x) \ge 3,
\vspace{-1pt}$$ which is contrary to Proposition \ref{indec-D}. The
proof of the lemma is completed.

\medskip

\begin{Lemma}\label{D-inf-reg} Let $Q$ be a canonical quiver of type $\D_{\infty}$, and let $w$ be a finite string with $3\le s(w)\le e(w)$. If there exist arrows $\alpha, \beta$ such that $\beta w\alpha$ or $\beta^{-1}w\alpha^{-1}$ is a string, then $M(w)$ is a regular representation.
\end{Lemma}

\noindent {\it Proof.} Write $a=s(w)$ and $b=e(w)$. Firstly, suppose that $Q$ has a
string $\beta w\alpha$ with $\alpha, \beta$ being arrows. Then
$\alpha$ is the arrow $(a-1)\to a$ and $\beta$ is the arrow $b\to
(b+1)$. Applying Lemma \ref{stableorbit2}(1) to the unique
infinite acyclic walk in $Q$ starting with $\beta$, we see that
either $M(w)$ is left stable or $\tau^n\hspace{-2pt} M(w)$ is
pseudo-projective for some integer $n\ge 0$. In particular, $M(w)$
is not preprojective.

For each integer $j\ge 1$, consider the indecomposable representation $N_{b-1, j}$ as defined in Notation \ref{notation}(1). We may assume that $N_{b-1, j}(x)=k^2$ for $2\le x\le b$, $N_{b-1, j}(b+1)=k$, and $N_{b-1, j}(\gamma)=\id$ for each arrow $\gamma: x\to y$ with $2\le x, y\le b$. Since $N_{b-1, j}$ is indecomposable, the map $N_{b-1, j}(\beta)$ is surjective with a non-zero kernel $\varphi: k\to k^2$.
On the other hand, $M(w)(x)=k$ for each vertex $x$ appearing in $w$, and $M(w)(\gamma)=\id$ for each arrow $\gamma$ such that $\gamma$ or $\gamma^{-1}$ is an edge in $w$. For $x\in Q_0$, we set $f(x)=\varphi$ if $w$ passes through $x$; and otherwise, $f(x)=0$.
By definition, $N_{b-1, j}(\beta)\, f(b)=0=f(b+1)\, M(w)(\beta)$, and  since $M(w)(a-1)=0$, we have $N_{b-1, j}(\alpha) f(a-1)=0=f(a) M(w)(\alpha)$. Moreover, if $\gamma: x\to y$ is an arrow such that $\gamma$ or $\gamma^{-1}$ is an edge in $w$, then $a\le x, y\le b$, and in this case, we have
$N_{b-1, j}(\gamma)\, f(x)=\varphi=f(y)\, M(w)(\gamma)$.
This verifies that $f=\{f(x)\mid x\in Q_0\}$ is a non-zero morphism in $\rep^+\hspace{-1pt}(Q)$ from $M(w)$ to $N_{b-1, j}$. Since the $N_{b-1, j}$ with $j\ge 1$ are pairwise non-isomorphic, by Proposition \ref{no_ext}(2), $M(w)$ is not preinjective. That is, $M(w)$ is regular.

Next, suppose that $Q$ has a string $\beta^{-1} w\alpha^{-1}$
with $\alpha, \beta$ being arrows. Then we deduce from Lemma \ref{stableorbit}(1) that $M(w)$ is not preinjective. Moreover,  in a dual
manner, we can show that there exists a non-zero morphism from $N_{b-1, j}$ to
$M(w)$ for each $j\ge 1$. Therefore, $M(w)$ is not preprojective by
Proposition \ref{no_ext}(1). That is, $M(w)$ is regular. The proof
of the lemma is completed.


\vspace{6pt}

We are ready to describe the connected components of $\Ga_{\rep^+\hspace{-0.5pt}(Q)}$ in the $\D_\infty$-case. By Theorems \ref{prep-component}
and \ref{preinjcomp}, we need only to concentrate on the regular components.

\vspace{3pt}


\begin{Theo}\label{D-inf-main} Let $Q$ be a quiver of type $\mathbb{D}_{\infty}$.
Then $\Ga_{\rep^+\hspace{-0.5pt}(Q)}\vspace{1pt}$ consists of the preprojective component,
at most one preinjective component, and exactly one regular component which is of shape 
$\mathbb{Z}\mathbb{A}_{\infty}$, $\mathbb{N}\mathbb{A}_{\infty}$, or
$\mathbb{N}^-\hspace{-2pt}\mathbb{A}_{\infty}$ in case $Q$ has no infinite path,
has left infinite paths, or has right infinite paths, respectively.
\end{Theo}

\noindent{\it Proof.} Observing that $Q^+$ is empty or connected, we deduce
easily the first two parts of the statement from Theorems \ref{prep-component} and \ref{preinjcomp}.
For proving
the last part, assume that $Q$ is canonical. If $Q$ has a vertex $x$ with $x>2$ which is a middle point
of some path, then $S_x$ is regular by Lemma
\ref{D-inf-reg}. Otherwise, every arrow $\alpha$ not attached to $2$
is a maximal path, and by Lemma \ref{D-inf-reg} again, $M(\alpha)$ is regular.
This shows that $\Ga_{\rep^+\hspace{-0.5pt}(Q)}$ has at least one
regular component $\Ga$. If $Q$ has no infinite path, then $\Ga$ is
of shape $\Z\A_{\infty}$ by Corollary \ref{no-infpath-cpt}. If $Q$ has
left infinite paths, then every left infinite acyclic walk in $Q$ is
an almost-path, and hence, $\Ga$ is of shape $\N\A_\infty$ by
Theorem \ref{res-reg-cpt}(2). If $Q$ has right infinite paths, then
every right infinite acyclic walk in $Q$ is an almost-path, and
therefore, $\Ga$ is of shape $\N^-\hspace{-1pt}\A_\infty$ by Theorem
\ref{res-reg-cpt}(1).

It remains to show that $\Ga$ is the only regular component of
$\Ga_{\rep^+\hspace{-0.5pt}(Q)}$. For this purpose, we fix a vertex
$a>2$ in such a way that $a$ is a source vertex if $Q$ has no
infinite path, and otherwise, $a$ is a middle point of an infinite
path. Denote by $\Sa$ the full subquiver of $Q$ generated by the
vertices $x\ge a$. Observe that $\Sa$ is a quiver of type
$\A_{\infty}$, which will become canonical if one replaces $x$ by
$x-a$. Let ${\Sa}_R$ and ${\Sa}_L$ be the sets of paths in $\Sa$ as
defined in Notation \ref{RL-sets}, each of them is equipped with a
source-translation written as $\sigma_{_{\hspace{-1pt}R}}$ for
$\Sa_R$ and $\sigma_{_{\hspace{-1pt}L}}$ for $\Sa_L$.

Consider first the case where $Q$ has no left infinite path. Then
$\Ga$ is of shape $\Z\A_\infty$ or $\N^-\hspace{-1pt}\A_\infty$.  In particular,
$\Ga$ is left stable. By Lemma \ref{ARrep+}, all but at most one
quasi-simple representations in $\Ga$ are finite dimensional. By Lemma \ref{D-qs-orbit},
$\Ga$ contains a finite dimensional quasi-simple representation $M$ such that ${\rm
supp}\,\tau^{i}M \subseteq \Sa\backslash\{a\}$ for all $i \ge 0$. In particular, the $\tau^iM$ are
finite dimensional representations in $\rep^+(\Sa)$. For each $i\ge 0$, consider the almost split sequence
$$\eta_i: \quad \xymatrixcolsep{15pt}\xymatrix{0 \ar[r]& \tau^{i+1} M \ar[r]& E_i \ar[r]& \tau^{i}M \ar[r]& 0}
$$
in $\rep^+\hspace{-1pt}(Q)$. We claim that $\eta_i$ is an almost split sequence in $\rep^+(\Sa)$.
Indeed, assume that $Q$ has a right infinite path. Then $\Sa$ is a right infinite path with $a$ being the source vertex. For $x,y \ge a$, we denote by $p_{\hspace{0.5pt}x, y}$ the path in $\Sa$
from $x$ to $y$. Then $M = M(p_{\,r,s})$ for
some  $s\ge r > a$. In view of Lemma \ref{AR-translate}, we see that
$\tau^i\hspace{-1pt}M = M(p_{\,r+i,s+i}),$ and in particular, $\eta_i$ lies in $\rep^+(\Sa)$ for all $i\ge 0$. Hence, our claim follows.
Assume that $Q$ has no infinite path. Then $a$ is a source vertex in $Q$. In particular, $\Sa$ contains the
predecessors of the successors of the vertices in ${\rm
supp}\hspace{0.4pt}\tau^iM$, for all $i\ge 0$. By
Proposition \ref{ARS-extension}(1), $\eta_i$ is an almost split sequence in
$\rep^+(\Sa)$, for all $i\ge 0$. This establishes our claim. In particular, $\tau_{_{\hspace{-1pt}\it
\Sigma}\hspace{0.3pt}}^iM = \tau^{i}M$, for all $i\ge 0$. Since
$E_0$ is indecomposable in $\rep^+(\Sa)$, by
Proposition \ref{5.3}, there exists a double-hook $(q, \beta, p)$ in $\Sa$ such that
$M= M(p)$. By Lemma \ref{sigma-ass}, either $p, q\in \Sa_R$ or $p, q\in \Sa_L$.
Since $M$ is left stable in
$\Ga_{\rep^+(\it\Sigma)}$, by Proposition
\ref{distinguishedorbits}, the second case occurs. Applying Lemma \ref{left-path-lemma}(1), we have
\vspace{-4pt}
$$M(\sigma^{-i}_{_{\hspace{-1pt}L}}(p))=\tau_{_{\hspace{-1pt}\it \Sigma}\hspace{0.3pt}}^i M(p)
=\tau_{_{\hspace{-1pt}\it \Sigma}\hspace{0.3pt}}^i M
=\tau^i M\in \Ga,$$ for all $i\ge 0$. Similarly, if $\Ga\hspace{0.2pt}'$ is
another regular component of $\Ga_{\rep^+\hspace{-0.5pt}(Q)}$, then there exists
some $p'\in \Sa_L$ such that $M(\sigma^{-i}_{_{\hspace{-1pt}L}}(p'))\in \Ga\hspace{0.2pt}'$ for all $i\ge 0$.
Since $p, p'$ both lie in $\Sa_L$, there exists some $t\in \Z$ such that $p'=\sigma^t_{_{\hspace{-1pt}L}}(p)$. As a consequence,
$\Ga$ intersects $\Ga'$, and hence $\Ga=\Ga'$.

Consider now the case where $Q$ has left infinite paths. Then, $\Ga$ is of shape $\N \A_\infty$. In particular,
$\Ga$ contains pseudo-projective but no  infinite dimensional representations. On the other hand, $Q$ has no right infinite path.
By the dual of what we have just proved, $\Ga_{\rep^-(Q)}$ has a unique regular component $\mathcal C$ which contains infinite dimensional representations. By Proposition
\ref{dualityregularcomp}(2), $\Ga$ is obtained from $\mathcal C$
by deleting the infinite dimensional representations. Thus $\Ga$ is
the unique regular component of $\Ga_{\rep^+\hspace{-0.5pt}(Q)}$.  The proof of the
theorem is completed.

{\center \section{Number of regular components}}

In case $Q$ is of infinite Dynkin type, as seen in the preceding section, the Auslander-Reiten quiver
$\Ga_{\rep^+\hspace{-1pt}(Q)}$ of $\rep^+\hspace{-1pt}(Q)$ has at most four connected components and at most two regular components.
The main objective of this section, to the contrary, is to show that $\Ga_{\rep^+\hspace{-1pt}(Q)}$ has
infinitely many regular components provided that $Q$ is not of finite or infinite Dynkin type.

\medskip

\begin{Lemma} \label{reg-reg} Let $\Sa$ be a finite subquiver of $Q$. If $M$ is a regular representation in $\rep(\hspace{-1pt}\it\Sigma)$, then it is a regular representation in $\rep^+\hspace{-1pt}(Q)$.
\end{Lemma}

\noindent{\it Proof.} Let $M$ be a regular representation in
$\Ga_{\rep(\hspace{-1.5pt}\it\Sigma)}$. It is well known that the
regular components of $\Ga_{\rep(\it\Sigma)}$ are stable tubes or of
shape $\Z\A_{\infty}$. Thus $\rep(\Sa)$ admits an infinite chain of
irreducible  epimorphisms \vspace{-4pt}
$$\xymatrixcolsep{15pt}\xymatrix{\cdots \ar[r] & M_i\ar[r]^-{f_i} & M_{i-1}\ar[r] & \cdots \ar[r] & M_1\ar[r]^-{f_1} & M}$$
and an infinite chain of irreducible monomorphisms \vspace{-3pt}
$$\xymatrixcolsep{15pt}\xymatrix{M \ar[r]^-{g_1} & N_1 \ar[r] & \cdots \ar[r] & N_{i-1}\ar[r]^-{g_i} & N_i\ar[r] & \cdots  .}$$
In particular, ${\rm Hom}(M_i, M)\ne 0$ and $\Hom(M, N_i)\ne 0$, for
every $i\ge 1$. By the first two statements of Proposition
\ref{no_ext}, $M$ is regular in $\rep^+\hspace{-1pt}(Q)$. The proof
of the lemma is completed.

\medskip

\begin{Lemma} \label{fin-supp} If $\,Q$ is infinite, then every regular component of $\Ga_{\rep^+\hspace{-1pt}(Q)}$ has at most finitely many representations supported by any given finite full subquiver of $Q$.\end{Lemma}

\noindent{\it Proof.} Let $Q$ be infinite with a finite full
subquiver $\Sa$. Assume that $\Ga$ is a regular component of
$\Ga_{\rep^+\hspace{-1pt}(Q)}\vspace{1pt}$ containing infinitely many
representations $M_j$ with $j\ge 0$ such that ${\rm supp}
\hspace{0.3pt} M_j\subseteq \Sa$. Since $\Sa$ is finite, we may
assume that ${\rm supp}\hspace{0.3pt} M_j=\Sa$ for all $j\ge 0$.
Setting $n_j$ to be the quasi-length of $M_j$, we deduce from
Proposition \ref{orbits}(2) that the $n_j$ are pairwise distinct. In
particular, we may assume that $n_j>0$ for all $j>0$. By Theorem
\ref{reg-cpt}, $\Ga$ is of shape $\Z\A_\infty$,
$\N\A_\infty$ or $\N^-\hspace{-1pt}\A_\infty$. Thus, for each $j>0$,
there exists in $\rep^+\hspace{-1pt}(Q)$ a chain of irreducible
epimorphisms of length $n_j-1$ from $M_j$ to a quasi-simple
representation $N_j$, and a chain of irreducible monomorphisms of
length $n_j-1$ from $\tau^{n_j-1}N_j$ to $M_j$. As a consequence,
the $N_j$ and the $\tau^{n_j-1}N_j$ with $j>0$ are all supported by
$\Sa$. Since the $n_j$ are pairwise distinct, the set ${\mathcal
Y}=\{N_j, \tau^{n_j-1} N_j \mid j\ge 0\}$ is infinite. Since $\Sa$
is finite, $\mathcal Y$ contains infinitely many quasi-simple
representations having the same support. This is contrary to
Proposition \ref{orbits}(2). The proof of the lemma is completed.

\medskip

Recall that a trivially valued translation quiver is called a {\it stable tube} of rank $n \hspace{0.5pt}(\hspace{1.5pt}>0)$ if it is of shape $\Z\A_\infty/\hspace{-3.5pt}<\hspace{-3.5pt}\tau^n\hspace{-3.5pt}>$, and a stable tube of rank one is called a {\it homogeneous tube}.
The following result seems to be well known. However, to the best of
our knowle\-dge, it is not explicitly stated anywhere. For this
reason, we include a proof which is suggested by Kerner.

\medskip

\begin{Prop} \label{fin-q-reg} Let $Q$ be a finite connected quiver without oriented cycles.

\vspace{-2pt}

\begin{enumerate}[$(1)$]

\item If $Q$ is of Euclidean type, then $\Ga_{\rep(Q)}$ has infinitely many homogeneous tubes.

\vspace{1pt}

\item If $Q$ is not of Dynkin type, then $\Ga_{\rep(Q)}$ has infinitely many
regular components.

\end{enumerate} \end{Prop}

\noindent{\it Proof.} (1) Let $Q$ be of Euclidean type. The regular
components of $\Ga_{\rep(Q)}$ are pairwise orthogonal stable tubes;
see \cite{DR}. Assume that $Q$ has only two vertices $a, b$. Then
$Q$ is the Kronecker quiver $K$ with exactly two arrows $\alpha,
\beta$ from $a$ to $b$. Note that the regular components of
$\Ga_{\rep(K)}$ are all homogeneous tubes; see \cite{DR}. Denote by
$\mathscr{P}$ the set of monic irreducible polynomials over $k$,
which is known to be infinite. For each $p\in \mathscr{P}$, define
$M_p\in \rep(K)$ by setting $M_p(a)=M_p(b)=k[x]/\hspace{-4pt} <\hspace{-3.5pt}p\hspace{-3.5pt}>$, $M_p(\alpha)=\id$, and $M_p(\beta)$
to be the multiplication by the class $\bar x\in
k[x]/\hspace{-3.5pt} <\hspace{-3pt}p\hspace{-2.5pt}>$. Having a
dimension vector of the form $(t, t)$, the representation $M_p$ is
regular; see \cite{DR}, and quasi-simple since $\End(M_p)\cong
K[x]/\hspace{-4pt} <\hspace{-3.5pt}p\hspace{-3.5pt}>$. Moreover, the
$M_p$ with $p\in \mathscr{P}$ are pairwise orthogonal, and hence lie in pairwise different components of
$\Ga_{\rep(K)}$.

\vspace{1pt}

Assume now that $Q$ has $n \,(\,>2)$ vertices and Statement (1)
holds for any Euclidean quiver of $n-1$ vertices. Then $\Ga_{\rep(Q)}$
has at least one stable tube $\mathcal T$ of rank $r >1$; see
\cite[Section 6]{DR}. Choose a quasi-simple representation $S$ in
$\mathcal T$. The perpendicular category $S^\perp$, that is
the full additive subcategory of $\rep(Q)$ gene\-rated by the
representations $L$ with $\Hom(S, L)=\Ext^1(S, L)=0$, is equivalent
to $\rep(Q\hspace{0.2pt}')$, where $Q\hspace{0.2pt}'$ is an Euclidean quiver of $n-1$ vertices;
see \cite[(10.1)]{GeL}. Thus, the
Auslander-Reiten quiver of $S^\perp$ has infinitely many homogenous
quasi-simple representations $S_i$ with $i\ge 1$. Since $\Ext^1(S_i,
S_i)\ne 0$, we see that $S_i$ lies in a stable tube $\mathcal{T}_i$
of $\Ga_{\rep(Q)}$. If $\mathcal{T}_i={\mathcal T}$ for some $i$ then, since $\End(S_i)$ is
divisible, the quasi-length of $S_i$ is at most $r$. Thus
we may assume $\mathcal{T}_i\ne {\mathcal T}$, for all $i\ge 1$. Consider, for each
$i \ge 1$, an almost split sequence
$$\eta_i: \;\; \xymatrixcolsep{14
pt}\xymatrix{0 \ar[r] & S_i \ar[r] & E_i \ar[r]^{g_i} & S_i \ar[r]
& 0}$$ in $S^\perp$. Let $f: L \to S_i$ be a  non-retraction
morphism in $\rep(Q)$ with $L\in\Ga_{\rep(Q)}\vspace{1pt}$. Then
$L$ lies in $\mathcal{T}_i$ or $L$ is preprojective. In the first
case, $L \in S^\perp$ since $\mathcal{T}_i\ne \mathcal{T}$, and
hence, $f$ factors through $g_i$. In the second case, consider a
pullback diagram \vspace{-2pt}
$$\xymatrixrowsep{12pt} \xymatrixcolsep{20pt}\xymatrix{
0 \ar[r] & S_i \ar[r] \ar@{=}[d] & E \ar[d] \ar[r] &  L \ar[d]^f \ar[r] & 0\\
0 \ar[r] & S_i \ar[r] & E_i \ar[r]^{g_i} & S_i \ar[r] & \;0}$$ in $\rep(Q)$.
Since $\Ext_{\rep(Q)}^1(L,S_i) \cong D \Hom_{\rep(Q)}(S_i, \tau L)=0$, the upper row splits. Thus $f$
factors through $g_i$. This shows that $\eta_i$ is an almost split sequence
in $\rep(Q)$. Consequently, the $\mathcal{T}_i$ with $i\ge 1$
are all homogenous tubes of $\Ga_{\rep(Q)}\vspace{1.5pt}$.

(2) Let $Q$ be of non-Dynkin type. By Statement (1), we may assume
that $Q$ is wild. Then every regular component of
$\Ga_{\rep^+\hspace{-1pt}(Q)}$ is of shape $\Z\A_\infty$; see \cite{R1}. Consider first the case where $Q$ has only two
vertices. We may assume that $Q$ consists of the
above-mentioned Kronecker quiver $K$ and possibly some extra arrows from
$a$ to $b$. Then $\rep(K)$ is a full additive
subcategory of $\rep(Q)$ generated by the representations
annihilated by the arrows other than $\alpha, \beta$. By Lemma
\ref{reg-reg}, the $M_p$ with $p\in \mathscr{P}$ are regular
representations in $\Ga_{\rep(Q)}$. Suppose that $M_q$ is not
quasi-simple in $\Ga_{\rep(Q)}$ for some $q\in \mathscr{P}$. Then
$\rep(Q)$ has an irreducible epimorphism $f: M_q\to N$. It is easy
to verify that $N$ is annihilated by the arrows other than $\alpha,
\beta$, and thus $f$ is an irreducible epimorphism in $\rep(K).$
This contradicts the fact that $M_q$ is quasi-simple in
$\Ga_{\rep(K)}$. Therefore, the $M_p$ with $p\in \mathscr{P}$ are
all quasi-simple in $\Ga_{\rep(Q)}$. Suppose on the contrary that
$\Ga_{\rep(Q)}$ has only finitely many regular components. Being
quasi-simple, the $M_p$ with $p\in \mathscr{P}$ are contained in
finitely many $\tau$-orbits of $\Ga_{\rep(Q)}$. Therefore,
$\mathscr{P}$ contains infinitely many $p_i$ with $i\ge 1$ such that
the $M_{p_i}$ lie in the the same $\tau$-orbit of $\Ga_{\rep(Q)}$. Write $M_i=M_{p_i}$, for all $i\ge 1$. Since $M_1$ is regular, there exists some $s>0$ such that $\Hom(M_1, \tau^iM_1)\ne 0$ for every $i\ge s$; see \cite[(1.3)]{Ker}. In particular, $M_j=\tau^{n_j}\hspace{-2pt}M_1$ with $n_j\le s$, for all $j\ge 1$. Then there exists some $r>1$ such that $n_r<-s$, that is, $-n_r>s$. This yields 
$$\Hom(M_r, M_1)=\Hom(\tau^{n_r}\hspace{-2pt}M_1, M_1)\cong \Hom(M_1,\, \tau^{-n_r}\hspace{-2pt}M_1)\ne 0, \vspace{-2pt}$$ contrary to the fact that the $M_p$ are pairwise orthogonal. Assume now that
$Q$ has $n \,(\,>2)$ vertices and Statement (2) holds for any non-Dynkin quiver of
$n-1$ vertices. Let $Q\hspace{0.2pt}'$ be a non-Dynkin connected full subquiver
of $Q$ with $n-1$ vertices. Then,
$\Ga_{\rep(Q\hspace{0.2pt}')}$ has infinitely many regular representations $N_i$
with $i\ge 1$. By Lemma \ref{reg-reg}, the $N_i$ are regular
representations in $\Ga_{\rep(Q)}$. Since the  $N_i$ are not sincere, they are
distributed in infinitely many regular components of
$\Ga_{\rep\hspace{-1pt}(Q)}$; see \cite[(1.3)]{Ker}. The proof of the proposition is
completed.

\medskip

\begin{Theo} \label{maintheonumber}
Let $Q$ be a connected strongly locally finite quiver. Then $\Ga_{\rep^+\hspace{-1pt}(Q)}\vspace{1pt}$ has only finitely many
regular components if and only if $Q$ is of finite or infinite Dynkin
type, and in this case, the number of regular components is at most two.

\end{Theo}

\noindent{\it Proof.} If $Q$ is of finite Dynkin type, then $\Ga_{\rep^+\hspace{-1pt}(Q)}$ has no regular component. If $Q$ is of infinite Dynkin type, then $\Ga_{\rep^+\hspace{-1pt}(Q)}\vspace{1pt}$ has at most two regular components by Theorems \ref{A-double-infity}, \ref{A-double-infity} and \ref{D-inf-main}. For proving the necessity, by Proposition \ref{fin-q-reg}, we only
need to consider the case where $Q$ is infinite but not of infinite
Dynkin type. Then $Q$ contains a connected finite full subquiver ${\it\Sigma}$ of non-Dynkin type.
By Proposition \ref{fin-q-reg},
${\Ga}_{\rep^+(\it\Sigma)}$ contains infinitely many regular
representations $M_i$ with $i\ge 1$. By Lemma \ref{reg-reg}, the
$M_i$ are regular representations in $\Ga_{\rep^+\hspace{-1pt}(Q)}\vspace{1pt}$. By
Lemma \ref{fin-supp}, they are distributed in infinitely many regular
components of $\Ga_{\rep^+\hspace{-1pt}(Q)}$. The proof of the theorem is
completed.

\medskip

To conclude this section, we shall show that any of the four types of
regular components may appear infinitely many times.

\medskip

\begin{Lemma}\label{supp-incl} Suppose that $Q$ is infinite and connected. Let $M$ be a representation in $\Ga_{\rep^+\hspace{-1pt}(Q)}$ such that
${\rm supp}\hspace{0.4pt}M  \subseteq {\rm
supp}(\tau^n\hspace{-1.5pt}M)$ for some $n>0$. If $Q$ has left
infinite paths, then $\tau^m\hspace{-1.5pt}M$ is pseudo-projective
for some $m\ge n$.
\end{Lemma}

\noindent{\it Proof.} Let $p$ be a left infinite path in $Q$. By
Proposition \ref{ARrep+}, $\tau^n\hspace{-1.5pt}M$ is finite
dimensional. Since $Q$ is connected, there exists a right infinite
acyclic walk $w$ which intersects ${\rm
supp}(\tau^n\hspace{-1.5pt}M)$ only at $s(w)$ and ends with the
inverse of a left infinite subpath of $p$. In particular, all but
finitely many edges in $w$ are inverses of arrows. Assume that $w$
starts with the inverse of an arrow, then $w^{-1}$ is an infinite
acyclic walk which ends with an arrow and intersects ${\rm
supp}(\tau^n\hspace{-1.5pt}M)$ only at $e(u^{-1})$. By Lemma
\ref{stableorbit}(2), ${\rm supp}\hspace{0.4pt}M$ contains some
vertex lying in $w^{-1}$ but different from $e(w^{-1})=S(w)$, which
is absurd since ${\rm supp}(\tau^n\hspace{-1.5pt}M) \supseteq {\rm
supp}\hspace{0.4pt}M$. Thus, $w$ starts with an arrow. By Lemma
\ref{stableorbit2}(2),  $\DTr^m\hspace{-1.5pt}M$ is infinite
dimensional for some $m>n$, that is, $\tau^{m-1}M$ is a
pseudo-projective representation in $\Ga$. The proof of the lemma is
completed.

\medskip

\begin{Theo} \label{inf-num-cpt}
Let $Q$ be a connected strongly locally finite quiver which is
infinite but not of infinite Dynkin type.

\begin{enumerate}[$(1)$]

\item If $Q$ has no infinite path, then $\Ga_{\rep^+\hspace{-1pt}(Q)}$ has infinitely many
regular components of shape
$\mathbb{Z}\mathbb{A}_{\infty}$.

\vspace{1pt}

\item If $Q$ has left infinite but no right infinite paths, then $\Ga_{\rep^+\hspace{-1pt}(Q)}$ has infinitely many
regular components of shape $\mathbb{N}\mathbb{A}_{\infty}$.
\vspace{1pt}

\item If $Q$ has right infinite but no left infinite paths, then $\Ga_{\rep^+\hspace{-1pt}(Q)}$ has infinitely many
regular components of shape $\mathbb{N}^-\hspace{-1.5pt}\mathbb{A}_{\infty}$.

\vspace{1pt}

\item If $Q$ has both left infinite paths and right infinite paths, then $\Ga_{\rep^+\hspace{-1pt}(Q)}$ has infinitely many regular components of wing type.

\end{enumerate}\end{Theo}

\noindent{\it Proof.} Statement (1) is an immediate consequence of
Theorem \ref{maintheonumber} and Corollary \ref{no-infpath-cpt}. For proving the
rest of the theorem, let ${\it\Sigma}$ be a wild finite
full subquiver of $Q$. Then ${\it
\Gamma}_{\rep({\it\Sigma})}$ contains infinitely many regular
representations $M_i$ with $i\ge 1$ such that $M_i$, $\tau_{_{\it\Sigma}}M_i$, and $\tau^-_{_{\it\Sigma}}M_i$ are all sincere representations of $\Sa$; see \cite[(1.3)]{Ker}. By Lemma \ref{reg-reg}, each $M_i$ lies in a regular component $\Ga_i$ of $\Ga_{\rep^+\hspace{-1pt}(Q)}\vspace{1pt}$. By Lemma \ref{fin-supp}, the $\Ga_i$ with $i\ge 1$ are pairwise distinct.

Assume that $Q$ has left infinite path. We claim that none of the
$M_i$ is left stable in  $\Ga_{\rep^+\hspace{-1pt}(Q)}$. Indeed,
let $i\ge 1$ be such that $\tau M_i\in \Ga_i.$ Then $\tau M_i$ is
finite dimensional. Choose $\Oa$ to be a finite connected full
subquiver of $Q$ containing the support of $M_i \hspace{-1pt}
\oplus\hspace{1pt}\tau \hspace{-1pt}M_i$. Then
$\tau_{_{\hspace{-1pt}\it\Omega}}\hspace{0.3pt} M_i = \tau M_i$.
Since $\Sa={\rm supp}\hspace{0.4pt}M_i\subseteq \Oa$, the path
algebra $k\Sa$ is a quotient of the path algebra $k\Oa$. By Lemma
5.2 stated in \cite[Chapter VIII]{ASS},
$\tau_{_{\hspace{-1pt}\it\Sigma}}\hspace{0.3pt}M_i$ is a
sub-representation of
$\tau_{_{\hspace{-1pt}\it\Omega}}\hspace{0.3pt}M_i$. This yields
\vspace{-2pt}
$${\rm supp} \hspace{0.4pt} M_i=
{\rm supp}(\tau_{_{\hspace{-1pt}\it \Sigma}}\hspace{0.3pt}M_i)
\subseteq {\rm
supp}(\tau_{_{\hspace{-1pt}\it\Omega}}\hspace{0.3pt}M_i)={\rm
supp} (\tau M_i).$$ By Lemma \ref{supp-incl}(2),
$\DTr^r\hspace{-1pt}M_i$ is infinite dimensional in $\rep(Q)$ for
some $r>0$. That is, $M_i$ is not left stable in
$\Ga_{\rep^+\hspace{-1pt}(Q)}$. This establishes our claim.
Dually, we may show that if $Q$ has right infinite paths, then
none of the $M_i$ is right stable in
$\Ga_{\rep^+\hspace{-1pt}(Q)}$.

Now suppose that $Q$ has left infinite but no right infinite
paths. Then $\Ga_{\rep^+\hspace{-1pt}(Q)}$ has no infinite
dimensional representation. By Theorem \ref{reg-cpt}, the $\Ga_i$
are all of shape $\Z\A_\infty$ or $\N\A_\infty$. By the above
claim, none of the $\Ga_i$ are left stable. Hence, the $\Ga_i$
with $i\ge 1$ are all of shape $\N\A_\infty$. This proves
Statement (2). Dually, Statement (3) holds true. Finally, suppose
that $Q$ has both left infinite paths and right infinite paths. As
seen above, each of the $M_i$ with $i\ge 1$ is neither left stable
nor right stable. By Theorem \ref{reg-cpt}(4), the $\Ga_i$ are all
finite wings. The proof of the theorem is completed.

\vspace{-5pt}

{\center \section{The bounded derived categories}}

The main objective of this section is to study the
Auslander-Reiten theory in $D^b(\rep^+\hspace{-1pt}(Q))$, the
derived category of the bounded complexes in
$\rep^+\hspace{-1pt}(Q)$. Making use of the previously obtained
results for $\rep^+\hspace{-1pt}(Q)$, we shall be able to give a
complete description of its Auslander-Reiten components of
$D^b(\rep^+\hspace{-1pt}(Q))$.

\medskip

We begin with an arbitrary hereditary abelian category
$\mathcal{H}$. Let $D^b(\mathcal{H})$ stand for the derived category
of the bounded complexes in $\mathcal{H}$. It is well known that
$D^b(\H)$ is a triangulated category whose translation functor is
the {\it shift functor} denoted by $[1].$ If $f: X\to Y$ is a
morphism in $D^b(\H)$ and $i\in \Z$, then we call $X[i]$ and $f[i]:
X[i]\to Y[i]$ the {\it shifts} by $i$ of $X$ and $f$, respectively.
As usual, we shall regard $\H$ as a full subcategory of $D^b(\H)$ by
identifying an object $X\in \H$ with the stalk complex $X[0]$. It is
important to observe that each object in $D^b(\mathcal H)$ is a
finite direct sum of the stalk complexes $X[i]$, where $X\in
\mathcal{H}$ and $i\in \Z$.  Moreover, for $X,Y \in \H$,
$\Hom_{D^b(\mathcal{H})}(X[i], Y[j])\ne 0$ only if $i\le j\le
i+1;\vspace{1pt}$ see \cite[(3.1)]{Lez}. Now, let
\vspace{-5pt}
$$\qquad\qquad \Delta: \quad\xymatrixcolsep{20pt}\xymatrix{X \ar[r]^-f & Y \ar[r]^-g & Z \ar[r]^-h &
X[1]}\vspace{-3pt}$$
be an exact triangle in $D^b(\H)$. One calls $X$ the {\it starting
term} and $Z$ the {\it ending term} of $\Delta$. We say that $\Da$ is an \emph{almost split triangle} if $f$ is
minimal left almost split and $g$ is minimal right almost split;
compare \cite{H}. For various equivalent conditions, we refer the
reader to \cite[(2.6)]{Kra}. One says that $D^b(\H)$ has left
(respectively, right) almost split triangles if every indecomposable
object in $D^b(\H)$ is the starting (respectively, ending) term of
an almost split triangle, and that $D^b(\H)$ has almost split
triangles if it has left and right almost split triangles.

\medskip

The following result tells us how the minimal almost split morphisms in $D^b(\H)$ are related to  those in $\H$.

\medskip

\begin{Lemma}\label{der-irr-mor} Let $X, Y,$ and $Z$ be objects in $\mathcal H$.

\begin{enumerate}[$(1)$]

\item If $\,(f, \eta)^T: X\to Y\oplus Z[1]$ is a minimal left almost split morphism in $D^b(\mathcal{H})$, then $f: X\to Y$ is a minimal left almost split morphism in $\mathcal{H}$.

\vspace{1pt}

\item If $\,(g, \zeta): Y\oplus Z[-1] \to X$ is a minimal right almost split morphism in $D^b(\mathcal{H})$, then $g: Y\to X$ is a minimal right almost split morphism in $\mathcal{H}$.

\vspace{1pt}

\item If $\,\xi: X\to Y[1]$ is an irreducible morphism in $D^b(\mathcal{H})$, then $X$ is injective and $Y$ is projective.

\end{enumerate}
\end{Lemma}

\noindent{\it Proof.} Let $\,(f, \eta)^T: X\to Y\oplus Z[1]$ be a
minimal left almost split morphism in $D^b(\mathcal{H})$. It is
evident that $f: X\to Y$ is left minimal and is not a section.
Suppose that $u: X\to M$ is a non-section morphism in $\mathcal H$.
Then $u$ factors through $(f, \eta)^T$ in $D^b(\mathcal{H})$. Since
$\Hom_{D^b(\mathcal{H})}(Z[1], M)=0\vspace{1.5pt}$, we see that $u$
factors through $f$ in $\mathcal H$. This proves Statement (1). In a
dual manner, we can establish Statement (2). Finally, suppose that
$\xi: X\to Y[1]$ is an irreducible morphism in $D^b(\mathcal{H})$.
Let
$$\xymatrixcolsep{18pt} \xymatrix{0\ar[r]& X \ar[r]^u & M \ar[r]^v & N \ar[r] & 0}\vspace{-2pt} $$
be a short exact sequence in $\mathcal H$. Since
$\Ext_{\mathcal{H}}^2(N, Y)=0\vspace{1.5pt}$, there exists some
morphism $\delta: M\to Y[1]$ in $D^b(\mathcal{H})$ such that $\xi=
\delta \circ u$.
Since $\Hom_{D^b(\mathcal{H})}(Y[1], M)=0\vspace{1pt}$, we see that
$\delta$ is not a retraction in $D^b(\mathcal{H})$. Thus $u$ is a
section in $D^b(\mathcal{H})$, and hence a section in $\mathcal H$.
Since $\mathcal H$ is abelian, this shows that $X$ is injective.
Dually, one can show that $Y$ is projective. The proof of the lemma
is completed.

\medskip

The following result relates the almost split triangles in $D^b(\H)$ to the almost split sequences in $\H$.
This was first established by Happel in the case where $\H$ is the category of finite dimensional representations of
a finite acyclic quiver; see \cite[(5.4)]{H}. Observe that our
approach  is very much different from Happel's.

\medskip

\begin{Theo}\label{ART} Let $\,\H \vspace{1pt}$ be a hereditary abelian category.

\vspace{-1pt}

\begin{enumerate}[$(1)$]

\item If $\,\xymatrixcolsep{15pt} \xymatrix{0\ar[r]& X \ar[r] & Y
\ar[r] & Z \ar[r] & 0}\vspace{-1pt}$ is an almost split sequence
in $\mathcal H$, then it induces an almost split triangle
$\xymatrixcolsep{16pt}\xymatrix{X \ar[r] & Y \ar[r] & Z\ar[r] &
X[1]}$ in $D^b(\mathcal{H})$.

\item If  $S$ is a simple object in $\H \vspace{1pt}$ with a
projective cover $P$ and an injective hull $I$, then
$D^b(\mathcal{H})$ has an almost split triangle as
follows$\,:$\vspace{-4pt}
    $$\xymatrixcolsep{15pt}\xymatrix{I\ar[r] & \left(I/S\right) \oplus
     \left({\rm rad} \hspace{0.5pt}P\right)[1]\ar[r] & P[1] \ar[r] & I[1].} \vspace{-4pt}$$

\item Every almost split triangle in
$D^b(\mathcal{H})\vspace{-1pt}$ is a shift of an almost split
triangle stated in the above two statements.

\end{enumerate} \end{Theo}

\noindent{\it Proof.} (1) Let $\eta:
\xymatrixrowsep{8pt}\xymatrixcolsep{14pt}\xymatrix{0\ar[r]& X
\ar[r]^-f & Y \ar[r]^-g & Z \ar[r] & 0}\vspace{-2pt}$ be an almost
split sequence in $\mathcal{H}$. Then it induces an exact triangle
$\Delta: \xymatrixcolsep{14pt}\xymatrix{X \ar[r]^-f & Y \ar[r]^-g &
Z \ar[r]^-\eta & X[1]}$ in $D^b(\mathcal{H})$. Since $g$ is right
minimal in $\mathcal H$, it is right minimal in $D^b(\mathcal{H})$.
We claim that each non-zero non-retraction morphism $\zeta: M\to Z$ in
$D^b(\mathcal{H})$ factors through $g$.
Indeed, we may assume that $\zeta=(h, \xi): M\oplus N[-1]\to Z,$
where $M, N\in \mathcal{H}$. Then $h$ is a non-retraction morphism
in $\mathcal H$, and hence it factors through $g$ in $\mathcal H$.
On the other hand, since $\Hom_{D^b(\mathcal H)}(L[-1], X[1])=0$, we
have $\eta \,\xi=0$, and thus $\xi$ factors through $g$ in
$D^b(\mathcal H)$. This proves that $g$ is minimal right almost
split in $D^b(\mathcal{H})$.
Hence, $\Delta$ is almost split; see \cite[(2.6)]{Kra}.

(2) Let $S$ be a simple object in $\mathcal H$ with a projective
cover $\varepsilon: P\to S$ and an injective hull $\iota: S\to I.$
In view of Lemma \ref{pc-ih}, we see that $P$ and $I$ are strongly
indecomposable. Setting  $h=\iota \,\varepsilon$, we get an exact
triangle \vspace{-2pt}
$$(*)\qquad \xymatrixrowsep{4pt}\xymatrixcolsep{18pt}\xymatrix{I[-1] \ar[r]^-f & M \ar[r]^-g & P \ar[r]^-h & I}\vspace{-5pt}$$
in $D^b(\mathcal{H})\vspace{1pt}$. Let  $\mu: X\to P$ be a
non-zero non-retraction morphism in $D^b(\mathcal{H})$. We may
assume that $\mu=(u, \delta): Y\oplus Z[-1]\to P$, where $Y, Z\in
\mathcal{H}$. Then $u: Y\to P$ is a non-retraction morphism in
$\mathcal{H}$. Thus $\varepsilon u=0$, and hence, $hu=0$. On the
other hand, since $\Hom_{D^b(\mathcal{H})}(Z[-1], I) \cong
\Ext^1_{\mathcal{H}}(Z,I)=0\vspace{1pt}$, we have $h \delta=0$.
This implies that $h\mu=0$, and consequently, $\mu$ factors
through $g$. Therefore, $g$ is right almost split in
$D^b(\mathcal{H})$. Since $I[-1]$ is strongly indecomposable,
$(*)$ is an almost split triangle in $D^b(\mathcal{H})$; see
\cite[(2.6)]{Kra}. Furthermore, we may assume that $M=N\,\oplus
\,L[-1]$, where $N, L\in \mathcal{H}$. Write $f= (f_1, \,
v[-1])^T$ and $g = (w, \, g_1)$, where $v: I\to L$ and $w: N\to P$
are morphisms in $\mathcal{H}$. By Lemma \ref{der-irr-mor}, $v$ is
minimal left almost split, and $w$ is minimal right almost split.
Hence, $N\cong{\rm rad}\,P$ and $L\cong I/S$.

\vspace{1pt}

(3) Let $\,\Delta: \xymatrixcolsep{14pt}\xymatrix{X \ar[r]^-f & Y
\ar[r]^-g & Z \ar[r]^-h & X[1]}$ be an almost split triangle in
$D^b(\mathcal{H})$. Up to a shift, we may assume that $X\in
\mathcal{H}.$ Since $g$ is right minimal, we may assume that
$Y=M\oplus N[1]$ with $M, N\in \mathcal{H}$. Write $f=(u, \zeta)^T:
X\to M\oplus N[1]\vspace{1pt}$ with $u: X\to M$ a morphism in
$\mathcal{H}$, and $g=(\xi, v): M\oplus N[1]\to Z$. By Lemma
\ref{der-irr-mor}(1), $u: X\to M$ is minimal left almost split in
$\mathcal{H}$.

Consider first the case where $X$ is not injective in $\mathcal{H}$.
Then $\mathcal{H}$ has a non-section monomorphism $X\to L$, which
factors through $u$. In particular, $u$ is a minimal left almost
split monomorphism in $\mathcal{H}$. It is then well known that
$\mathcal H$ has an almost split sequence $\,\zeta: \,
\xymatrixcolsep{15pt} \xymatrix{0\ar[r]& X \ar[r]^u& M\ar[r]^w &
N\ar[r] & 0}\vspace{-3pt}$; see \cite[(2.13), (2.14)]{AuR}. By
Statement (1), $\xymatrixcolsep{15pt} \xymatrix{X \ar[r]^u&
M\ar[r]^w & N\ar[r]^-\zeta & X[1]}$ is an almost split triangle
$D^b(\mathcal{H})$, which is isomorphic to $\Delta$. In other words,
$\Delta$ is of the form as stated in Statement (1).

Consider next the case where $X=I$, an injective object in $\mathcal{H}$.
Then $u: I\to M$ is a minimal left almost split epimorphism in
$\mathcal{H}$. Let $q: S\to I$ be the kernel of $u$. By Lemma \ref{pc-ih}(2), $S$ is
simple with $q$ being its injective hull, and $M \cong I/S$. Suppose that $Z\in \mathcal{H}$. Then $v=0$ and
$\xi u=-v\zeta=0$. Since $u: X\to M$ is an epimorphism in
$\mathcal{H}$, we get $\xi=0$, and hence $g=0$. As a consequence,
$h: Z\to  I[1]$ is a section, which is impossible since
$\Hom_{D^b(\mathcal{H})}(I[1], Z)=0\vspace{1pt}$. This shows that
$Z\not\in \H$. Since $h\ne 0$ and $Z$ is indecomposable, $Z=P[1]$
for some $P\in \mathcal{H}$. Then, $h=s[1]$ and $v=j[1]$, where $s:
P\to I$ and $j: N\to P$ are morphisms in $\mathcal H$. Now
\vspace{-6pt}$$\,\Delta[-1]: \; \xymatrixcolsep{24pt}\xymatrix{I[-1]
\ar[r]^-{f[-1]} & Y[-1] \ar[r]^-{g[-1]} & P \ar[r]^-s &
I,}\vspace{-3pt}$$ is an almost split triangle in $D^b(\H)$, where
$g[-1]=(\xi[-1], j): M[-1]\oplus N \to P$. By Lemma
\ref{der-irr-mor}(2), $j: N\to P$ is minimal right almost split in
$\mathcal{H}$. If $P$ is not projective, then we can show that
$\Delta[-1]$ is isomorphic to an almost split triangle induced from
an almost split sequence in $\mathcal H$ ending with $P$. In
particular, $I[-1]$ is isomorphic to an object in $\mathcal H$,
which is absurd. Thus $P$ is projective. Therefore, $j: N\to P$ is a
minimal right almost split monomorphism. Let $\varepsilon: P\to T$
be the cokernel of $j$. By Lemma \ref{pc-ih}(1), $T$ is simple with
$\varepsilon$ being its projective cover and $N\cong {\rm
rad}\hspace{0.4pt}P$. Moreover, since $s \circ g[-1]=0$ and $f s=0$,
we have $s j=0$ and $u s =0$. This yields a factorization $s=q p
\,\varepsilon$, where $p: T\to S$ is a non-zero morphism in $\mathcal
H$. Since $T, S$ are simple, $p$ is an isomorphism. That is,
$\Delta$ is of the form as stated in Statement (2). The proof of the
theorem is completed.

\medskip

Combining Theorem \ref{ART} and Corollary \ref{AR-ab} yields
immediately the following consequence. This is, in the Hom-finite
case, a result of Reiten and Van Den Bergh, which is stated
without a complete proof in \cite[(I.3.2)]{RVDB}.

\medskip

\begin{Cor} \label{der-arc} If $\H$ is a hereditary abelian
category, then

\vspace{-2pt}

\begin{enumerate}[$(1)$]

\item $D^b(\H)$ has left almost split triangles if and only if
$\H$ is left Auslander-Reiten  and the socle of any indecomposable
injective object has a projective cover$\,;$

\vspace{1pt}

\item $D^b(\H)$ has right almost split triangles if and only if
$\H$ is right Auslander-Reiten and the top of any indecomposable
projective object has an injective hull.

\end{enumerate} \end{Cor}

\medskip

From now on, we shall specialize our previous results to the case
where $\mathcal H$ is an abelian full subcategory of $\rep(Q)$.
First of all, combining Theorems \ref{AR-sequence} and \ref{ART},
we get immediately the following description of certain almost
split triangles in $D^b(\rep(Q))$.

\medskip

\begin{Theo} \label{Der-ART} Let $Q$ be a strongly locally finite quiver, and let $M$ be an indecomposable representation in $\rep(Q)$.

\begin{enumerate}[$(1)$]

\vspace{-1pt}

\item If $M$ is a non-projective object in
$\rep^+\hspace{-1pt}(Q)$, then
$D^b\hspace{-1pt}(\rep(Q)\hspace{-1pt})$ has an almost split
triangle $\xymatrixcolsep{15pt}\xymatrix{\DTr M \ar[r] & N \ar[r]
&  M \ar[r] & (\DTr\hspace{0.5pt}M)[1].}$

\vspace{0pt}

\item  If $M$ is a non-projective object in $\rep^-(Q)$, then
$D^b\hspace{-1pt}(\rep(Q)\hspace{-1pt})$ has an almost split
triangle $\xymatrixcolsep{15pt}\xymatrix{M \ar[r] & N \ar[r] &
\TrD M \ar[r] & M[1].}$

\vspace{0pt}

\item If $x$ is a vertex in $Q$, then $D^b\hspace{-1pt}(\rep(Q)\hspace{-1pt})$ has an almost split triangle
\vspace{-1pt}
$$\xymatrixcolsep{15pt}\xymatrix{I_x\ar[r] & I_x/S_x \oplus\,
 \left({\rm rad} \hspace{0.5pt}P_x\right)[1] \ar[r] & P_x[1]  \ar[r] & I_x[1].}\vspace{-3pt}$$

\end{enumerate}\end{Theo}

\medskip

The rest of the section is devoted to our main objective, that is,
to study the Auslander-Reiten theory in
$D^b(\rep^+\hspace{-1pt}(Q))$. We begin with a complete description
of its almost split triangles.

\medskip

\begin{Theo} \label{der-art} Let $Q$ be a strongly locally finite quiver,
and let $M$ be an indecomposable representation in $\rep^+\hspace{-1pt}(Q)$.

\vspace{0pt}

\begin{enumerate}[$(1)$]

\item If $M$ is neither projective nor pseudo-projective, then
$D^b\hspace{-1pt}(\rep^+\hspace{-1pt}(Q)\hspace{-1pt})$ has an
almost split triangle $\xymatrixcolsep{15pt}\xymatrix{\DTr M
\ar[r] & N \ar[r] & M \ar[r] & (\DTr\hspace{0.5pt} M)[1]}.$

\vspace{0pt}

\item If $M$ is finite dimensional and not injective, then
$D^b\hspace{-1pt}(\rep^+\hspace{-1pt}(Q)\hspace{-1pt})$ has an
almost split triangle $\xymatrixcolsep{15pt}\xymatrix{M \ar[r] & N
\ar[r]& \TrD M \ar[r] & M[1]}$.

\vspace{0pt}

\item  If $x$ is a vertex in $Q^+$, then
$D^b\hspace{-1pt}(\rep^+\hspace{-1pt}(Q)\hspace{-1pt})$ has an
almost split triangle \vspace{-2pt}
$$\xymatrixcolsep{15pt}\xymatrix{I_x\ar[r] & I_x/S_x \oplus\,
\left({\rm rad} \hspace{0.5pt}P_x\right)[1] \ar[r] & P_x[1] \ar[r]
& I_x[1].}\vspace{-4pt}$$

\item Every almost split triangle in $D^b(\rep^+\hspace{-1pt}(Q))$
is a shift of an almost split triangle stated in the above three
statements.

\vspace{1pt}

\end{enumerate}\end{Theo}

\noindent{\it Proof.} The first three statements follow
immediately from Theorems \ref{AR-sequence} and \ref{ART}. Now
consider an almost split triangle $\Delta:
\xymatrixcolsep{15pt}\xymatrix{L \ar[r] & M \ar[r]& N \ar[r] &
L[1]}$ in $D^b(\rep^+\hspace{-1pt}(Q))$. Up to a shift, we may
assume that $L\in \rep^+\hspace{-1pt}(Q)$. If $\Delta$ is induced
from an almost split sequence in $\rep^+\hspace{-1pt}(Q)$, then it
is of the form stated in Statement (1) or (2). Otherwise, by
Theorem \ref{ART}, $N=P_x[1]$ for some $x\in Q_0$ and $L$ is the
injective hull of $S_x$ in $\rep^+(Q)$. By Proposition
\ref{cor1.14shiping}(3), $x\in Q^+$,  and hence $L\cong I_x$. That
is, $\Delta$ is of the form as stated in Statement (3). The proof
of the theorem is completed.

\medskip

Next, we want to describe the irreducible morphisms in
$D^b(\rep^+\hspace{-1pt}(Q))$. Being non-zero, they are of the
form $f: M[i]\to N[i]$ or $\zeta: M[i]\to N[i+1]$, where $M, N\in
\rep^+\hspace{-1pt}(Q)$ and $i\in \Z$. It is evident that we need
only to study the irreducible morphisms of the second kind.

\medskip

\begin{Lemma}\label{der+irr-mor} Let $M, N$ be representations in $\rep^+\hspace{-1pt}(Q)$. If $M$ is indecompo\-sable, then
$D^b(\rep^+\hspace{-1pt}(Q))\vspace{1pt}$ has an irreducible morphism
$\eta: M\to N[1]$ if and only if there exists some  $x\in Q^+$ such that $M\cong I_x$ and
$N$ is a direct summand of $\,{\rm rad}\hspace{0.4pt} P_x$. \end{Lemma}

\noindent{\it Proof.} Suppose that $M$ is indecomposable. The
sufficiency follows easily from Theorem \ref{der-art}(3). Let
$\eta: M\to N[1]$ be an irreducible morphism in
$D^b(\rep^+\hspace{-1pt}(Q))\vspace{1pt}$. We claim that $M$ is
finite dimensional. Indeed, by Lemma \ref{der-irr-mor}(3) and
Proposition \ref{cor1.14shiping},  $N$ has some $P_y$ with $y\in
Q_0$ as a direct summand. Then $D^b(\rep^+\hspace{-1pt}(Q))$ has
an irreducible morphism $\zeta: M\to P_y[1]$; see
\cite[(3.2)]{Bau}. By Theorem \ref{Der-ART}(3),  $D^b(\rep(Q))$
has an almost split triangle \vspace{-3pt}
$$
\xymatrixcolsep{20pt}\xymatrix{I_y\ar[r] & I_y/S_y  \oplus \left({\rm rad \hskip 0.5pt}P_y\right)[1] \ar[r]^-{(\theta, \,q)} & P_y[1] \ar[r] & I_y[1].}\vspace{-2pt}$$
Since $\zeta$ is not a retraction in $D^b(\rep(Q))$, we have $\zeta=\theta f + q \xi,$ where
$f: M \to I_y/S_y$ is a  morphism in $\rep(Q)$ and $\xi: M\to \left({\rm rad}\hspace{0.4pt}P_y\right)[1]$ is a morphism in $D^b(\rep(Q))$. The composition of $\theta$ and $f$ is given by a pullback diagram

\vspace{-5pt}
$$\xymatrixrowsep{16pt} \xymatrixcolsep{20pt}\xymatrix@1{
\theta f: &0 \ar[r] & P_y \ar[r] \ar@{=}[d] & \, L \ar[d] \ar[r] & \,M \ar[d]^f \ar[r] & \,0\\
\theta: &0 \ar[r] & P_y \ar[r] & \,U \ar[r] & \,I_y/S_y \ar[r] &
\,0}$$ in $\rep(Q)$, where $L\in \rep^+\hspace{-1pt}(Q)$. Let $\Sa$
be the successor-closed subquiver of $Q$ generated by the
vertices in ${\rm supp}\hspace{0.4pt}L$. Restricting the preceding
pullback diagram to $\Sa$ yields a pullback diagram
$$\xymatrixrowsep{18pt} \xymatrixcolsep{20pt}\xymatrix{
\theta f: &0 \ar[r] & P_y \ar[r] \ar@{=}[d] & \,L \ar[d] \ar[r] & \,M \ar[d]^{f_{_{\hspace{-1pt}\it \Sigma}}} \ar[r] & \,0\\
\theta_{_{\hspace{-0.5pt}\it \Sigma}}: &0 \ar[r] & P_y \ar[r] &
\,U_{_{\hspace{-1pt}\it \Sigma}} \ar[r] &
\,(I_y/S_y)_{_{\hspace{-1pt}\it\Sigma}} \ar[r] &
\,0}\vspace{2pt}$$ in $\rep^+\hspace{-1pt}(Q)$. That is, we have a
factorization $\zeta = (\theta_{_{\hspace{-0.5pt}\it \Sigma}},
\,q) \, (f_{_{\hspace{-1pt}\it \Sigma}}, \,\xi)^T \vspace{1pt}$ in
$D^b(\rep^+\hspace{-1pt}(Q))$. Since
$(\theta_{_{\hspace{-0.5pt}\it \Sigma}}, \, q):
(I_y/S_y)_{_{\hspace{-1pt}\it \Sigma}} \oplus {\rm
rad}\hspace{0.4pt} P_y[1]\to P_y[1] \vspace{1.5pt}$ is clearly not
a retraction, the morphism $(f_{_{\hspace{-1pt}\it \Sigma}},
\,\xi)^T: M\to (I_y/S_y)_{_{\hspace{-1pt}\it \Sigma}} \oplus ({\rm
rad}\hspace{0.4pt} P_y)[1]$ is a section. Then,
$f_{_{\hspace{-0.4pt}\it \Sigma}}$ is a section, and hence $M$ is
a direct summand of $(I_y/S_y)_{_{\hspace{-0.4pt}\it \Sigma}}$. On
the other hand,  since $\Sa$ is top-finite by Lemma
\ref{Lemma1.6_shiping}, $(I_y/S_y)_{_{\hspace{-0.4pt}\it \Sigma}}$
is finite dimensional, and so is $M$. This establishes our claim.
Then, by Proposition \ref{cor1.14shiping} and Lemma
\ref{der-irr-mor}(3), we may assume that $M=I_x$ with $x\in Q^+$.
By Theorem \ref{der-art}(3), $D^b(\rep^+\hspace{-1pt}(Q))$ has an
almost split triangle \vspace{-1pt}
$$
\xymatrixcolsep{20pt}\xymatrix{I_x\ar[r] & I_x/S_x  \oplus
\left({\rm rad \hskip 0.5pt}P_x\right)[1] \ar[r] & P_x[1] \ar[r] &
I_x[1].}\vspace{-2pt}$$ Since $\eta: I_x\to N[1]$ is an
irreducible morphism in $D^b(\rep^+\hspace{-1pt}(Q))$, one has a
retraction $(\xi, w): I_x/S_x  \oplus \left({\rm rad \hskip
0.5pt}P_x\right)[1] \to N[1]$. Since
$\Hom_{D^b(\rep^+\hspace{-1pt}(Q))}(N[1], \, I_x/S_x)=0$, we see
that $w: \left({\rm rad \hskip 0.5pt}P_x\right)[1] \to N[1]$ is a
retraction. As a consequence, $N$ is a direct summand of ${\rm
rad}\hspace{0.4pt}P_x$. The proof of the lemma is completed.

\medskip

We are ready to describe the Auslander-quiver
$\Ga_{D^b(\rep^+\hspace{-1pt}(Q))}\vspace{1pt}$ of
$D^b(\rep^+\hspace{-1pt}(Q))$. Since $\rep^+\hspace{-1pt}(Q)$ is
hereditary and abelian, the vertices in
$\Ga_{D^b(\rep^+\hspace{-1pt}(Q))}\vspace{1.5pt}$ can be chosen to
be the complexes of the form $M[i]$, and then the arrows are of
the form $M[i]\to N[i]$ or $M[i]\to N[i+1]$, where $M, N\in
\Ga_{\rep^+\hspace{-0.5pt}(Q)}$ and $i\in \Z$.

\medskip

\begin{Lemma} \label{proparrowsDQ}
Let $Q$ be a strongly locally finite quiver.
\begin{enumerate}[$(1)$]

\item If $\alpha: M\to N$ is an arrow in
$\Ga_{\rep^+\hspace{-0.5pt}(Q)}$, then it is also an arrow in
$\Ga_{D^b(\rep^+\hspace{-1pt}(Q))}$.

\item If $\beta: x\to y$ is an arrow in $Q$, where $x\in Q^+$,
then it induces an arrow $[\beta]: I_x\to P_y[1]$ in
$\Ga_{D^b(\rep^+\hspace{-1pt}(Q))}\vspace{1.5pt}$.

\item Each arrow $\Ga_{D^b(\rep^+\hspace{-1pt}(Q))}\vspace{1.5pt}$
is a shift of an arrow stated in the above two statements. In
particular, $\Ga_{D^b(\rep^+\hspace{-1pt}(Q))}\vspace{1.5pt}$ has
a symmetric valuation.

\end{enumerate} \end{Lemma}

\noindent {\it Proof.} First of all, if $X$ is a representation in
$\Ga_{\rep^+\hspace{-0.5pt}(Q)}\vspace{1pt}$, then  ${\rm
End}_{\hspace{0.5pt}\rep^+\hspace{-1pt}(Q)}(X)\vspace{1pt}$ and
${\rm End}_{D^b(\rep^+\hspace{-1pt}(Q))}(X)\vspace{1.5pt}$ have
the same residue algebra which we denote by
$k_{_X}\vspace{1.5pt}$.

(1) Let $M, N\in \Ga_{\rep^+\hspace{-0.5pt}(Q)}\vspace{1pt}$ such that
the number of arrows from $M$ to $N$ in
$\Ga_{\rep^+\hspace{-0.5pt}(Q)}$ is $d_{_{MN}} \,\hspace{0.5pt} >
\hspace{-1pt} 0.$ Since $\Hom_{D^b(\rep^+\hspace{-1pt}(Q))}(X[1],
N)=0\vspace{1pt}$ for any representation $X$ in
$\rep^+\hspace{-1pt}(Q)$, we have ${\rm
rad}^2_{\rep^+\hspace{-1pt}(Q)}(M, N)={\rm
rad}^2_{D^b(\rep^+\hspace{-1pt}(Q))}(M, N)\vspace{1.5pt}$, and
consequently, ${\rm irr}_{D^b(\rep^+\hspace{-1pt}(Q))}(M,N) = {\rm
irr}_{\rep^+\hspace{-1pt}(Q)}(M,N).$ This yields \vspace{-2pt}
$$\begin{array}{rcl}
{\rm dim}_{k_M}\, {\rm irr}_{D^b(\rep^+\hspace{-1pt}(Q))}(M,N) & = & {\rm dim}_{k_M}\,  {\rm irr}_{\rep^+\hspace{-1pt}(Q)}(M,N) \vspace{1pt} \\
&=& {\rm dim}_{k_N}\, {\rm irr}_{\rep^+\hspace{-1pt}(Q)}(M,N) \vspace{1pt} \\
& = & {\rm dim}_{k_N}\,  {\rm irr}_{D^b(\rep^+\hspace{-1pt}(Q))}(M,N), \vspace{-1pt} \\
\end{array}$$
from which we see that the valued arrow $M\to N$ in
$\Ga_{D^b(\rep^+\hspace{-1pt}(Q))}$ has a symmetric valuation
$(d_{_{MN}}, d_{_{MN}})$, and it is replaced by the $d_{_{MN}}$
unvalued arrows in $\Ga_{\rep^+\hspace{-0.5pt}(Q)}$ from $M$ to $N$.

(2) Let $x\in Q^+$ and $y\in Q_0$ be such that the number of
arrows in $Q$ from $x$ to $y$ is $d_{xy}>0$. Then $P_y$ is a
direct summand of ${\rm rad}\hspace{0.4pt}P_x$. By Lemma
\ref{der+irr-mor}(3), $\Ga_{D^b(\rep^+\hspace{-1pt}(Q))}$ has a
valued arrow $\alpha_{xy}: I_x\to P_y[1]$. Since $k_{I_x} \cong k
\cong k_{P_y[1]}$ by Proposition \ref{prop1.3_shiping}, we see
that $\alpha_{xy}$ has a symmetric valuation $(d, d)$. On the
other hand, $d_{xy}$ is the maximal integer such that
$P_y^{d_{xy}}$ is a direct summand of ${\rm
rad}\hspace{0.4pt}P_x$. By Lemma \ref{der+irr-mor}(3), $d_{xy}$ is
the maximal integer such that $D^b(\rep^+\hspace{-1pt}(Q))$ has an
irreducible morphism $\zeta: I_y\to (P_y^{d_{xy}})[1]$. Hence
$d=d_{xy}$. Therefore, the valued arrow $\alpha_{xy}: I_x\to
P_y[1]$ is replaced by $d_{xy}$ unvalued arrows from $I_x$ to
$P_y[1]$, which are indexed by the arrows in $Q$ from $x$ to $y$.

Finally, every arrow in
$\Ga_{D^b(\rep^+\hspace{-1pt}(Q))}\vspace{1pt}$ is a shift of an
arrow $\gamma$ which is of the form $M\to N$ or $M\to N[1]$, where
$M, N\in \Ga_{\rep^+\hspace{-0.5pt}(Q)}$. By Lemma
\ref{der+irr-mor}, $\gamma$ is as stated in (1) or (2). The proof
of the lemma is completed.

\medskip

As an immediate consequence of Theorem \ref{der-art}, the
following result describes the Auslander-Reiten translation of
$\Ga_{D^b(\rep^+\hspace{-1pt}(Q))}$ which we write as $\tau_{_{\hspace{-0.4pt} D}}$.

\medskip

\begin{Lemma} \label{der+AR-trans} If $M$ is a representation lying in
$\Ga_{\rep^+\hspace{-0.5pt}(Q)}$, then

\vspace{-2pt}

\begin{enumerate}[$(1)$]

\item $\taud M$ is defined if and only if either $\tau M$ is
defined with $\tau_{_{\hspace{-0.4pt} D}} M=\tau M$, or $M=P_x$
for some $x\in Q^+$ with $\tau_{_{\hspace{-0.4pt} D}} M=I_x[-1]$.

\vspace{1pt}

\item $\tau_{_{\hspace{-0.4pt} D}}^-M$ is defined if and only if
either $\tau^-M$ is defined with $\tau_{_{\hspace{-0.4pt}
D}}^-M=\tau^-M$ or $M=I_x$ for some $x\in Q^+$ with
$\tau_{_{\hspace{-0.4pt} D}}^-M=P_x[1]$.

\end{enumerate}\end{Lemma}

\medskip

If $Q$ is connected, then $\Ga_{\rep^+\hspace{-0.5pt}(Q)}$ has a
unique preprojective component $\mathcal{P}_Q$. By Lemmas
\ref{proparrowsDQ} and \ref{der+AR-trans},
$\Ga_{D^b(\rep^+\hspace{-1pt}(Q))} \vspace{1pt}$ has a connected
component $\mathcal{C}_Q$ which is obtained by gluing
$\mathcal{P}_Q$ together with the shifts by $-1$ of the preinjective
components of $\Ga_{\rep^+\hspace{-0.5pt}(Q)}$ in the following way:
for each pair $(x, y)$ with $x\in Q^+$ and $y$ an immediate
successor of $x$ in $Q$, one draws $d_{xy}$ arrows from $I_x[-1]$ to
$P_y$, where $d_{xy}$ is the number of arrows from $x$ to $y$ in
$Q$. We call $\mathcal{C}_Q$ the {\it connecting component} of
$\Ga_{D^b(\rep^+\hspace{-1pt}(Q))}$ and describe its shape in the
following result.

\medskip

\begin{Prop}\label{conn} Let $Q$ be a connected strongly locally finite
quiver. Then the connecting component $\mathcal{C}_Q$ of
$\Ga_{D^b(\rep^+\hspace{-1pt}(Q)}$ embeds in $\mathbb{Z}Q^{\rm
\,op}$. Furthermore,

\vspace{-1.5pt}

\begin{enumerate}[$(1)$]

\item if $Q$ has no infinite path, then $\mathcal{C}_Q$ is of
shape $\mathbb{Z} Q^{\rm\, op};$

\vspace{1pt}

\item if $Q$ has right infinite but no left infinite paths, then
$\mathcal{C}_Q$ is of shape $\N^-\hspace{-2pt}\Da$ with $\Da$ the
right-most section of the preprojective component of
$\Ga_{\rep^+\hspace{-0.5pt}(Q)};\vspace{1pt}$

\vspace{0pt}

\item if $Q$ has left infinite but no right infinite paths, then
$\mathcal{C}_Q$ is isomorphic to a right stable translation
subquiver of  $\mathbb{Z} Q^{\rm\, op}$.

\end{enumerate}
\end{Prop}

\noindent{\it Proof.} By Theorem \ref{prep-component}, the
preprojective component $\mathcal{P}_Q$ of
$\Ga_{\rep^+\hspace{-0.5pt}(Q)}$ has a left-most section $P_Q$
generated by the $P_x$ with $x\in Q_0$ and isomorphic to $Q^{\rm
\, op}$. If $M$ is a preinjective representation in
$\Ga_{\rep^+\hspace{-0.5pt}(Q)}$, then it lies in the $\tau$-orbit
of some injective representation $I_y$ with $y\in Q^+$. By Lemma
\ref{der+AR-trans}(2), $M[-1]$ lies in the
$\tau_{_{\hspace{-0.4pt} D}}$-orbit of $P_y$. This shows that
$P_Q$ is a section of $\mathcal{C}_Q$, and consequently,
$\mathcal{C}_Q$ embeds in $\Z Q^{\hspace{0.4pt} \rm op}$.

(1) Suppose that $Q$ has no infinite path. By Theorem
\ref{preinjcomp}, $\Ga_{\rep^+\hspace{-0.5pt}(Q)}$ has a unique
preinjective component $\mathcal{I}$ of shape
$\N^-\hspace{-2pt}Q^{\rm \,op}$, while $\mathcal{P}_Q$ is of shape
$\N Q^{\rm \,op}$ by Theorem \ref{prep-component}(1). Therefore,
$\mathcal{C}_Q\cong \mathbb{Z} Q^{\rm\, op}\vspace{1pt}$.

(2) Suppose that $Q$ has right infinite but no left infinite
paths.  By Theorem \ref{preinjcomp},
$\Ga_{\rep^+\hspace{-0.5pt}(Q)}$ has a unique preinjective
component $\mathcal{I}_Q$ of shape $\N^-\hspace{-2pt}Q^{\rm
\,op}$, and $\mathcal{P}_Q$ has a right-most section $\Da$ by
Theorem \ref{prep-component}(2). Hence $\mathcal{C}_Q$ is of shape
$\N^-\hspace{-2pt}\Da$.

(3) Suppose that $Q$ has left infinite but no right infinite paths.
By Theorem \ref{prep-component}, the $P_x$ with $x\in Q_0$ are right
stable in $\Ga_{\rep^+\hspace{-0.5pt}(Q)}\vspace{1pt}$, and hence
right stable in $\Ga_{D^b(\rep^+\hspace{-0.5pt}(Q))}\vspace{1pt}$.
Containing a section of right stable vertices, $\mathcal{C}_Q$ is
right stable as a translation quiver. Thus, $\mathcal{C}_Q$ is
isomorphic to a right stable translation subquiver of $\mathbb{Z}
Q^{\rm\, op}$. The proof of the theorem is completed.

\medskip

By Lemmas \ref{proparrowsDQ} and \ref{der+AR-trans}, a regular
component of $\Ga_{\rep^+\hspace{-0.5pt}(Q)}$ remains to be a
connected component of $\Ga_{D^b(\rep^+\hspace{-1pt}(Q))}$. This
yields immediately the following result.

\medskip

\begin{Theo} \label{Der-cpt} If $Q$ is a connected strongly locally finite quiver, then the
connected components of $\Ga_{D^b(\rep^+\hspace{-1pt}(Q)}$ are the
shifts of the regular components of
$\Ga_{\rep^+\hspace{-0.5pt}(Q)}$ and those of the connecting
component $\mathcal{C}_Q$.
\end{Theo}

\medskip

\noindent{\sc Remark.} The shapes of the connected components of $\Ga_{D^b(\rep^+\hspace{-1pt}(Q)}$
are described by Theorem \ref{reg-cpt} and Proposition \ref{conn}.

\medskip

\noindent {\sc Example.}  Let $Q$ be the following quiver
$$ \footnotesize
\xymatrixrowsep{12pt}\xymatrixcolsep{18pt}\xymatrix{
&&& 0 \ar[d] &&&\\
\cdots \ar[r] & -2 \ar[r] & -1 \ar[r] & 1 \ar[r] & 2 \ar[r] & 3
\ar[r] & \cdots } \vspace{3pt} $$ Since the representations $P_x$
with $x\in Q_0$ are all infinite dimensional, the preprojective
component of $\Ga_{\rep^+\hspace{-0.5pt}(Q)}\vspace{1pt}$ is of
shape $Q^{\rm \,op}$. On the other hand, since $Q^+ = \{0\}$, we
see that $\Ga_{\rep^+\hspace{-0.5pt}(Q)}\vspace{1pt}$ has a unique
preinjective component $\{I_0\}$. Therefore, the connecting
component of $\Ga_{D^b(\rep^+\hspace{-1pt}(Q))}$ is as
follows$\,:$

\vspace{-3pt}

$$\hspace{50pt}\footnotesize
\xymatrixrowsep{12pt}\xymatrixcolsep{18pt}\xymatrix{
&& & I_0[-1] \ar[dr] \ar@{<.}[rr]  & & \, P_0   & &&&\\
& \cdots \ar[r] & P_3 \ar[r] & P_2 \ar[r] &P_1 \ar[r] \ar[ur]  &
P_{-1} \ar[r] & P_{-2}\ar[r] & \cdots } \vspace{4pt}$$
Moreover, by Theorem \ref{inf-num-cpt}(4),
$\Ga_{\rep^+\hspace{-0.5pt}(Q)}$ has infinitely many regular
components of wing type, and so does
$\Ga_{D^b(\rep^+\hspace{-1pt}(Q))}$.

\medskip

In case $Q$ is finite, Happel's result says that
$D^b(\rep^+\hspace{-1pt}(Q))$ has almost split triangles; see
\cite{H}. In the infinite case, we are able to find the precise
conditions on $Q$ such that $D^b(\rep^+\hspace{-1pt}(Q))$ has
(left, right) almost split triangles.

\medskip

\begin{Theo} \label{fin-theo} If $Q$ is a strongly locally finite quiver,
then $D^b(\rep^+\hspace{-1pt}(Q))$ has $($left, right$\hspace{0.4pt})$
almost split triangles if and only if $Q$ has no
$($right, left$\hspace{0.4pt})$ infinite path.
\end{Theo}

\noindent{\it Proof.} Firstly, we show that
$D^b(\rep^+\hspace{-1pt}(Q))$ has  left  almost split triangles if
and only if $Q$ has no right infinite path. Indeed, the necessity
follows immediately from Corollary \ref{der-arc}(1) and Theorem
\ref{RARcat}(1). Suppose that $Q$ has no right infinite path. Then
$\rep^+\hspace{-1pt}(Q)=\rep^b(Q)$. By Theorem \ref{RARcat}(1),
$\rep^+\hspace{-1pt}(Q)$ is left Auslander-Reiten. If $I$ is an
injective object in $\rep^+\hspace{-1pt}(Q)$, by Proposition
\ref{cor1.14shiping}(2), $I=I_x$ for some $x\in Q_0$. Then the
socle of $I$ is $S_x$ which has a projective cover $P_x$ in
$\rep^+\hspace{-1pt}(Q)$. By Corollary \ref{der-arc}(1),
$D^b(\rep^+\hspace{-1pt}(Q))$ has left almost split triangles.

Next, we claim that $D^b(\rep^+\hspace{-1pt}(Q))$ has right almost
split triangles if and only if $Q$ has no left infinite path.
Indeed, by Corollary \ref{der-arc}(2),
$D^b(\rep^+\hspace{-1pt}(Q))$ has right almost split triangles if
and only if $\rep^+\hspace{-1pt}(Q)$ is right Auslander-Reiten and
every simple representation $S_x$ with $x\in Q_0$ has an injective
hull in $\rep^+\hspace{-1pt}(Q)$, where the second condition is
equivalent to $Q=Q^+$ by Proposition \ref{cor1.14shiping}(3), that
is, $Q$ has no left infinite path. Now the claim follows from
Theorem \ref{RARcat}(2).

Finally, it follows from the above two statements that
$D^b(\rep^+\hspace{-1pt}(Q))$ has almost split triangles if and
only if $Q$ has no infinite path. The proof of the theorem is
completed.

\bigskip\bigskip

{\sc Acknowledgement.}  The first author is supported in part by Fondo Sectorial SEP-CONACYT via the project 43374F
while the last two authors are grateful for support in part from the Natural
Science and Engineering Research Council of Canada.

\end{document}